# Singular and maximal radon transforms: Analysis and geometry

By Michael Christ, Alexander Nagel, Elias M. Stein, and Stephen Wainger*

**Table of Contents**

*Research supported by NSF grants.

## Part 3. Analytic theory

## Part 4. Appendix

## 0. Introduction

The purpose of this paper is to prove the $L^p$ boundedness of singular Radon transforms and their maximal analogues. These operators differ from the traditional singular integrals and maximal functions in that their definition at any point $x$ in $\mathbb{R}^n$ involves integration over a $k$-dimensional submanifold of $\mathbb{R}^n$, depending on $x$, with $k < n$. The role of the underlying geometric data which determines the submanifolds and how they depend on $x$, makes the analysis of these operators quite different from their standard analogues. In fact, much of our work is involved in the examination of the resulting geometric situation, and the elucidation of an attached notion of curvature (a kind of "finite-type" condition) which is crucial for our analysis.

We begin by describing our results, first somewhat imprecisely in order to simplify the statements. We assume that for each $x \in \mathbb{R}^n$ there is a smooth $k$-dimensional submanifold[1] $M_x$, with $x \in M_x$, so that $M_x$ varies smoothly with $x$. Also, for each $x$ we denote by $d\sigma_x$ an integration measure on $M_x$ with smooth density; and $K_x = K_x(y)$ a $k$-dimensional Calderón-Zygmund kernel defined, for $y \in M_x$, which has its singularity at $y = x$. We also assume that the mappings $x \to d\sigma_x$, and $x \mapsto K_x$, are smooth. Then we form the singular Radon transform

$$(0.1) \qquad T(f)(x) = \int_{M_x} f(y) K_x(y) \, d\sigma_x(y),$$

---

[1] In this introduction the sets $M_x$ are assumed to be manifolds for the sake of simplicity, but our main results are formulated in somewhat greater generality.

and the corresponding maximal operator,

$$(0.2) \qquad M(f)(x) = \sup_{r>0} \frac{1}{r^k} \left| \int_{M_x \cap B(x,r)} f(y) \, d\sigma_x(y) \right|$$

where $B(x,r)$ is the ball of radius $r$ centered at $x$.

Among our main results are Theorems 11.1 and 11.2, which state that if $\{M_x\}$ satisfies the curvature condition below, then the operators $T$ and $M$ are bounded on $L^p$, for $1 < p < \infty$.[2]

*The curvature condition.* One way to make the above statements precise is in terms of a parametric representation of the submanifolds $\{M_x\}$. For this purpose we assume given a $C^\infty$ function $\gamma$, defined in a neighborhood of the origin in $\mathbb{R}^n \times \mathbb{R}^k$, taking values in $\mathbb{R}^n$, with $\gamma(x,0) \equiv x$. Then we set $M_x = \{\gamma(x,t), \, t \in U\}$ where $U$ is a suitable neighborhood of the origin in $\mathbb{R}^k$. It is also useful to think of $\gamma$ as a family of (local) diffeomorphisms of $\mathbb{R}^n$, $\{\gamma_t\}$, parametrized by $t$, and given by $\gamma_t(x) = \gamma(x,t)$. Next, starting with a standard Calderón-Zygmund kernel $K$ on $\mathbb{R}^k$, a suitable $C^\infty$ cut-off function $\psi$, and a small positive constant $a$, we redefine (0.1) in a precise form as the principal-value integral

$$(0.1)' \qquad T(f)(x) = \psi(x) \int_{|t| \leq a} f(\gamma(x,t)) \, K(t) \, dt.$$

$M(f)$ can be handled similarly.

The curvature condition needed can now be stated in a number of equivalent ways:

(i) A first form is in terms of a noncommutative version of Taylor's formula. This formula is valid for all families $\{\gamma_t\}$ of diffeomorphisms as above, and is interesting in its own right: it states that there exist (unique) vector fields $\{X_\alpha\}$, with $\alpha = (\alpha_1, \ldots \alpha_k) \neq 0$, so that asymptotically $\gamma_t(x) \sim \exp\left(\sum_\alpha \frac{t^\alpha}{\alpha!} X_\alpha\right)(x)$, as $t \to 0$. The curvature condition $(\mathcal{C}_{\mathfrak{g}})$ is then that the Lie algebra generated by the $X_\alpha$ should span the tangent space to $\mathbb{R}^n$. In the special Euclidean-translation-invariant situation, when $\gamma_t(x) = x - \gamma(t)$, the condition is that the vectors $\left(\frac{\partial}{\partial t}\right)^\alpha \gamma(t)\Big|_{t=0}$ span $\mathbb{R}^n$. The general condition, unlike that special case, is diffeomorphism invariant. Moreover it is highly suggestive, bringing to mind Hörmander's condition guaranteeing sub-ellipticity. However our proofs require, in addition, the equivalent formulation below.

---

[2] It is known that without some curvature conditions, these conclusions may fail utterly.

(ii) An equivalent condition, $(\mathcal{C}_J)$, is stated in terms of repeated compositions of the mapping $\gamma_t$. One defines $\Gamma(\tau) = \Gamma(\tau, x)$ by

$$\Gamma(\tau, x) = \gamma_{t^1} \circ \gamma_{t^2} \cdots \circ \gamma_{t^n}(x)$$

with $\tau = (t^1, t^2, \ldots t^n) \in \mathbb{R}^{nk}$. The condition states that for some $n \times n$ sub-determinant $J$ of $\partial \Gamma/\partial \tau$, and some multi-index $\alpha$, we have $\left(\frac{\partial}{\partial \tau}\right)^\alpha J(\tau) \neq 0$ at $\tau = 0$.

The way the mapping $\tau \to \Gamma(\tau)$ arises can be understood as follows. We decompose the operator $T$ in $(0.1)'$ as $T = \sum_{j=0}^{\infty} T_j$, in the standard dyadic way writing

$$T_j(f)(x) = \int f(\gamma_t(x)) \, \psi_j(x, t) \, dt,$$

with $\psi_j$ supported where $|t| \approx 2^{-j}$. The key point (for the $L^2$ theory) is the almost orthogonality estimate

(0.3) $$\|T_i T_j^*\| + \|T_i^* T_j\| \leq C 2^{-\epsilon|i-j|},$$

and in fact, the more elaborate form from which it is deduced

(0.4) $$\|(T_j T_j^*)^N T_i\| \leq C 2^{-\epsilon'|i-j|}, \text{ when } i \geq j.$$

Now such products as $(T_j)^n$ can be written as

$$(T_j)^n(f)(x) = \int f(\Gamma(\tau, x)) \, \Phi(\tau, x) \, d\tau$$

and $(T_j T_j^*)^N$ can be similarly expressed.

(iii) A third equivalent condition, $(\mathcal{C}_M)$, has a very simple statement in the case $\gamma_t$ is real-analytic: it is that there is no submanifold (of positive co-dimension) which is (locally) invariant under the $\gamma_t$. In the $C^\infty$ case the requirement becomes the noninvariance up to infinite order, as in Definition 8.3 below.

There are a number of other equivalent ways of stating the basic curvature condition, but unlike the geometric formulations (i)–(iii) above, the one that follows is analytic in nature. We consider a variant of the operator (0.1) (or $(0.1)'$), where the singular kernel $K$ is replaced by a $C^\infty$ density, with small support in $t$. Then the condition is that this operator (which in appearance is now more like a standard Radon transform) is smoothing of some positive degree, either in the Sobolev-space sense, or as a mapping from $L^p$ to $L^q$, with $q > p$.

*Background.* The first example of the operator (0.1) arose when the method of rotations was applied to the singular-integrals associated to the heat equation. The operator obtained was the Hilbert transform on the parabola, where $n = 2$, and $\gamma(x,t) = x - (t, t^2)$ (Fabes [10]). The $L^2$ theorem contained there was generalized in Stein and Wainger [44]; see also Alpár [1], Halasz [17] and Kaufman [21]. The method of rotations also suggested the study of the maximal operator (0.2) in connection with Poisson integrals on symmetric spaces (Stein [40]). The initial $L^p$ results were then obtained by Nagel, Riviere, and Wainger [26], [27] and the general Euclidean-translation-invariant theory, when $k = 1$, was worked out in Stein and Wainger [45]; see also Stein [42]. All these results relied in a crucial way on the use of the Plancherel formula; here curvature entered via the method of stationary phase and the estimation of certain Fourier transforms.

With this, attention turned to the problem of "variable" manifolds, i.e. the non-translation invariant case, where new tools were needed. In Nagel, Stein, Wainger [28] such an $L^2$ result is obtained in the special case of certain curves in the plane; it was an early indication that orthogonality (e.g. the consideration of $TT^*$ instead of $T$) may be decisive. The efforts then focused on the setting of nilpotent Lie groups, with the results of Geller and Stein [11] for the Heisenberg group, various extensions by Müller [22], [23], [24] and culminated with Christ [5], where the general group-invariant case for $k = 1$ was established; this last was generalized to higher $k$ in Ricci and Stein [34].

Meanwhile, prompted by the connections with the $\bar{\partial}$-Neumann problem for pseudo-convex domains, Phong and Stein [31] worked out the theory under the assumption of nowhere vanishing *rotational curvature* (which is the "best possible" situation, and which also arises naturally in the theory of the Fourier integral operator); see also Greenleaf and Uhlmann [14].

*Methods used.* Here we want to highlight three techniques which are very useful in our work. First, in order to exploit the curvature condition as expressed in (i), we lift matters to a higher-dimensional setting, where the corresponding vector fields are "free." One of the consequences of this lifting is that we now have local dilations, which essentially allow us to re-scale crucial estimates to unit scale (e.g. in effect to reduce (0.4) to the case $j = 0$). In applying this lifting, we have imitated a general approach used in Rothschild and Stein [36]. It should be noted, however, that the lifting technique is not used in establishing the equivalence of the various curvature conditions. Nor is it used in our proof of the smoothing property of nonsingular generalized Radon transforms with $K \in C^\infty$. A second key idea is the fact that the Cotlar-lemma estimates (0.3) can be reduced to inequalities like (0.4). This method already occurred in Christ [5]. Thirdly, in proving the inequalities (0.4), we need to know that the integral kernel of the operator $(T_0 T_0^*)^N$ has some smoothness, and for this we utilize the curvature condition in the form (ii). This is accomplished by a theorem guaranteeing that certain transported measures have relatively smooth Radon-Nikodym derivatives, generalizing earlier forms of this principle in Christ [5], and Ricci and Stein [34].

*Organization of the paper.* We have divided this work into three parts. Part 1 contains the background material. In it we recall a number of facts needed, and we also state reformulations of some theorems in the literature. Therefore detailed proofs are not given here. Part 2 is devoted to the definitions and study of the various curvature conditions, and to the proof of their equivalence. Numerous illustrative examples are also provided. Part 3 contains the proofs of the $L^p$ estimates.

A word of explanation about the writing of this paper may be in order. The main results grew out of work begun over a dozen years ago. At that time the four of us joined forces, basing our work in part on the manuscript Christ [6], and ideas developed by the three other authors. A draft containing all the main results of the present work was prepared about a year later, and over the next few years the results were described in several lectures given by the authors; in addition they were the subject of graduate courses at Princeton in 1991 and 1996 given by one of us (E.M.S.). However, a final version was not prepared until recently, one reason being that our efforts were viewed as part of a larger project which we had hoped to complete. With that project still not done, it seemed best not to delay publication any longer.

We are grateful to the referee for useful corrections.

## Part 1. Background and preliminaries

This first part is devoted to recalling or elaborating several known ideas needed for the proofs of the theorems in parts 2 and 3 below. The proofs of the assertions made here will for the most part be omitted, because the statements are either well-known, or can be established by minor modification of the existing proofs in the literature.

### 1. Vector fields and the exponential mapping

We begin by recalling some basic facts concerning vector fields and the exponential mapping. A vector field is given in local coordinates by

$$X = \sum_{j=1}^{n} a_j(x) \frac{\partial}{\partial x_j} = a(x) \cdot \nabla_x$$

with $a(x) = (a_1(x), \ldots, a_n(x))$, where the $a_j$ are real-valued $C^\infty$ functions defined on some open subset $O$ of $\mathbb{R}^n$.

Associated to the vector field $X$ is the flow $\varphi_t = \varphi_t(x)$. It satisfies the differential equation
$$\frac{d}{dt}\varphi_t = a(\varphi_t)$$
and the initial condition
$$\varphi_0(x) = x, \quad \text{for every } x \in O.$$
The existence and uniqueness theorems for ordinary differential equations guarantee that the mappings $x \mapsto \varphi_t(x)$, defined for sufficiently small $t$, satisfy

(1.1) $$(\varphi_t \circ \varphi_s)(x) = \varphi_{t+s}(x)$$

if $x \in O_1$, with $\bar{O}_1 \subset O$, when $t$ and $s$ are sufficiently small. (For later purposes we should note that the above theorems guarantee the existence of $\varphi_t(x)$, up to $t = 1$, if $x \in O_1$ with $\bar{O}_1$ a compact subset of $O$, and with the $C^0$ norm of $a$ sufficiently small.) For any $f \in C^\infty$ defined near $O$,
$$\left.\frac{d f(\varphi_t(x))}{dt}\right|_{t=0} = X(f)(x),$$
and more generally,

(1.2) $$\frac{d f(\varphi_t(x))}{dt} = X(f)(\varphi_t(x)), \quad x \in O_1$$

if $t$ is sufficiently small; this is merely a restatement of the defining equation $d\varphi_t/dt = a(\varphi_t)$.

These facts suggest that we write $\varphi_t = \exp(tX) = e^{tX}$, also $\varphi_t(x) = (\exp tX)(x)$. Using (1.2) repeatedly, we obtain a version of Taylor's formula, namely

(1.3) $$f(\exp(tX)(x)) = \sum_{k=0}^{N} \frac{(X^k(f))(x)}{k!} + O(t^{N+1}), \text{ as } t \to 0$$

for every $x \in O_1$, whenever $f \in C^{(N+1)}(O)$.

Next, let $X_1, \ldots, X_p$ be a finite collection of vector fields defined in $O$. For $u = (u_1, \ldots u_p) \in \mathbb{R}^p$ sufficiently small, and keeping in mind the parenthetical remark made earlier, we can define $\exp(u_1 X_1 + \ldots u_p X_p)$ to be $\exp(tX)$, with $t = 1$, where $X = u_1 X_1 + u_2 X_2 \ldots + u_p X_p$. As a result the mapping
$$(x, u) \mapsto \exp(u_1 X_1 \ldots + u_p X_p)(x)$$
is smooth jointly in $x$ and $u$, as long as $x \in O_1$ and $u$ is sufficiently small. There are two consequences one should note. First, the generalization of (1.3), namely

(1.4) $$f(\exp(u_1 X_1 \ldots + u_p X_p)x) = \sum_{k=0}^{N} \left(\sum_{j=1}^{p} u_j X_j\right)^k (f)(x)/k! + O(|u|^{N+1})$$

as $u \to 0$, whenever $f \in C^{(N+1)}(O)$.

Second, suppose $X_1, \ldots X_n$ are $n$ vector fields, linearly independent at each point $x \in O$. Then for each $y \in O_1$, the mapping $(u_1, u_2, \ldots u_n) \mapsto \exp(u_1 X_1 + \ldots + u_n X_n)(y)$ is a local diffeomorphism of a small neighborhood of the origin in $\mathbb{R}^n$ onto a corresponding neighborhood of $y$. As a result, $(u_1, \ldots u_n)$ can be taken to be a coordinate system for the point $x = \exp(u_1 X_1 + \ldots u_n X_n)(y)$, which is centered at $y$. These are the exponential coordinates, determined by $X_1, \ldots, X_n$ and centered at $y$.

Finally, we recall that if $X$ and $Y$ are a pair of vector fields defined in $O$, so is their commutator bracket $[X, Y] = XY - YX$. This bracket makes the vector fields defined on $O$ into a Lie algebra, to which topic we now turn.

## 2. Lie algebras

For our purposes a Lie algebra $\mathcal{L}$ is a (non-associative) algebra over $\mathbb{R}$, whose product, denoted by $[X, Y]$, with $X, Y \in \mathcal{L}$, satisfies $[X, Y] = -[Y, X]$ and

$$[X, [Y, Z]] + [Y, [Z, X]] + [Z, [X, Y]] = 0, \qquad \text{for all } X, Y, Z \in \mathcal{L}.$$

We consider several examples.

*Example* 2.0.　The Lie algebra of vector fields defined over $O \subset \mathbb{R}^n$, as above in Section 1.

*Example* 2.1.　Let $\mathcal{A}$ be an associative algebra over $\mathbb{R}$ with product $X \cdot Y$. If we define the bracket $[X, Y] = X \cdot Y - Y \cdot X$, then $\mathcal{A}$ becomes a Lie algebra, which we denote by $\mathcal{L}(\mathcal{A})$.

*Example* 2.2.　In example 1, we start with $\mathcal{A}$ the *free* associative algebra, freely generated by $p$ generators, $X_1, \ldots, X_p$. The algebra $\mathcal{A}$ can be characterized by its universal properties as follows: First $\mathcal{A}$ is generated by $X_1, \ldots, X_p$, that is, no proper subalgebra of $\mathcal{A}$ contains $\{X_1, \ldots X_p\}$. Second, if $\mathcal{A}'$ is *any* associative algebra, and $\Phi$ is a mapping from the set $\{X_1, \ldots X_p\}$ to $\mathcal{A}'$, then $\Phi$ can be extended (uniquely) to a homomorphism from $\mathcal{A}$ to $\mathcal{A}'$. The algebra $\mathcal{A}$ can be realized as $\bigoplus_{k=1}^{\infty} V^{(k)}$, where $V^{(k)}$ is the tensor product of $k$ copies of $V$, and $V$ is the vector space spanned by $X_1, X_2, \ldots X_p$. (Further details can be found in [20].)

From the associative algebra $\mathcal{A}$ we form the Lie algebra $\mathcal{L}(\mathcal{A})$ as in Example 1, and then pass to $\mathcal{L}_0(\mathcal{A})$ which is the Lie subalgebra of $\mathcal{L}(\mathcal{A})$ generated by $X_1, \ldots, X_p$. This will be called the *free* Lie algebra (with $p$ generators), and will be written $\mathcal{F}_p = \mathcal{L}_0(\mathcal{A})$.

$\mathcal{F}_p$ has the following universal property:

PROPOSITION 2.1. *Suppose that $\mathcal{L}'$ is any Lie algebra, and $\Phi$ is a mapping from the set $\{X_1, \ldots, X_p\}$ into $\mathcal{L}'$. Then $\Phi$ can be (uniquely) extended to a Lie algebra homomorphism of $\mathcal{F}_p$ to $\mathcal{L}'$.*

The proof that $\mathcal{F}_p$ has this universal property can be found in [20]. Once $p$ is fixed, $\mathcal{F}_p$ is uniquely determined (up to Lie algebra isomorphism) by this property. This uniqueness is closely related to the following assertion, which is itself a simple consequence of the proposition.

COROLLARY 2.2. *Let $V$ be the vector subspace of $\mathcal{F}_p$ spanned by $X_1, \ldots X_p$. Then any linear isomorphism $\Phi$ of $V$ can be (uniquely) extended to a Lie algebra automorphism of $\mathcal{F}_p$.*

We now come to (nonisotropic) dilations of $\mathcal{F}_p$. We fix $p$ strictly positive integers, $a_1, a_2, \ldots a_p$; these will be the exponents of the dilations. For each $r > 0$, we consider the mapping $\Phi_r$ defined on $V$ by $\Phi_r(\sum_{j=1}^p c_j X_j) = \sum_{j=1}^p c_j r^{a_j} X_j$. Then by the corollary, $\Phi_r$ extends in a unique way to an automorphism of $\mathcal{F}_p$; this extension will also be denoted by $\Phi_r$. One also notes that $\Phi_{r_1} \circ \Phi_{r_2} = \Phi_{r_1 r_2}$. We call a repeated commutator $[X_{i_1}, [X_{i_2} \ldots [X_{i_{k-1}}, X_{i_k}] \ldots]$ involving generators $X_{i_1}, X_{i_2}, \ldots, X_{i_k}$, a commutator of length $k$; such a commutator undergoes multiplication by the factor $r^{a_{i_1} + a_{i_2} \cdots + a_{i_k}}$ under the action of $\Phi_r$, and is said to be homogeneous of degree $a_{i_1} + a_{i_2} \ldots + a_{i_k}$. Since the linear span of commutators of all lengths is $\mathcal{F}_p$, we can write $\mathcal{F}_p$ as a direct sum

$$(2.1) \qquad \mathcal{F}_p = \bigoplus_{\ell=1}^{\infty} \mathcal{F}_p^{\ell},$$

where $\mathcal{F}_p^{\ell}$ denotes the subspace of all elements that are homogeneous of degree $\ell$ under $\Phi_r$. The decomposition (2.1) makes $\mathcal{F}_p$ into a graded Lie algebra, with $[\mathcal{F}_p^{\ell_1}, \mathcal{F}_p^{\ell_2}] \subset \mathcal{F}_p^{\ell_1 + \ell_2}$. We should note that of course the free Lie algebra $\mathcal{F}_p$ is infinite-dimensional.

*Example* 2.3. Taking into account the gradation above for $\mathcal{F}_p$, for any positive integer $m$ we define $\mathcal{I}^m = \bigoplus_{\ell > m} \mathcal{F}_p^{\ell}$. Since $\mathcal{I}^m$ is spanned by elements homogeneous of degree $> m$, it is clear that $\mathcal{I}^m$ is a (Lie algebra) ideal. For any $m \geq \max(a_1, \ldots a_p)$ we may therefore form the quotient Lie algebra $\mathcal{N} = \mathcal{F}_p / \mathcal{I}_p^m$. We let $Y_1, \ldots Y_p$ denote the images of $X_1, \ldots X_p$ respectively under the natural projection. Then $Y_1, \ldots, Y_p$ are linearly independent elements of $\mathcal{N}$ which generate the Lie algebra $\mathcal{N}$. The automorphism $\Phi_r$ of $\mathcal{F}_p$ induces a corresponding automorphism $\tilde{\Phi}_r$ of $\mathcal{N}$, with $\tilde{\Phi}_r(Y_j) = r^{a_j} Y_j$. Note that $\mathcal{N}$ is naturally identifiable with $\bigoplus_{\ell \leq m} \mathcal{F}_p^{\ell}$. We shall also use the notation $\mathcal{N}_m^{a_1, \ldots a_p}$ for $\mathcal{N}$, to indicate its dependence on the exponents $a_1, a_2, \ldots a_p$, and the order $m$.

*Definition* 2.3. $\mathcal{N}_m^{a_1,\ldots a_p}$ is called the *relatively free nilpotent Lie algebra of order $m$ with $p$ generators*.

We will also refer to $\mathcal{N}$ as being *free up to order $m$*. $\mathcal{N}$ has a universal mapping property: If $\mathfrak{g}$ is an arbitrary Lie algebra containing distinguished elements $X_1,\ldots X_p$ of which any iterated commutator having degree $> m$ vanishes, then there exists a unique Lie algebra homomorphism from $\mathcal{N}$ to $\mathfrak{g}$ mapping $Y_j$ to $X_j$ for each $1 \leq j \leq p$.[3] This degree is defined in the same weighted sense as for $\mathcal{N}_m^{a_1,\ldots a_p}$.

For any *ordered $k$-tuple of integers* $I = \{i_1, i_2, \ldots i_k\}$ with $1 \leq i_j \leq p$, we let $Y_I$ denote as above the commutator $[Y_{i_1} [Y_{i_2} \ldots [Y_{i_{k-1}}, Y_{i_k}]\ldots]]$. Because of the homogeneity induced by the automorphisms $\tilde{\Phi}_r$, we assign to $I$ the degree

$$|I| = a_{i_1} + a_{i_2} \ldots + a_{i_k}.$$

Thus, we can choose a basis of $\mathcal{N}_m^{a_1,\ldots a_p}$ to consist of $\{Y_I\}$, as $I$ ranges over an appropriate subset of the collection of all multi-indices satisfying $|I| \leq m$. Given $m$ and $p$, we choose this collection $\{I\}$ and the resulting basis $\{Y_I\}$ once and for all. We shall refer to the chosen fixed collection $\{I\}$ as *basic*. Note that each element $Y \in \mathcal{N}$ can be written as $Y = \sum_{I \text{ basic}} c_I Y_I$, and the dilations

$$Y = \sum_{I \text{ basic}} c_I Y_I \to \sum_{I \text{ basic}} c_I r^{|I|} Y_I, \quad r > 0,$$

are automorphisms of $\mathcal{N}$. The "relatively free" Lie algebra $\mathcal{N}$ will be one of our chief tools in what follows.

## 3. The Baker-Campbell-Hausdorff formula

Let $\mathcal{A}$ be an associative algebra over $\mathbb{R}$. We define a formal power series to be an expression in the indeterminate $t$ of the form

$$A(t) = \sum_{k=0}^{\infty} a_k t^k$$

where the coefficients $a_k$ are elements of $\mathcal{A}$. No restriction is placed on the sequence $\{a_k\}_{k=0}^{\infty}$; therefore giving the formal power series $A(t)$ is the same as prescribing the (arbitrary) sequence $\{a_k\}$.

These series can be added and multiplied in the standard way. Also if $a \in \mathcal{A}$, we can define $\exp(ta)$ to be the formal power series given by $\sum_{k=0}^{\infty} a^k t^k/k!$. Note that if $A(t)$ is a formal power series without constant term, then $\exp(A(t)) = \sum_{k=0}^{\infty} (A(t))^k/k!$ is itself a formal power series $\sum_{k=0}^{\infty} b_k t^k$, where each $b_k$ is a (noncommutative) polynomial in $a_1, \ldots a_k$.

---

[3] An alternative method for constructing $\mathcal{N}$ is to show first that given any $p, m$, there exists $M < \infty$ such that for any nilpotent Lie algebra $\mathfrak{g}$ generated by $p$ of its elements $X_1,\ldots X_p$, such that any iterated commutator having degree $> m$ vanishes, the dimension of $\mathfrak{g}$ is at most $M$. A second step is to show that any such $\mathfrak{g}$ having maximal dimension has the required universal mapping property.

Given any two elements $a, b \in \mathcal{A}$, we define a Lie polynomial in $a$ and $b$ to be a finite linear combination of repeated commutators of the form $[a_{i_1}, [a_{i_2} \ldots [a_{i_{k-1}}, a_{i_k}] \ldots ]]$, where each $a_i$ is either $a$ or $b$. The Baker-Campbell-Hausdorff formula then states:

PROPOSITION 3.1. *Suppose $a$ and $b$ belong to $\mathcal{A}$. Then there is a formal power series*

$$C(t) = \sum_{k=1}^{\infty} c_k t^k$$

*so that*

(3.1) $$\exp(ta) \cdot \exp(tb) = \exp(C(t)).$$

*Here $c_k = c_k(a, b)$ is a homogeneous Lie polynomial of degree $k$ in $a$ and $b$.*

A proof of this theorem may be found in [20]. The assertion that $c_k$ is homogeneous of degree $k$ means that $c_k(ra, rb) = r^k c_k(a, b)$, when $r \in \mathbb{R}$. The first few polynomials in the formula are $c_1(a, b) = a + b$, $c_2(a, b) = \frac{1}{2}[a, b]$, $c_3 = \frac{1}{12}[a, [a, b]] + \frac{1}{12}[b, [a, b]]$.

As a consequence of the above formal identity we can obtain the following analytic version for vector fields. Let $X_1, X_2, \ldots X_p, Y_1, \ldots, Y_p$ be a collection of real smooth vector fields defined on some open set $O \subset \mathbb{R}^n$. For $u = (u_1, \ldots u_p)$ and $v = (v_1, \ldots, v_p)$ in $\mathbb{R}^p$ write $u \cdot X = u_1 X_1 + u_2 X_2 \ldots + u_p X_p$, $v \cdot Y = v_1 Y_1 + \ldots + v_p Y_p$. Note that by what was said in Section 1, we can define the local diffeomorphisms $\exp(u \cdot X)$ and $\exp(v \cdot Y)$ when $u$ and $v$ are sufficiently small.

COROLLARY 3.2. *As local diffeomorphisms, for each $N > 0$,*

(3.2) $$\exp(vY) \cdot \exp(u \cdot X) = \exp\left(\sum_{k=1}^{N} c_k(u \cdot X, v \cdot Y)\right) + O(|u| + |v|)^{N+1}$$

*as $|u| + |v| \to 0$.*

Note that $c_k(u \cdot X, v \cdot Y)$ is itself a vector field whose coefficients depend on $u$ and $v$ as homogeneous polynomials of degree $k$; this is because $c_k(a, b)$ is a Lie polynomial in $a, b$ of degree $k$. Thus, each of the exponentials in (3.2) is well-defined as long as $u$ and $v$ are sufficiently small.

To prove the corollary we write for $a, b \in \mathcal{A}$,

$$S_N(a, b) = \sum_{k=1}^{N} c_k(a, b),$$

and note that

$$(3.3) \quad \left(\sum_{k=0}^{N} \frac{t^k a^k}{k!}\right)\left(\sum_{k=0}^{N} \frac{t^k b^k}{k!}\right) = \sum_{k=0}^{N} \frac{(S_N(ta, tb))^k}{k!} + E(t)$$

where $E(t)$ is a polynomial in $t$, whose coefficients of $t^k$ vanish for all $k \leq N$. Indeed (3.3) is a direct consequence of the formal identity (3.1) once we observe that

$$\exp(ta) = \sum_{k=0}^{N} \frac{t^k a^k}{k!} + E_1(t), \quad \exp(tb) = \sum_{k=0}^{N} \frac{t^k b^k}{k!} + E_2(t)$$

and

$$\exp(C(t)) = \sum_{k=0}^{N} \frac{(S_N(ta, tb))^k}{k!} + E_3(t),$$

where the $E_j(t)$ are formal power series whose coefficients of $t^k$ vanish for all $k \leq N$. Now (3.3) shows that

$$(3.4) \quad \left(\sum_{k=0}^{N} \frac{a^k}{k!}\right)\left(\sum_{k=0}^{N} \frac{b^k}{k!}\right) = \sum_{k=0}^{N} \frac{(S_N(a,b))^k}{k!} + R$$

where $R$ is a (noncommutative) polynomial in $a$ and $b$ whose terms are each homogeneous of degree $> N$. Finally, (3.2) follows from the Taylor formula (1.4), when $a = u_1 X_1 \ldots + u_p X_p$, and $b = v_1 Y_1 + \ldots + v_p Y_p$.

We shall also need a more extended version of Corollary 3.2 which can be deduced in the same way from (3.4). We use the notation $u^\alpha = u_1^{\alpha_1} u_2^{\alpha_2} \ldots u_p^{\alpha_p}$, if $u \in \mathbb{R}^p$. Let

$$P(u, X) = \sum_{0 < |\alpha| \leq m} u^\alpha X_\alpha, \quad Q(v, Y) = \sum_{0 < |\alpha| \leq m} v^\alpha Y_\alpha$$

be polynomials in $u, v$, without constant term, whose coefficients are vector fields.

COROLLARY 3.3. *For each $N > 0$,*

$$(3.5) \quad \exp(Q(v, Y)) \cdot \exp(P(u, X))$$
$$= \exp\left(\sum_{k=1}^{N} c_k(P(u, X), Q(v, Y))\right) + O((|u| + |v|)^{N+1}),$$

*as $|u| + |v| \to 0$.*

## 4. The Lie group corresponding to $\mathcal{N}$

By the use of the Baker-Campbell-Hausdorff formula we can describe a Lie group corresponding to the Lie algebra $\mathcal{N}_m^{a_1, \ldots a_p}$. As is well-known from the general theory of Lie groups, there is a unique connected, simply connected

Lie group $N = N_m^{a_1,\ldots a_p}$ whose Lie algebra is $\mathcal{N}_m^{a_1,\ldots a_p}$; since the latter is nilpotent, the corresponding exponential mapping is a diffeomorphism, and the underlying space of $N$ may be identified with $\mathbb{R}^d$, where $d =$ dimension $(\mathcal{N}_m^{a_1,\ldots a_p})$. For these facts see [18],[33],[46].

The above assertion may be described more explicitly as follows. In the identification of $N$ with $\mathbb{R}^d$, the group identity is the origin in $\mathbb{R}^d$, and the Lie algebra $\mathcal{N}$ consists of the left-invariant vector fields on $N$ (i.e. on $\mathbb{R}^d$). The exponential map leads to the identification of each $u \in N$ with $\exp\left(\sum_{I \text{ basic}} u_I Y_I\right)(0)$, which we write more simply as $\exp\left(\sum_{I \text{ basic}} u_I Y_I\right)$. Here $u = (u_I)_{I \text{ basic}}$ are coordinates for $\mathbb{R}^d$.

The multiplication law in $N$ is a consequence of the formula (3.1). It takes the form

$$(4.1) \quad \exp\left(\sum_{I \text{ basic}} v_I Y_I\right) \cdot \exp\left(\sum_{I \text{ basic}} u_I Y_I\right) = \exp\left(\sum_{I \text{ basic}} P_I(u,v) Y_I\right).$$

where $P_I(u,v)$ is a polynomial in $u$ and $v$ which is homogeneous of degree $|I|$ in the following sense. Recall the dilations defined on $\mathcal{N}_m^{a_1,\ldots a_p}$. They induce corresponding dilations $\delta_r : \mathbb{R}^d \mapsto \mathbb{R}^d$:

*Definition* 4.1. For any $x = \exp(\sum_{I \text{ basic}} u_I Y_I)$ and $r > 0$,

$$\delta_r(x) = \exp\left(\sum_{I \text{ basic}} r^{|I|} u_I Y_I\right).$$

Then $P_I$ is *homogeneous* in the sense that

$$P_I(\delta_r(u), \delta_r(v)) = r^{|I|} P_I(u,v).$$

*Definition* 4.2. The *norm function* $\rho$ and *quasi-distance* $d$ on $N$ are

$$\rho(u) = \sum_I |u_I|^{1/|I|}, \quad \text{and} \quad d(x,y) = \rho(x^{-1}y).$$

These are linked with the dilation structure through the identity

$$\rho(\delta_r u) = r\rho(u) \quad \text{for all } u \in N, \ r \in \mathbb{R}^+.$$

## 5. Free vector fields

Next we treat the notion of a collection of real vector fields, $\{X_1, \ldots X_p\}$ being "free", relative to exponents $a_1, a_2, \ldots a_p$ and the order $m$. We assume that $X_1, \ldots X_p$ are real, smooth vector fields defined (on some open set $O$) in $\mathbb{R}^d$. Here $d = \dim(\mathcal{N}_m^{a_1,\ldots a_p})$. For any $k$-tuple $I = \{i_1, i_2, \ldots i_k\}$ with $1 \leq i_j \leq p$, we write as before $X_I$ for the corresponding $k$-fold commutator, and $|I| = a_{i_1} + a_{i_2} \ldots + a_{i_k}$.

*Definition* 5.1. A *collection of vector fields* $\{X_1, \ldots X_p\}$ defined in an open subset $O$ of $\mathbb{R}^d$ is said to be *free relative to the exponents* $a_1, a_2, \ldots a_p$, and the order $m$, if $d$ equals the dimension of $\mathcal{N}_m^{a_1, \ldots a_p}$ and

(5.1)    $\{X_I(x)\}_{|I| \leq m}$ spans the tangent space of $\mathbb{R}^d$, for each $x \in O$.

Note that if (5.1) holds, then the collection $\{X_I\}_{I \text{ basic}}$ already spans and thus forms a basis, because any linear relation among the $Y_I$ in $\mathcal{N}$ with $|I| \leq m$, implies the corresponding linear relation among the $X_I$.

Now a fundamental idea we shall use (as in [36],[12],[19]) is that a collection of free vector fields (relative to $a_1, \ldots a_p$ and $m$) can in many ways be well approximated by the Lie algebra $\mathcal{N}_m^{a_1, \ldots a_p}$ and its action as a collection of left-invariant vector fields on the group $N$. In particular, on $N$ there are the following objects of importance: the multiplication law of the group (in the form of the mapping $(x,y) \to x^{-1} \cdot y$ of $N \times N \to N$); the automorphic dilations on $N$, coming from the dilations on $\mathcal{N}$; and the corresponding norm function and quasi-distance on $N$, as defined in Section 4.

The analogue of the first of these in our general setting will be the mapping $\Theta$ defined as follows. For each $x \in \mathbb{R}^d$ (more precisely, for $x \in O$), we consider the mapping $y \to \Theta_x(y)$ from a neighborhood of $x$ to a neighborhood of the origin, given by:

*Definition* 5.2. $\Theta_x(y) = (u_I)$ where $y = \exp(\sum_{I \text{ basic}} u_I X_I)(x)$.

By the properties of the exponential mapping described in Section 1, the mapping $y \to \Theta_x(y)$ is a diffeomorphism of a neighborhood of $x \in O$ with a neighborhood of the origin. Note that since $\exp\left(-\sum_{I \text{ basic}} u_I Y_I\right)(y) = x$, we have

(5.2)    $\Theta_x(y) = -\Theta_y(x).$

Consider now the special case when $\{X_i\}_{i=1}^p$ equals the collection $\{Y_i\}_{i=1}^p$ of left-invariant vector fields on the group $N$ discussed in Section 4. Then by left invariance,

$$\exp\left(\sum_{I \text{ basic}} u_i Y_I\right)(x) = x \cdot \exp\left(\sum_{I \text{ basic}} u_I Y_I\right).$$

So, via the identification of $N$ with $\mathcal{N}$, we see that $\Theta_x(y) = x^{-1} \cdot y$ in this case.

Recall the dilations, norm function, and left-invariant quasi-distance on $N$ defined in Section 4. In our more general context, this leads us to the following two definitions. First, the quasi-distance

$$(5.3) \quad d(x,y) = \rho(\Theta_x(y)) = \sum_{I \text{ basic}} |u_I|^{1/|I|}, \quad \text{if} \quad y = \exp\left(\sum u_I X_I\right)(x).$$

We claim that whenever $y$ and $z$ are in a sufficiently small neighborhood of $x$,

$$(5.4) \quad \begin{cases} \text{(i)} & d(x,y) \geq 0, \text{ and } d(x,y) = 0 \text{ only when } x = y \\ \text{(ii)} & d(x,y) = d(y,x) \\ \text{(iii)} & d(x,z) \leq c(d(x,y) + d(y,z)). \end{cases}$$

Now $d(x,y) = 0$ only if $x = y$, since the exponential mapping is a local diffeomorphism; so property (i) is clear. Property (ii) follows directly from the anti-symmetry property (5.2) of the $\Theta$ mapping. Now turning to property (iii), we write

$$y = \exp \sum_{I \text{ basic}} (u_I X_I)(x), \quad \text{and} \quad z = \exp \sum_{I \text{ basic}} (v_I X_I)(y).$$

Then

$$d(x,y) = \sum_{I \text{ basic}} |u_I|^{1/I}, \quad d(y,z) = \sum_{I \text{ basic}} |v_I|^{1/|I|}.$$

Express

$$(5.5) \quad z = \exp\left(\sum_{I \text{ basic}} w_I X_I\right)(x),$$

and note the alternative representation

$$(5.6) \quad z = \exp\left(\sum_{I \text{ basic}} v_I X_I\right) \cdot \exp\left(\sum_{I \text{ basic}} u_I X_I\right)(x).$$

We shall apply the Baker-Campbell-Hausdorff formula to compare (5.5) with (5.6). In order to do this we make two remarks regarding commutators involving $u_{I_1} X_{I_1}, u_{I_2} X_{I_2}, \ldots$ etc. First, when $|I_1| + |I_2| \leq m$,

$$(5.7) \quad [X_{I_1}, X_{I_2}] = \sum_{I \text{ basic}} c_I X_I$$

where the $c_I$ are the same constants determined by the identical relation that holds in the Lie algebra $\mathcal{N}_m^{a_1, \ldots a_p}$, namely $[Y_{I_1}, Y_{I_2}] = \sum_{I \text{ basic}} c_I Y_I$; this is because $\mathcal{N}$ is free up to order $m$. When $|I_1| + |I_2| > m$, (5.7) no longer holds (in $\mathcal{N}$ the corresponding right-hand side is zero). In this case, we note that each coefficient of the vector field $[u_{I_1} X_{I_1}, v_{I_2} X_{I_2}]$ is $O(|u_{I_1}||v_{I_2}|)$, which is $O(\rho(u)^{|I_1|} \rho(v)^{|I_2|}) = O(\rho(u)^m + \rho(v)^m)$.

Applying the corollary in Section 3 to (5.6), and taking into account the product formula (4.1), we see that

$$(5.8)$$
$$z = \exp\left(\sum_{I \text{ basic}} P_I(u,v) X_I\right)(x) + O\left(\rho(u)^m + \rho(v)^m\right) + O\left(|u|^{m+1} + |v|^{m+1}\right)$$

where $|u| = \sum_{I \text{ basic}} |u_I|$. Since $|u| \leq c\rho(u)$ and $|v| \leq c\rho(v)$, the second $O$ term can be subsumed in the first when $u$ and $v$ are small. Given the diffeomorphic character of the exponential mapping occurring in (5.8), it follows that

$$z = \exp\left(\sum_{I \text{ basic}} w_I X_I\right)(x), \quad \text{where} \quad w_I = P_I(u,v) + O(\rho(u)^m + \rho(v)^m).$$

Now by homogeneity, $P_I(u,v) = O(\rho(u)^{|I|} + \rho(v)^{|I|})$. Thus, clearly $\sum_{I \text{ basic}} |w_I|^{1/|I|} \leq c(\rho(u) + \rho(v))$ and since $d(x,z) = \sum_{I \text{ basic}} |w_I|^{1/|I|}$, the triangle inequality (5.4)(iii) is proved.[4]

Finally, we come to the appropriate notion of dilations in our context. These are (local) dilations $\delta_r^x$, centered at $x$:

*Definition* 5.3. For $y = \exp(\sum_{I \text{ basic}} u_I X_I)(x)$, sufficiently close to $x$, and $r$ sufficiently small,

$$(5.9) \qquad \delta_r^x(y) = \exp\left(\sum_{I \text{ basic}} u_I r^{|I|} X_I\right)(x).$$

An equivalent expression for these local dilations is

$$\delta_r^x(y) = \Theta_x^{-1} \delta_r \Theta_x (y).$$

Here "sufficiently small" parameters $r$ might be quite large; what is required is merely that $|u_I| r^{|I|} = O(1)$, so that the right-hand side of (5.9) will be defined.

Note that by (5.3),

$$d(x, \delta_r^x(y)) = r\, d(x,y).$$

If we let $B(x,r)$ denote the ball $= \{y : d(x,y) < r\}$ then clearly $\delta_r^x(B(x,s)) = B(x,rs)$. Also, if $|B(x,r)|$ denotes the measure of $B(x,r)$,

$$(5.10) \qquad |B(x,r)| \approx r^Q, \quad \text{as} \quad r \to 0$$

where $Q = \sum_{I \text{ basic}} |I|$ is the homogeneous dimension of $\mathcal{N}_m^{a_1, \ldots a_p}$. The symbol $\approx$ means that the ratio $|B(x,r)|/r^Q$ tends to a positive constant as $r \to 0$, uniformly for $x$ in any compact subset of $O$.

Indeed, the mapping $y \to \Theta_x(y) = (u_I)$ is a local diffeomorphism of $y$ near $x$, to points near the origin in the $u$-space; and $y \in B(x,r)$ exactly when $\sum_{I \text{ basic}} |u_I|^{1/|I|} < r$. Of course

$$\left|\left\{u : \sum_{I \text{ basic}} |u_I|^{1/|I|} < r\right\}\right| = cr^Q,$$

---

[4]For the study of (5.4) in a more general setting, see also [29].

for an appropriate constant $c > 0$, as a simple homogeneity argument shows. This establishes (5.10), and with it the doubling property

(5.11) $$|B(x, 2r)| \leq c|B(x,r)|$$

for all sufficiently small $r > 0$.

A last point of significance is a fact which also implies the triangle inequality (5.4), and which will be quite useful in Part III below. It is the assertion that

(5.12) $$\rho(\Theta_x(y_2) - \Theta_x(y_1)) \leq C\{d(y_1, y_2) + d(y_1, y_2)^{1/m} d(x, y_1)^{1-1/m}\}.$$

(See also the analogous statement in [36, §12].)

To prove (5.12) we may assume that $y_1$ is close to $y_2$, and also close to $x$. Then we can write

$$y_1 = \exp\left(\sum_I u_I X_I\right)(x), \quad y_2 = \exp\left(\sum_I v_I X_I\right)(y_1)$$

and alternatively $y_2 = \exp\left(\sum_I w_I X_I\right)(x)$. Here, and below, the sums $\sum_I$ are taken over the basic $I$'s. So we have

$$u = (u_I) = \Theta_x(y_1), \quad w = \Theta_x(y_2),$$

and

$$d(y_1, y_2) = \rho(v), \quad d(x, y_1) = \rho(u).$$

Also

$$\exp\left(\sum_I v_I X_I\right) \cdot \exp\left(\sum_I u_I X_I\right)(x) = \exp\left(\sum_I w_I X_I\right)(x).$$

We apply to this the Baker-Campbell-Hausdorff formula in the same way as in the argument leading to (5.8) and obtain

(5.13) $$w_I = u_I + v_I + Q_I(u,v) + R_I(u,v).$$

Here $u_I + v_I + Q_I(u,v) = P_I(u,v)$ is the term arising in the multiplication formula (4.1) for the group $N$, so $P_I$ is homogeneous of degree $|I|$; the error term $R_I$ is $O(\rho(u)^{m+1} + \rho(v)^{m+1})$. Observe next that when $v = 0$, we have $y_2 = y_1$, which means that $Q_I(u, 0) \equiv 0$ and $R_I(u, 0) \equiv 0$. Thus, writing $Q_I$ as a sum of homogeneous monomials, we see that

$$|Q_I(u,v)| \leq C \sum_{k+\ell=|I|,\ \ell \geq 1} \rho(u)^k \rho(v)^\ell \leq C'\left(\rho(v)\rho(u)^{|I|-1} + \rho(v)^{|I|}\right).$$

Since $R_I(u, v)$ is likewise a smooth function of $u$ and $v$ which vanishes when $v = 0$, each monomial in its Taylor expansion is $O(\rho(u)^k \rho(v)^\ell)$ for some $k, \ell$ satisfying $\ell \geq 1$ and $k + \ell > m$, and hence we get

$$R_I(u,v) = O\left(|v||u|^{m-1} + |v|^m\right) = O\left(\rho(v)\rho(u)^{|I|-1} + \rho(v)^{|I|}\right),$$

because $|I| \leq m$. Hence

$$(\Theta_x(y_2) - \Theta_x(y_1))_I = w_I - u_I = v_I + O(\rho(v)\,\rho(u)^{|I|-1} + \rho(v)^{|I|}).$$

Thus

$$\rho(\Theta_x(y_2) - \Theta_x(y_1)) = \sum_I |w_I - u_I|^{1/|I|} \leq C\{\rho(v) + \rho(v)^{1/m}\,\rho(u)^{1-1/m}\}.$$

Since $\rho(v) = d(y_1, y_2)$ and $\rho(u) = d(x, y_1)$, the assertion (5.12) is proved.

## 6. Freeing vector fields

Let us revert to a neighborhood $O$ of the origin in $\mathbb{R}^n$. Suppose we are given $p$ smooth real vector fields, $X_1, \ldots X_p$, defined on $O$. We write

$$(6.1) \qquad X_i = \sum_{j=1}^n a_i^j(x) \frac{\partial}{\partial x_j}, \quad i = 1, \ldots p.$$

Suppose that we have by some procedure assigned a degree $a_j$ to each $X_j$, where the $a_j$ are strictly positive integers. With these exponents we keep to the notation used previously, and write $X_I$ for the repeated commutator involving $X_{i_1}, X_{i_2}, \ldots X_{i_k}$, where $I = (i_1, i_2, \ldots i_k)$, with $|I| = a_{i_1} + a_{i_2} \ldots + a_{i_k}$.

For the present we make the following assumption on $X_1, \ldots X_p$:

$$(6.2) \qquad \{X_I\}_{|I| \leq m} \text{ spans the tangent space of } \mathbb{R}^n \text{ for each } x \in O.$$

Let $d$ be the dimension of $\mathcal{N}_m^{a_1, \ldots a_p}$. In $\mathbb{R}^d$ we adopt coordinates $u = (x, z) \in \mathbb{R}^n \times \mathbb{R}^{d-n}$.

PROPOSITION 6.1. *Under the assumption (6.2), the $X_i$ can be extended to vector fields $\tilde{X}_i$ defined in a neighborhood of the origin in $\mathbb{R}^d$, taking the form*

$$(6.3) \qquad \tilde{X}_i = \sum_{j=1}^n a_i^j(x) \frac{\partial}{\partial x_j} + \sum_{k=1}^{d-n} b_i^k(x, z) \frac{\partial}{\partial z_k},$$

*so that $\{\tilde{X}_i : i = 1, \ldots, p\}$ is free relative to exponents $a_1, \ldots a_p$ and order $m$, in $\mathbb{R}^d$.*

This is the "lifting" theorem proved in [36] for the case $a_1 = a_2 \ldots = a_p = 1$. (See also the alternate derivations in [12] and [19].) The case where $a_1 = a_2 \ldots = a_{p-1} = 1$ and $a_p = 2$ already occurred in [36]. One can adapt any of those proofs, with minor changes, to establish the more general formulation given here. A different proof will be given in an appendix, Section 22.

# 7. Transporting measures

We shall consider the following situation: $\Phi$ will be a given $C^\infty$ mapping from a closed finite ball $\bar{B}$ in $\mathbb{R}^d$ to $\mathbb{R}^n$, with $d \geq n$. We shall designate the typical point in $\mathbb{R}^d$ by $\tau$, and that in $\mathbb{R}^n$ by $y$. Consider a measure $\psi(\tau)d\tau$ in $\mathbb{R}^d$, whose density $\psi$ is in $C^1(\bar{B})$, and which has compact support in $B$. We let $d\mu = \Phi_*(\psi d\tau)$ be the transported measure in $\mathbb{R}^n$; that is, $\mu$ is defined by the integration formula

$$\int_{\mathbb{R}^n} f(y)\, d\mu(y) = \int_{\bar{B}} f(\Phi(\tau))\, \psi(\tau) d\tau.$$

Our goal is to show that under appropriate conditions $\mu$ is absolutely continuous with respect to Lebesgue measure, and that its Radon-Nikodym derivative $h$, defined by $d\mu(y) = h(y)dy$, possesses the following degree of smoothness.

*Definition* 7.1. For $0 < \delta \leq 1$, $L^1_\delta(\mathbb{R}^n)$ is the Banach space consisting of all functions $h \in L^1(\mathbb{R}^n)$ that satisfy

(7.1) $$\int_{\mathbb{R}^n} |h(y-z) - h(y)|\, dy \leq A|z|^\delta, \quad \text{for all } z \in \mathbb{R}^n.$$

The norm on $L^1_\delta$ is defined to be $\|h\|_{L^1}$ plus the smallest constant $A$ for which (7.1) holds. With these definitions we can state our result.

PROPOSITION 7.2. *Let $J$ be the determinant of some $n \times n$ sub-matrix of the Jacobian matrix $\partial \Phi / \partial \tau$ of $\Phi$. Assume that for some $\alpha$,*

(7.2) $$\partial^\alpha_\tau J(\tau) \neq 0 \quad \text{for every } \tau \in \bar{B}.$$

*Then the transported measure $d\mu = \Phi_*(\psi d\tau)$ is absolutely continuous, and its Radon-Nikodym derivative $h$ belongs to $L^1_\delta$ for all $\delta < (2k)^{-1}$, where $k = |\alpha|$. Moreover, the $L^1_\delta$ norm of $h$ can be controlled in terms of the $C^{k+2}(\bar{B})$ norm of $\Phi$, a lower bound for $\partial^\alpha_\tau J(\tau)$ in $\bar{B}$, the $C^1$ norm of $\psi$, and the numbers $\delta$ and $k$.*

The proposition is an easy variant, in the $C^\infty$ context, of a parallel result formulated for real-analytic mappings in [34]; see also [5]. We will be able to follow closely the proof given in [34, §2], once we have the following lemma.

LEMMA 7.3. *Suppose $F$ is a real-valued function in $C^{(k+1)}(\bar{B})$, and for some $\alpha$, there is $|\partial^\alpha_\tau F(\tau)| \geq b > 0$, throughout $\bar{B}$. Let $k = |\alpha|$. Then*

(7.3) $$\int_{\bar{B}} |F(\tau)|^{-\sigma}\, d\tau \leq A < \infty,$$

*for any $\sigma < 1/k$. The constant $A$ in (7.3) depends only on the norm $\|F\|_{C^{k+1}(\bar{B})}$, the bound $b$, the volume of $B$, and $\sigma$, $k$.*

Versions of this lemma appeared in [23] and [24].

*Proof.* We need to consider only the case $k \geq 1$. We may use Lemma 3.4 in [5], or argue as follows. As is well known, if $0 < \sigma < 1$,

$$|u|^{-\sigma} = c_\sigma \int_{-\infty}^{\infty} e^{i\lambda u} |\lambda|^{-1+\sigma} d\lambda,$$

for an appropriate constant $c_\sigma$. Therefore

(7.4) $\quad \int_{\bar{B}} |F(\tau)|^{-\sigma} d\tau = c_\sigma \int_{-\infty}^{\infty} I(\lambda) |\lambda|^{-1+\sigma} d\tau$, where $I(\lambda) = \int_{\bar{B}} e^{i\lambda F(\tau)} d\tau$.

Now it is known that $|I(\lambda)| \leq C'|\lambda|^{-1/k}$ (where $C'$ depends only on the $C^{k+1}$ norm of $F$ in $\bar{B}$, $b$, and the radius of $B$). For this see [43, Ch. 8], where the result is stated and proved for integrals of the form $\int e^{i\lambda F(\tau)} \psi(\tau) d\tau$ with $\psi$ having compact support in a ball. However, the result holds as well (with nearly identical proof) if $\psi(\tau) d\tau$ is replaced by $\chi_{\bar{B}} d\tau$, where $\chi_{\bar{B}}$ is the characteristic function of $\bar{B}$. Inserting the decay estimate in (7.4) for $|\lambda| \geq 1$, and the trivial estimate $|I(\lambda)| \leq |\bar{B}|$ for $|\lambda| < 1$, gives (7.3) and the lemma. □

By the lemma with $F = J$, it follows that $Z = \{\tau \in \bar{B} : J(\tau) = 0\}$ has Lebesgue measure zero. We shall now cover $\bar{B}/Z$ by a collection of balls $\{B_j\}$ as follows. We let $\underline{c}$ denote a small constant, which we shall fix momentarily. For each $\tau \in B/Z$ we consider the Euclidean ball centered at $\tau$, with radius $\underline{c}|J(\tau)|$, which we write as $B(\tau, \underline{c}|J(\tau)|)$. We choose $\underline{c}$ so small that

(7.5) $\quad 9/10 \leq \left|\dfrac{J(\tau')}{J(\tau)}\right| \leq 11/10$, whenever $\tau' \in B(\tau, 8\underline{c}|J(\tau)|)$.

The existence of such a constant is guaranteed by the fact that $J(\tau)$ is a Lipschitz function, which in turn follows from the fact that $\Phi$ is in $C^{k+2}(\bar{B})$.

Next, let $B(\tau_j, \underline{c}|J(\tau_j)|)$ be a maximal disjoint collection of such balls, and let
$B_j = B(\tau_j, r_j)$, $B_j^* = B(\tau_j, 2r_j)$, where $r_j = 4\underline{c}|J(\tau_j)|$. We claim that

(7.6) $\quad \bigcup_j B_j \supset B/Z.$

In fact, if $\tau \in B/Z$, either it is one of the centers $\tau_j$ (in which case it is covered), or it is not. In the latter case, the ball $B(\tau, \underline{c}|J(\tau)|)$ must intersect one of the balls $B(\tau_j, \underline{c}|J(\tau_j)|)$, for otherwise that collection would not be maximal. Hence either $\tau \in B(\tau_j, 3\underline{c}|J(\tau_j)|)$ or $\tau_j \in B(\tau, 3\underline{c}|J(\tau)|)$. In either case it follows from (7.5) that $|J(\tau)| \leq \frac{11}{9}|J(\tau_j)|$, and hence

$$\tau \in B(\tau, \underline{c}|J(\tau)|) \subset B(\tau_j, 4\underline{c}|J(\tau_j)|) = B_j.$$

Note that some of the balls $B_j$ may go outside $B$, but this does not matter. Next, observe that the $B_j^*$ have the bounded intersection property: there exists $M < \infty$ such that no point is covered by more than $M$ of the balls $B_j^*$. In fact, all balls $B_j$ for which $\tau \in B_j^*$ have comparable radii, are disjoint, and are contained in a ball of comparable radius, because of (7.5).

Once the balls $B_j$ have been determined, we choose in the standard way a corresponding partition of unity $\{\eta_j\}$, $C^\infty$ functions supported in $B_j^*$, with $\sum_j \eta_j = 1$ on $B/Z$, and satisfying $|\nabla \eta_j| \le Cr_j^{-1}$. We then let $d\mu_j = \Phi_*(\psi(\tau)\eta_j(\tau)d\tau)$ and prove (as in [34, §2]) that $d\mu_j = h_j(y)dy$ with

$$\int |h_j(y)|dy \le Cr_j^d, \quad \int |\nabla h_j(y)|dy \le Cr_j^{d-2}.$$

Since $h = \sum_j h_j$, it then follows[5] that $h \in L^1_\delta$ with $\delta = \sigma/2$, if $\int_{\bar{B}} |J(\tau)|^{-\sigma} d\tau < \infty$. Therefore, Lemma 7.3 establishes the proposition.

# Part 2. Geometric theory

## 8. Curvature: Introduction

This section summarizes the main features of the geometric aspects of this paper. Proofs, details and elaborations are presented in Sections 9 and 10. Readers interested primarily in the later estimates might choose to skip over those two sections on a first reading.

The following notation and language will be used. If $F, F_\alpha$ are $C^\infty$ functions of $(x,t)$ and of $x$, respectively, we write

$$F \sim \sum_\alpha t^\alpha F_\alpha$$

to mean that the right-hand side is the Taylor series of $F$ with respect to $t$, at $t = 0$. Taylor expansions with respect to other variables will also be denoted by the symbol $\sim$. $\mathbb{N}$ denotes the set $\{0, 1, 2, \ldots\}$ of nonnegative integers. When we say that a function is defined in $\mathbb{R}^m$ or maps $\mathbb{R}^m$ to $\mathbb{R}^n$, we will mean that it is defined merely in some open subset of $\mathbb{R}^m$, and that its range is contained in $\mathbb{R}^n$.

*Definition* 8.1. The set of all iterated commutators of a collection $\{Y_\beta\}$ of $C^\infty$ vector fields in an open set $U$ is defined to be the smallest set $S$ of vector fields in $U$ such that (i) each $Y_\beta$ belongs to $S$, and (ii) if $V, W \in S$ then $[V, W] \in S$.

The Lie algebra generated by $\{Y_\beta\}$ is defined to be the smallest $C^\infty$ submodule, containing every iterated commutator of $\{Y_\beta\}$, of the set of all $C^\infty$ vector fields on $U$.

---

[5]See [34, pp. 62–63]. (A subscript $j$ is missing in equation (8) on page 62.)

8.1. *Three notions of curvature.* Suppose that $\gamma$ is a $C^\infty$ function defined in some neighborhood of $(x_0, 0)$ in $\mathbb{R}^n \times \mathbb{R}^k$, always satisfying the hypothesis

(8.1) $$\gamma(x, 0) \equiv x.$$

Define
$$\gamma_t(x) \equiv \gamma(x, t).$$

Attention will always be restricted to small $t$.

Such a mapping $\gamma$ may be regarded in either of two ways. First, $x \mapsto \gamma_t(x)$ is a family of diffeomorphisms of $\mathbb{R}^n$, depending smoothly on the parameter $t \in \mathbb{R}^k$. Second, under the auxiliary hypothesis that the differential $\partial \gamma / \partial t$ has rank $k$ when $t = 0$, $t \mapsto \gamma(x)$ parametrizes for each $x$ a $k$-dimensional submanifold $M_x \subset \mathbb{R}^n$ containing $x$, by

(8.2) $$M_x = \{\gamma(x, t) : t \in U\}$$

where $U \subset \mathbb{R}^k$ is a small neighborhood of the origin.

The inverse diffeomorphism, for each $t$, will be denoted by either $\gamma_t^{-1}(x)$ or $\gamma^{-1}(x, t)$. The rank condition is natural and will be satisfied in most of our examples, but our theory does not require it; we permit $k$ to be an arbitrary positive integer, not necessarily less than the ambient dimension $n$.

A fundamental notion is that of an invariant submanifold:

*Definition* 8.2. A *submanifold* $M \subset \mathbb{R}^n$ is *locally invariant under $\gamma$ at $x_0$* if there exists a neighborhood $V$ of $(x_0, 0)$ in $M \times \mathbb{R}^k$ such that $\gamma(x, t) \in M$ for every $(x, t) \in V$.

This notion of invariance would be appropriate and sufficient for our purpose if we were working solely within the class of real analytic mappings. However because we wish to allow mappings which are merely of class $C^\infty$, and because our theory is more naturally formulated in terms of Taylor expansions than of convergent power series, the following weaker notion is more germane.

*Definition* 8.3. A *submanifold* $M$ of $\mathbb{R}^n$ containing $x_0$ is *invariant* under $\gamma$ to infinite order at $x_0$ if for all $(x, t) \in M \times \mathbb{R}^k$ sufficiently close to $(x_0, 0)$,

$$\text{distance}\,(\gamma(x, t), M) = O\Big(\text{distance}\,(x, x_0) + |t|\Big)^N \quad \text{as } x \to x_0 \text{ and } t \to 0$$

for every positive integer $N$.

When $\gamma, M$ are respectively a real analytic mapping and a real analytic submanifold, local invariance at $x$ is equivalent to invariance to infinite order at $x$. But quantities such as $\exp(-1/|t|)$ and $\exp(-1/|x-x_0|)$, which do not affect the Taylor expansion of $\gamma$, may prevent a submanifold from being genuinely invariant.

*Definition* 8.4. $\gamma$ satisfies curvature condition $(\mathcal{C}_M)$ at $x_0$ if there exists no $C^\infty$ submanifold of $\mathbb{R}^n$, of positive codimension, that is invariant under $\gamma$ to infinite order at $x_0$.

A second curvature condition is formulated in terms of certain vector fields, whose existence and uniqueness are guaranteed by the next theorem.

THEOREM 8.5. *Let $\gamma$ be any $C^\infty$ mapping from a neighborhood of $(x_0, 0)$ $\in \mathbb{R}^n \times \mathbb{R}^k$ to $\mathbb{R}^n$, satisfying $\gamma(x, 0) \equiv x$. Then there exists a unique collection $\{X_\alpha : 0 \neq \alpha \in \mathbb{N}^n\}$ of $C^\infty$ vector fields, all defined in some common neighborhood $U$ of $x_0$, such that*

$$(8.3) \qquad \gamma(x,t) = \exp(\sum_{0 < |\alpha| < N} t^\alpha X_\alpha/\alpha!)(x) + O(|t|^N)$$

*for each positive integer $N$, for all $x \in U$, as $|t| \to 0$.*

Formal Taylor series identities of the form (8.3) will frequently be indicated by the notation "$\sim$":

$$(8.4) \qquad \gamma(x,t) \sim \exp(\sum t^\alpha X_\alpha/\alpha!)(x).$$

Since both the left-hand side and first term on the right in (8.3) are $C^\infty$ functions of $x, t$, so is the term denoted by $O(|t|^N)$. Therefore

$$\frac{\partial^{a+b}}{\partial x^a \partial t^b}\left[\gamma(x,t) - \exp(\sum_{|\alpha| < N} t^\alpha X_\alpha/\alpha!)(x)\right] = O(|t|^{\max(N-b,0)})$$

for any $a, b$, uniformly in $x$ in a neighborhood of $x_0$. Similar considerations apply to other terms denoted by $O(|t|^N)$ later on.

A note of caution is in order: (8.4) does not mean that the integral curves of the vector fields $\sum_{|\alpha|<N} t^\alpha X_\alpha/\alpha!$ are subsets of the manifolds parametrized by $\gamma$; this is false for typical mappings $\gamma$. Consider the example $\gamma(x_1, x_2; t) = (x_1 + t, x_2 + t^2)$. Here $x \in \mathbb{R}^2$ and $t \in \mathbb{R}^1$, so the index $\alpha$ belongs to $\{1, 2, \ldots\}$. The vector fields are $X_1 = \partial/\partial x_1$ and $X_2 = 2\partial/\partial x_2$, while $X_j \equiv 0$ for all $j > 2$. For each $x, t$, the integral curve $s \mapsto \exp(s \sum_\alpha t^\alpha X_\alpha/\alpha!)(x)$ is a parametrized line segment joining $x$ to $\gamma(x, t)$, whereas the curve $t \mapsto \gamma(x, t)$ is a parabola. The endpoint $s = 1$ of the line segment is the point $\gamma(x, t)$ of the parabola.

*Definition* 8.6. $\gamma$ satisfies curvature condition $(\mathcal{C}_\mathfrak{g})$ at $x_0$ if the vector fields $X_\alpha$ together with all their iterated commutators span the tangent space to $\mathbb{R}^n$ at $x_0$.

A third curvature condition is phrased in terms of iterates of the map $t \mapsto \gamma(x, t)$. For any $1 \leq j \leq n$ define $\Gamma^1(x, t) = \gamma(x, t)$ and

$$(8.5) \qquad \Gamma^j(x, t^1, \ldots t^j) = \gamma\left(\Gamma^{j-1}(x, t^1, \ldots t^{j-1}), t^j\right)$$

for $(t^1, \ldots t^j) \in \mathbb{R}^{kj}$ sufficiently close to 0. Among these we single out the $n^{\text{th}}$ iterate

(8.6) $$\Gamma(x, \tau) = \Gamma^n(x, \tau)$$

for $\tau \in \mathbb{R}^{kn}$. The domain of the map $\tau \mapsto \Gamma(x, \tau)$ is a small neighborhood of $0 \in \mathbb{R}^{nk}$; its range is contained in a small neighborhood of $x \in \mathbb{R}^n$.

Write $\tau = (\tau_1, \ldots \tau_{kn})$, where the coordinates belong to $\mathbb{R}^1$ and are ordered in any fixed manner. To each $n$-tuple $\xi = (\xi_1, \ldots \xi_n)$ of elements of $\{1, 2, \ldots kn\}$ is associated the Jacobian determinant

(8.7) $$J_\xi(x, \tau) = \det\left(\frac{\partial \Gamma(x, \tau)}{\partial(\tau_{\xi_1}, \ldots \tau_{\xi_n})}\right)$$

of the $n \times n$ submatrix of the differential of $\Gamma$ with respect to $\tau$. In the next definition, $\partial_\tau^\beta$ represents an arbitrary partial derivative with respect to the full variable $\tau \in \mathbb{R}^{kn}$, not merely $(\tau_{\xi_1}, \ldots, \tau_{\xi_n})$.

*Definition* 8.7.   $\gamma$ satisfies curvature condition $(\mathcal{C}_J)$ at $x_0$ if there exist an $n$-tuple $\xi$ and a multi-index $\beta$ such that

(8.8) $$\left.\partial_\tau^\beta J_\xi(x_0, \tau)\right|_{\tau=0} \neq 0.$$

The notation is a mnemonic device: condition $(\mathcal{C}_M)$ involves a submanifold $M$, $(\mathcal{C}_\mathfrak{g})$ involves a Lie algebra $\mathfrak{g}$, and $(\mathcal{C}_J)$ involves Jacobians $J$.

Each of these three conditions plays a role in our analysis of the operators associated to a mapping $\gamma$. (1) $(\mathcal{C}_J)$ will be used to show that certain associated operators are smoothing. (2) $(\mathcal{C}_\mathfrak{g})$ will ultimately be used to construct a coordinate system in which $(\mathcal{C}_J)$ can be exploited in a systematic way with uniform dependence on certain parameters; in particular, in this coordinate system, $\gamma$ will possess a form of scale– and basepoint–invariance. (3) $(\mathcal{C}_M)$ is useful in a negative sense; given that it fails, it is easy to establish that various other curvature conditions or operator estimates fail to hold.

8.2. *Theorems*. The following theorem is one of the main results of this paper. It will be used in the proofs of our analytic conclusions concerning the mapping properties of certain operators.

THEOREM 8.8.   *For every $\gamma$, the three conditions $(\mathcal{C}_\mathfrak{g}), (\mathcal{C}_M), (\mathcal{C}_J)$ are mutually equivalent at every point.*

*Definition* 8.9.   $\gamma$ is said to be curved to finite order at a point $x_0$, or equivalently to satisfy condition $(\mathcal{C})$ at $x_0$, if it satisfies these three equivalent conditions at $x_0$.

Three additional curvature conditions, $(\mathcal{C}_Y)$, $(\mathcal{C}_\Lambda)$, and $(\mathcal{C}_J)'$, will be introduced in Sections 9 and 10, and will be shown[6] also to be equivalent to $(\mathcal{C})$, as will two analytic properties of integral operators naturally associated to mappings $\gamma$.

The next proposition is a direct consequence of the definitions, but is nonetheless fundamental.

PROPOSITION 8.10. $(\mathcal{C})$ *is invariant under diffeomorphism.*

By this we mean the following.[7] Let $\phi$ be any diffeomorphism of $\mathbb{R}^n$, and $\psi : \mathbb{R}^{n+k} \mapsto \mathbb{R}^k$ any $C^\infty$ map satisfying $\psi(x, 0) \equiv 0$, whose differential with respect to $t$ is invertible at $t = 0$. Given a mapping $\gamma$, set $\tilde\gamma(x, t) = \phi^{-1} \circ \gamma(\phi(x), \psi(x, t))$. Then $\tilde\gamma$ satisfies $(\mathcal{C})$ at $x$ if and only if $\gamma$ satisfies $(\mathcal{C})$ at $\phi(x)$.

A consequence of Proposition 8.10 is that the validity of $(\mathcal{C})$ at $x_0$ does not depend in any way on the geometry or curvature of an individual manifold $M_x$ as defined in (8.2);[8] it is instead a property of the family of manifolds $\{M_x\}$. For it is always possible to choose a diffeomorphism $\phi$ that maps $M_{x_0}$ to an affine subspace $\tilde M$, and then $t \mapsto \tilde\gamma(x_0, t)$ parametrizes $\tilde M$ near $x_0$.

Under any of several mild supplementary hypotheses, $(\mathcal{C})$ has other equivalent formulations, in terms of the analytic properties of associated operators

$$(8.9) \qquad Tf(x) = \int f(\gamma(x,t)) \, K(x,t) \, dt$$

where $K \in C^\infty(\mathbb{R}^n \times \mathbb{R}^k)$. Now $K$ is always assumed to be supported in a small neighborhood of $(x_0, 0)$, and $f$ in a small neighborhood of $x_0$. For each nonnegative $s \in \mathbb{R}$ denote by $H_s$ the usual Sobolev space of all functions defined in some Euclidean space having $s$ derivatives in $L^2$.

THEOREM 8.11. *If $\gamma$ satisfies curvature condition $(\mathcal{C})$ at $x_0$ then there exist $s > 0$ and a neighborhood $U$ of $(x_0, 0)$ such that for every $K \in C^\infty$ supported in $U$, the operator $T$ maps $L^2(\mathbb{R}^n)$ to $H_s$.*

*Conversely, suppose that $\gamma$ does not satisfy curvature condition $(\mathcal{C})$ at $x_0$. Suppose in addition either that the differential of the map $t \mapsto \gamma(x_0, t)$ has maximal rank $k$ at $t = 0$, or that $K \geq 0$. Then there exists a neighborhood $V$ of $(x_0, 0)$ such that whenever $K \in C^\infty$ is supported in $V$ and $K(x_0, 0) \neq 0$, $T$ maps $L^2$ to $H_s$ for no exponent $s > 0$.*

---

[6] A supplementary hypothesis is imposed on $\gamma$ in the formulation of $(\mathcal{C}_\Lambda)$.

[7] Our theory may be generalized to mappings which do not necessarily satisfy $\gamma(x, 0) \equiv x$. The appropriate generalization of diffeomorphism invariance then involves $\phi_1 \circ \gamma(\phi_2(x), \psi(x, t))$, where $\phi_1$ and $\phi_2$ are unrelated.

[8] In this remark we assume for simplicity that the differential of $\gamma$ with respect to $t$ has rank $k$ at $t = 0$, so that each $M_x$ is a manifold.

The next result is a variant involving only the scale of spaces $L^r$.

THEOREM 8.12. *Suppose that $\gamma$ satisfies curvature condition $(\mathcal{C})$ at $x_0$. Then there exists a neighborhood $U$ of $(x_0, 0)$ such that for each $p \in (1, \infty)$ there exists an exponent $q > p$ such that for any $K \in L^\infty$ supported in $U$, the operator $T$ defined by (8.9) maps $L^p(\mathbb{R}^n)$ to $L^q(\mathbb{R}^n)$.*

*Conversely, suppose that there exist $K \in C^0$ satisfying $K(x_0, 0) \neq 0$ and exponents $1 \leq p < q \leq \infty$ such that $T$ maps $L^p(\mathbb{R}^n)$ to $L^q(\mathbb{R}^n)$. Suppose further either that $K \geq 0$, or that the differential of the map $t \to \gamma(x_0, t)$ has rank $k$ at $t = 0$. Then $\gamma$ satisfies curvature condition $(\mathcal{C})$ at $x_0$.*

## Examples

*Example* 8.1. Suppose that $\gamma(x, t) = x - h(t)$ for some $h : \mathbb{R}^k \mapsto \mathbb{R}^n$ defined in a neighborhood of 0 and satisfying $h(0) = 0$. Then $\gamma$ is curved to finite order if and only if the set of all vectors $\{\partial^\alpha h/\partial t^\alpha(0)\}$ spans $\mathbb{R}^n$; this does not depend on $x_0$.

Indeed, the vector fields in the exponential representation for $\gamma$ are

$$X_j \equiv -\partial^j h/\partial t^j \Big|_{t=0},$$

for $j = 1, 2, 3, \ldots$. These have constant coefficients, so that all their commutators vanish identically; consequently $(\mathcal{C}_{\mathfrak{g}})$ holds if and only if they span $\mathbb{R}^n$. For the analysis of this class of examples see [45].

*Example* 8.2. Let $Z_1, \ldots Z_k$ be any collection of $C^\infty$ vector fields, defined in an open subset of $\mathbb{R}^n$, whose iterated commutators span the tangent space to the ambient space $\mathbb{R}^n$ at a point $x_0$. Then $\gamma(x, t) = \exp(\sum t_j Z_j)(x)$ satisfies $(\mathcal{C})$ at $x_0$. Our theory does not require that the $Z_j$ be linearly independent at $x_0$, nor that they span a subspace of constant dimension, nor that $k$ be less than $n$.

*Example* 8.3. There exist mappings curved to finite order of the form $\gamma(x, t) = x + v(x) \cdot t$ where $v : \mathbb{R}^n \mapsto \mathbb{R}^n$, even though each manifold $t \mapsto \gamma(x, t)$ is flat. One such example is $\gamma(x, t) = (x_1 + t, x_2 + x_1 t)$ in $\mathbb{R}^2$, with $t \in \mathbb{R}^1$.

*Example* 8.4.
$$\gamma(x, t) = (x_1 + t, x_2 + 2x_1 t + t^2)$$

appears on a superficial examination to have the features of both of the preceding two examples. Yet it is not curved to finite order at any point, and indeed, may be brought into the form $\tilde{\gamma}(x, t) \equiv (x_1 + t, x_2)$ by conjugation with the diffeomorphism $x \mapsto (x_1, x_2 - x_1^2)$.

It may be tempting to attempt to formulate a curvature condition by examining all partial derivatives of $\gamma$ with respect to $(x, t)$ at $(x_0, 0)$. This example demonstrates that only certain combinations of partial derivatives have any significance. For example, $\partial^2\gamma/\partial x_1\partial t$ is nonzero in this example, yet there is no curvature in our sense. For the first nontrivial example of a polynomial expression in terms of these partial derivatives that does have an invariant meaning, see (9.9).

*Example* 8.5. Let $G$ be a Lie group and $h : \mathbb{R}^k \mapsto G$ a smooth map satisfying $h(0) = 1$, the identity element of $G$. Define $\gamma(g,t) = g \cdot h(t)$. Write $h(t) \sim \exp(\sum_\beta t^\beta Z_\beta/\beta!)$ where each $Z_\beta$ belongs to the Lie algebra $\mathfrak{g}$ of $G$, and exp denotes the exponential mapping from a neighborhood of $0 \in \mathfrak{g}$ to a neighborhood of the identity element of $e \in G$. Each $Z_\beta$ is thus identified with a left invariant vector field on $G$. Then $\gamma$ satisfies $(\mathcal{C}_\mathfrak{g})$ at one point of $G$ if and only if it does so at every point; and this holds if and only if the $Z_\beta$ generate the entire Lie algebra $\mathfrak{g}$.

Indeed, the vector fields $X_\alpha$ in the exponential representation for $\gamma$ are precisely the $Z_\alpha$. These generate a certain subalgebra of locally left invariant vector fields in a neighborhood of $e$, that is, a certain subalgebra of $\mathfrak{g}$. Thus $(\mathcal{C}_\mathfrak{g})$ holds at $e$ if and only this subalgebra is all of $\mathfrak{g}$. The rank of this subalgebra is constant, so that $(\mathcal{C}_\mathfrak{g})$ holds at $e$ if and only if it holds at every point of $G$.

*Example* 8.6. The real analytic mapping $\gamma(x,t) = (x_1 + t, x_2 + x_2 t)$ is curved to finite order at some points but not all. It is curved to finite order at a point $(x_1, x_2)$ if and only if $x_2 \neq 0$. The manifold $\{x_2 = 0\}$ is invariant under $\gamma$.

*Example* 8.7. The mapping $\gamma(x,t) = (x_1+t, x_2+\exp(-t^{-2}))$ is not curved to finite order at any point. The manifolds $\{x_2 = \text{constant}\}$ are invariant under $\gamma$ to infinite order at every point, but there exists no manifold of positive codimension that is truly locally invariant under $\gamma$.

Similarly for $\gamma(x,t) = (x_1 + t, x_2 + t\exp(-x_1^{-2}))$, the manifolds $\{x_2 = \text{constant}\}$ are invariant to infinite order at each point where $x_1 = 0$, but there exist no locally invariant submanifolds.

*Example* 8.8. Phong and Stein [31] have discussed a notion of rotational curvature for the hypersurface case $k = n - 1$, assuming the differential of $\gamma$ with respect to $t$ to be injective. Given a point $x^0 \in \mathbb{R}^n$, it is always possible to introduce coordinates $x = (x', x_n) \in \mathbb{R}^{n-1} \times \mathbb{R}$ with origin at $x^0$, and a smooth function $\phi : \mathbb{R}^{n+(n-1)} \mapsto \mathbb{R}$ satisfying $\nabla\phi(0,0) = 0$, such that for any points $x, y \in \mathbb{R}^n$ near the origin, $y = \gamma(x,t)$ for some $t \in \mathbb{R}^k$ near 0 if and only if $x_n = y_n + \phi(x, y')$. The rotational curvature condition of [31] is then equivalent to the nonvanishing of the determinant of the matrix

$$\det\left(\frac{\partial^2\phi}{\partial x_i \partial y_j}\right)_{1\leq i,j\leq n-1}(0,0) \neq 0.$$

On the other hand, our condition ($\mathcal{C}$) is equivalent to the requirement merely that there exist indices $i, j \leq n-1$ for which the Taylor expansion of $\partial^2 \phi / \partial x_i \partial y_j$ with respect to $(x', y') = (x_1, \ldots x_{n-1}, 0, y_1, \ldots y_{n-1}, 0)$ does not vanish to infinite order at the origin; see (9.18) and Proposition 9.15. Thus rotational curvature in the sense of [31] is in a natural sense the strongest possible version of ($\mathcal{C}$), in the hypersurface case.

*Example* 8.9. A condition closely related to ($\mathcal{C}$) in a special case has been studied by Baouendi, Ebenfelt and Rothschild [2]. Let $X \subset \mathbb{C}^n$ be a real analytic, generic CR manifold. In local coordinates in a neighborhood $\Omega \subset \mathbb{C}^n$ of a point $z_0 \in X$, $X = \{z \in \Omega : \rho(z, \bar{z}) = 0\}$, where $\rho : \Omega \times \overline{\Omega} \mapsto \mathbb{C}$ is holomorphic, satisfies $\rho(z, \bar{w}) \equiv \overline{\rho(w, \bar{z})}$, and its components satisfy $\partial \rho_1 \wedge \partial \rho_2 \wedge \ldots \wedge \partial \rho_d \neq 0$. Here $\partial \rho_j(z, \bar{z}) = \sum_i \frac{\partial h}{\partial z_i} dz_i$ where $h(z) = \rho_j(z, \bar{z})$.

To each $z \in X$ are associated the Segre sets $N_j \subset \mathbb{C}^n$ defined by $N_0(z) = \{z\}$, and $N_{j+1}(z) = \{w \in \Omega : \rho(w, \bar{\zeta}) = 0 \text{ for some } \zeta \in N_j(z)\}$; in these definitions $z, w, \zeta$ are assumed to lie sufficiently close to $z_0$.

Define $\Lambda \subset \Omega \times \Omega$ to be $\{(z, w) : \rho(z, \bar{w}) = 0\}$, and parametrize $\Lambda$ as $\{(\gamma(z, t), z)\}$, where $t \in \mathbb{C}^{n-d} = \mathbb{C}^k$, and so that $\gamma(z, 0) = z$ for all $z \in X$ (though not necessarily for all $z \in \Omega$; in this respect [2] differs from our framework) for some analytic function $\gamma$. Then in our notation, $N_j$ corresponds to the image of $\mathbb{C}^{jk}$ under the map $(t_1, \ldots t_j) \in \mathbb{C}^{jk} \mapsto \Gamma^j(z, (t_1, \ldots t_j))$. One of the geometric conditions considered in [2] is that some $N_j$ should contain an open subset of $\mathbb{C}^n$. In our language this means that the Jacobian determinant of $\Gamma(z_0, \tau)$ with respect to $\tau$ does not vanish identically. Because of the hypothesis of analyticity, this is equivalent to the nonvanishing of some coefficient of its Taylor expansion about $\tau = 0$, that is, to ($\mathcal{C}_J$) at points of $X$.

## 9. Curvature: Some details

9.1. *The exponential representation.* In this subsection we establish the existence and uniqueness of the vector fields $X_\alpha$ in the exponential representation (8.3).

We will often work with $\mathbb{R}^n$-valued formal Taylor series $\exp(\sum_\alpha t^\alpha X_\alpha / \alpha!)(x)$ with respect to $t \in \mathbb{R}^k$, in which the $X_\alpha$ are $C^\infty$ vector fields defined in an open subset $U \subset \mathbb{R}^n$. These are to be interpreted as follows. The mappings $F_N(t) = \exp(\sum_{|\alpha| \leq N} t^\alpha X_\alpha / \alpha!)(x)$ are well defined for all $|t| < \delta(N) > 0$, for all $x$ belonging to any relatively compact subset of $U$. With respect to a coordinate system in $\mathbb{R}^n$, $F_N$ may then be expanded in a Taylor series $T_N(t) = \sum_\beta c_\beta^N(x) t^\beta / \beta!$ with respect to $t$, whose coefficients $c_\beta^N$ are $C^\infty$ functions of $x$. For any multi-index $\beta$, $c_\beta^N$ is independent of $N$ for all sufficiently large $N$.

Defining $c_\beta(x)$ to be this stable limit, we obtain a formal series

$$(9.1) \qquad T_\infty \sim \sum_\beta c_\beta(x) t^\beta/\beta! \ .$$

The coefficients of this series are $\mathbb{R}^n$-valued $C^\infty$ functions of $x \in U$.

*Definition* 9.1. Let $\{X_\alpha\}$ be a collection of $C^\infty$ vector fields defined in a neighborhood $U$ of $x_0 \in \mathbb{R}^n$. Then $\exp\left(\sum_\alpha t^\alpha X_\alpha/\alpha!\right)(x)$ is defined to be the formal Taylor series (9.1) in $t$.

We will also work with exponential mappings associated to vector fields depending on a parameter $t$, such as $\exp(\sum_{|\alpha|<N} t^\alpha X_\alpha/\alpha!)(x)$; such exponentials are of course well defined $C^\infty$ functions of $x,t$ for small $t$, not merely formal power series.

*Proof of existence.* Define a mapping $H : \mathbb{R}^{n+k} \mapsto \mathbb{R}^{n+k}$ by

$$(9.2) \qquad H(x,t) = (\gamma(x,t), t).$$

Because $\gamma(x,0) \equiv x$, the differential of $H$ at any point $(x,0)$ is invertible, whence $H$ defines a diffeomorphism of a small neighborhood $U_1$ of $(x_0,0)$ in $\mathbb{R}^{n+k}$ with a second such neighborhood $U_2$. For small $s \in \mathbb{R}^k$ define

$$(9.3) \qquad \varphi_s(x,t) = H(\gamma_t^{-1}(x), s+t) = (\gamma_{s+t}\gamma_t^{-1}(x), s+t).$$

Then $\varphi_{s_1+s_2} \equiv \varphi_{s_1} \circ \varphi_{s_2}$. As remarked in Section 1, this means that there exist vector fields $\{V_j : 1 \le j \le k\}$ in $\mathbb{R}^{n+k}$, defined near $(x_0,0)$, satisfying $\varphi_s(x,t) \equiv \exp(\sum_j s_j V_j)(x,t)$. In particular,

$$(9.4) \qquad \exp(\sum_{j=1}^k s_j V_j)(x,0) = \varphi_s(x,0) = H(x,s) = (\gamma(x,s), s).$$

Differentiating with respect to $s$ gives

$$(9.5) \qquad V_j = H_*(\partial/\partial t_j),$$

the push forward of $\partial/\partial t_j$ by $H$. This could equally well have been taken as an alternative definition of $V_j$.

Adopting coordinates $(y,s) \in \mathbb{R}^{n+k}$ in $U_2$, we may express each $V_j$ as $Z_j + \partial_{s_j}$ where $Z_j$ is a vector field in $U_2$ whose coefficients are functions of $(y,s)$, but which is everywhere in the span of $\partial_{y_1}, \ldots \partial_{y_n}$. Substituting (9.4) into the identity

$$(\gamma(x,t), 0) = \exp(-\sum_j t_j \partial/\partial s_j)(\gamma(x,t), t)$$

and applying the Baker-Campbell-Hausdorff formula yield a Taylor series identity in $t$ of the form

$$(9.6) \quad (\gamma(x,t), 0) = \exp(-\sum_j t_j \partial/\partial s_j) \exp(\sum_j t_j [Z_j + \partial/\partial s_j])(x, 0)$$
$$\sim \exp\Big( \sum_{|\beta|>0} t^\beta Z_\beta/\beta! \Big)(x, 0),$$

where each $Z_\beta$ is a finite linear combination, with constant coefficients, of the vector fields $Z_i$ and $\partial/\partial s_j$ and of their commutators.

Each such commutator, excepting the $\partial/\partial s_j$ themselves, is in the span of the coordinate vector fields $\partial/\partial y_i$, with coefficients that are $C^\infty$ functions of $y, s$. The coefficient of $\partial/\partial s_j$ is $t_j - t_j \equiv 0$. Therefore each $Z_\alpha$ takes the form

$$Z_\alpha(y, s) = \sum_{j=1}^n a_{j,\alpha}(y, s) \partial/\partial y_j.$$

Define vector fields $X_\alpha$ on $\mathbb{R}^n$ by

$$(9.7) \qquad X_\alpha(y) = \sum_j a_{j,\alpha}(y, 0) \partial/\partial y_j.$$

Then for any $N$, $\exp(\sum_{|\alpha|\leq N} t^\alpha X_\alpha/\alpha!)(x)$ is identically equal to the $\mathbb{R}^n$ coordinate of the quantity $\exp(\sum_{|\alpha|\leq N} t^\alpha Z_\alpha/\alpha!)(x, 0)$, by the uniqueness of solutions to the initial value problem for ordinary differential equations. Thus (9.6) may be rewritten as

$$\gamma(x, t) \sim \exp(\sum_{|\alpha|>0} t^\alpha X_\alpha/\alpha!)(x). \qquad \square$$

*Proof of uniqueness.* The reasoning is by induction on $N$. Suppose that for some vector fields $Y_\alpha$,

$$\exp(\sum_{|\alpha|\leq N} t^\alpha X_\alpha/\alpha!)(x)$$
$$= \exp\Big( \sum_{|\alpha|<N} t^\alpha X_\alpha/\alpha!)(x) + \sum_{|\alpha|=N} t^\alpha Y_\alpha/\alpha! \Big)(x) + O(|t|^{N+1})$$

for all $x$ in some neighborhood of $x_0$. Composing with $\exp(-\sum_{|\alpha|\leq N} t^\alpha X_\alpha/\alpha!)$ on both sides yields

$$x + O(|t|^{N+1})$$
$$= \exp(-\sum_{|\alpha|\leq N} t^\alpha X_\alpha/\alpha!) \exp\Big( \sum_{|\alpha|<N} t^\alpha X_\alpha/\alpha! + \sum_{|\alpha|=N} t^\alpha Y_\alpha/\alpha! \Big)(x).$$

The Baker-Campbell-Hausdorff formula may be used to simplify the right-hand side, yielding

$$x + O(|t|^{N+1}) = \exp\left(\sum_{|\alpha|=N} t^\alpha (Y_\alpha - X_\alpha)/\alpha!\right)(x) + O(|t|^{N+1}).$$

Fix any multi-index $\delta$ satisfying $|\delta| = N$, apply $\partial^\delta/\partial t^\delta$ to both sides, and evaluate at $t = 0$ to obtain $0 = \delta!\,(Y_\alpha - X_\alpha)(x)$ for all $x$ near $x_0$. □

*Remark* 9.1. From this proof of uniqueness can also be derived an alternative construction of the vector fields $X_\alpha$, by induction on $|\alpha|$. For $|\alpha| = 1$, $X_\alpha = \partial \gamma(x,t)/\partial t_j$ where $\alpha = (0, \ldots, 1, \ldots 0)$ with a single 1 in the $j^{\text{th}}$ position. Once the $X_\alpha$ have been constructed for all $|\alpha| \leq N$, consider

$$G(x,t) = \exp(-\sum_{|\alpha|\leq N} t^\alpha X_\alpha/\alpha!)(\gamma(x,t)).$$

It can be shown that for each multi-index satisfying $|\beta| = N+1$, $\partial^\beta G/\partial t^\beta(x,0)$ defines a vector field on $\mathbb{R}^n$, which is taken to be $X_\beta$.

Another approach to the computation of the $X_\alpha$ is as follows. Fix coordinates, write $\gamma(x,t) = (\gamma_1, \ldots \gamma_n)$, and denote by $e_j$ the $j^{\text{th}}$ coordinate function in $\mathbb{R}^n$. According to (??), the Taylor expansion of the component $\gamma_j$ with respect to $t$ about $x$ is given by the formal series

$$\sum_{k=0}^{\infty} \frac{1}{k!}\left[\sum_\beta t^\beta X_\beta/\beta!\right]^k (e_j)(x).$$

Expanding this formal double series in $t$, the coefficient of $t^\alpha$ is

$$\sum_{m=1}^{|\alpha|} \sum_{\beta_1+\cdots+\beta_m=\alpha} \frac{1}{m!\beta!} X_{\beta_1} \cdots X_{\beta_m}(e_j)(x).$$

This says that

$$\frac{1}{\alpha!}\frac{\partial^\alpha \gamma_j}{\partial t^\alpha}(x,0) = \frac{1}{\alpha!} X_\alpha(e_j) + \sum_{m=2}^{|\alpha|} \sum_{\beta_1+\cdots+\beta_m=\alpha} \frac{1}{m!\beta!} X_{\beta_1} \cdots X_{\beta_m}(e_j)(x).$$

The final conclusion is a recursive formula.

LEMMA 9.2. *For every $\alpha$,*

$$(9.8) \quad X_\alpha(e_j) = \frac{\partial^\alpha \gamma_j}{\partial t^\alpha}(x,0) - \sum_{m=2}^{|\alpha|} \sum_{\beta_1+\ldots+\beta_m=\alpha} \frac{\alpha!}{m!\beta!} X_{\beta_1} \cdots X_{\beta_m}(e_j)(x).$$

Implicit in this computation are proofs of uniqueness and existence of $\{X_\alpha\}$. One consequence is that $X_\alpha(x)$ depends only on those Taylor coefficients $\partial^{\beta_1+\beta_2}\gamma/\partial x^{\beta_1}\partial t^{\beta_2}(x,0)$ of total order $|\beta| \leq |\alpha|$. One has of course

$$X_\alpha(x) = \sum_{j=1}^n \frac{\partial \gamma_j}{\partial t^\alpha}(x,0)\frac{\partial}{\partial x_j} \qquad \text{for } |\alpha|=1.$$

In the simplest case $k=1$ so that $t \in \mathbb{R}$, the second order vector field is

$$(9.9) \quad X_2(x) = \sum_{j=1}^n \left[\frac{\partial^2 \gamma_j}{\partial t^2}(x,0) - 2\sum_{m=1}^n \frac{\partial \gamma_m}{\partial t}(x,0)\frac{\partial^2 \gamma_j}{\partial x_m \partial t}(x,0)\right]\frac{\partial}{\partial x_j}.$$

9.2. *Diffeomorphism invariance.* For any $\gamma$, there exists a unique, smooth mapping $\gamma^{-1}$ satisfying $\gamma(\gamma^{-1}(x,t),t) \equiv x \equiv \gamma^{-1}(\gamma(x,t),t)$, by the implicit function theorem.

LEMMA 9.3. *The inverse mapping $\gamma^{-1}$ satisfies*

$$(9.10) \qquad \gamma^{-1}(x,t) \sim \exp(-\sum_\alpha t^\alpha X_\alpha/\alpha!)(x)$$

*for the same vector fields $X_\alpha$ occurring in the exponential representation for $\gamma$.*

A consequence is that $\gamma$ satisfies condition $(\mathcal{C}_\mathfrak{g})$ at $x_0$ if and only if $\gamma^{-1}$ does so.

*Proof.* Consider the identity

$$\exp(\sum_{|\alpha|\leq N} t^\alpha X_\alpha/\alpha!)\exp(\sum_{|\alpha|\leq N} -t^\alpha X_\alpha/\alpha!)(x) \equiv x\ .$$

By letting $N \to \infty$ we conclude that

$$\gamma\Big(\exp(-\sum_\alpha t^\alpha X_\alpha/\alpha!)(x),t\Big) \sim x$$

as Taylor series in $t$. Since the equation $\gamma(\gamma^{-1}(x,t),t)$ uniquely determines the Taylor expansion of $\gamma^{-1}$ about $t=0$, the expansion with respect to $t$ of $\gamma^{-1}$ equals the formal series $\exp(-\sum_\alpha t^\alpha X_\alpha/\alpha!)(x)$. □

To any diffeomorphism $\phi$ from an open subset of $\mathbb{R}^n$ to a second open subset is associated the push-forward mapping $\phi_*$, which transforms vector fields defined on the domain of $\phi$ to vector fields on its range. It is defined by $\phi_*(X)(\phi(x)) = D\phi(X(x))$, where $D\phi(x)$ denotes the differential of $\phi$ at $x$.

LEMMA 9.4. *Let $\gamma : \mathbb{R}^{n+k} \mapsto \mathbb{R}^n$ be $C^\infty$ and satisfy $\gamma(x,0) \equiv x$, and let $\{X_\alpha\}$ be the vector fields in the exponential representation (8.3) for $\gamma$. Let $\phi$ be a diffeomorphism of $\mathbb{R}^n$ with itself, and define $\tilde{\gamma}(x,t) = \phi(\gamma(\phi^{-1}x,t))$. Then*

$$(9.11) \qquad \tilde{\gamma}(x,t) \sim \exp\Big(\sum_\alpha t^\alpha \cdot \phi_* X_\alpha/\alpha!\Big)(x).$$

Coupling this result with the uniqueness of the exponential representation, we conclude that the vector fields in that representation are themselves invariant under diffeomorphism, in a natural sense.

*Proof.* For any vector field $Y$ and any diffeomorphism,

$$\phi(\exp(Y)(x)) = \exp(\phi_* Y)\big(\phi(x)\big),$$

as follows from the definition of the exponential mapping by replacing $Y$ by $sY$ and differentiating with respect to $s$. Applying this to $Y = \sum_{|\alpha| \le N} t^\alpha X_\alpha/\alpha!$ for each $N$ yields the stated conclusion. □

*Proof of diffeomorphism invariance.* We demonstrate here invariance of $(\mathcal{C}_M)$, and invariance of $(\mathcal{C}_J)$ and of $(\mathcal{C}_\mathfrak{g})$ in the special case $\psi(x,t) \equiv t$. Some of these facts will be used in the proof of equivalence of the curvature conditions. Full diffeomorphism invariance of all three will then be a consequence of their equivalence.

For $(\mathcal{C}_\mathfrak{g})$, consider $\tilde\gamma(x,t) = \phi^{-1} \circ \gamma(\phi(x),t)$. Let $\{X_\alpha\}$ be the vector fields associated to $\gamma$. By (9.11), $\tilde\gamma$ satisfies $(\mathcal{C}_\mathfrak{g})$ at $x_0$ if and only if the Lie algebra generated by $\{\phi_*^{-1} X_\alpha\}$ spans the tangent space to $\mathbb{R}^n$ at $x_0$. The push forward operation is a Lie algebra homomorphism, so this Lie algebra equals the image under $\phi_*^{-1}$ of the Lie algebra generated by $\{X_\alpha\}$. Thus $\tilde\gamma$ satisfies $(\mathcal{C}_\mathfrak{g})$ at $x_0$ if and only if $\gamma$ does so at $\phi(x_0)$.

Consider condition $(\mathcal{C}_M)$ first in the simpler case $\tilde\gamma(x,t) = \phi^{-1} \circ \gamma(\phi(x),t)$. If a submanifold $M$ is invariant to infinite order under $\gamma$ at $x_0$, then $\phi^{-1}(M)$ is invariant to infinite order under $\tilde\gamma(x,t)$ at $\phi^{-1}(x_0)$.

Next consider $\tilde\gamma(x,t) = \gamma(x,\psi(x,t))$. Since $|\psi(x,t)| \sim |t|$, again condition $(\mathcal{C}_M)$ is manifestly invariant. Combining these two cases by composition yields invariance of $(\mathcal{C}_M)$ in the general case.

To prove that $(\mathcal{C}_J)$ holds for $\gamma$ if and only if it holds for $\tilde\gamma_t(x) = \phi^{-1}\gamma_t(\phi x)$, note that the $n$-fold iterate $\tilde\Gamma$ of $\tilde\gamma$ satisfies $\tilde\Gamma(x,\tau) = \phi^{-1}\Gamma(\phi x, \tau)$. Thus $\tilde J_\xi(x,\tau)$ is the product of the determinant of $D\phi^{-1}$ at $\phi(x)$ with $J_\xi(\phi x, \tau)$. Consequently $\partial_\tau^\beta \tilde J_\xi(x,0)$ equals $\partial_\tau^\beta J_\xi(\phi x, 0)$ times a nonvanishing function of $x$; if one is nonzero, so is the other. □

9.3. *Curvature condition* $(\mathcal{C}_Y)$. A fourth curvature condition is phrased in terms of a different collection of vector fields. For $1 \le j \le k$ define

(9.12) $$Y_j(x,t) = \frac{\partial}{\partial s_j} \gamma_{t+s}\gamma_t^{-1}(x)\Big|_{s=0};$$

each of these is a vector field in $\mathbb{R}^n$ evaluated at $x$, depending smoothly on the parameter $t \in \mathbb{R}^k$. Define vector fields $Y_{\alpha,j}(x)$ to be the Taylor coefficients

(9.13) $$Y_j(x,t) \sim \sum_{|\alpha| \ge 0} \frac{t^\alpha}{\alpha!} Y_{\alpha,j}(x).$$

*Definition* 9.5. $\gamma$ satisfies curvature condition $(\mathcal{C}_Y)$ at $x_0$ if the vector fields $Y_{\alpha,j} : 1 \leq j \leq k, |\alpha| \geq 0$, together with all their iterated commutators, span the tangent space to $\mathbb{R}^n$ at $x_0$.

PROPOSITION 9.6. *Curvature condition $(\mathcal{C}_Y)$ is equivalent to $(\mathcal{C}_\mathfrak{g})$.*

*Proof.* For any set $S$ of vector fields, denote by $\mathfrak{L}S$ the Lie algebra generated by $S$, as a $C^\infty(\mathbb{R}^n)$ module. Fix any point $x_0$. Fix $j$ and let $e_j \in \mathbb{R}^k$ be the $j^{\text{th}}$ coordinate vector. Let $s = \sigma e_j$ where $\sigma \in \mathbb{R}$.

By the Baker-Campbell-Hausdorff formula, in the sense of formal Taylor series in $t, s$,

$$\gamma_{t+s}\gamma_t^{-1}(x) \sim \exp(\sum_\alpha (s+t)^\alpha X_\alpha/\alpha!) \exp(-\sum_\alpha t^\alpha X_\alpha/\alpha!)(x)$$
$$\sim \exp(\sigma W)(x) + O(|\sigma|^2)$$

where

$$W \sim \sum_\alpha \alpha_j t^{\alpha-e_j}(X_\alpha + V_{\alpha,j})/\alpha!,$$

and each $V_{\alpha,j}$ belongs to $\mathfrak{L}\{X_\beta : |\beta| < |\alpha|\}$. Only multi-indices $\alpha$ whose $j^{\text{th}}$ components are $\geq 1$ appear in the Taylor expansion of $W$.

Thus $Y_{\delta,j}$ equals a nonzero coefficient times $X_{\delta+e_j} + V_{\delta+e_j}$, for each $\delta, j$. Therefore for each $N$,

$$\mathfrak{L}\{Y_{\delta,j} : |\delta| \leq N, \ 1 \leq j \leq k\} = \mathfrak{L}\{X_\alpha : |\alpha| \leq N+1\}.$$

In particular, the Lie algebra generated by all $Y_{\delta,j}$ always coincides with that generated by all $X_\alpha$. Therefore $(\mathcal{C}_Y) \Leftrightarrow (\mathcal{C}_\mathfrak{g})$. $\square$

9.4. *Two lemmas.* Invariance of a submanifold may be conveniently reformulated in a special local coordinate system. Consider coordinates $(x', x'') = x \in \mathbb{R}^p \times \mathbb{R}^{n-p}$, for some $0 < p < n$. Corresponding to this decomposition is an expression

(9.14) $$\gamma(x,t) = \gamma\big((x',x''),t\big) = \big(\gamma'(x,t), \gamma''(x,t)\big),$$

where $\gamma'$ maps $\mathbb{R}^{n+k}$ to $\mathbb{R}^p$ and $\gamma'' : \mathbb{R}^{n+k} \mapsto \mathbb{R}^{n-p}$. Consider the submanifold

$$M = \{x : x'' = 0\}.$$

LEMMA 9.7. *Let $\gamma$ take the form (9.14). Then $M$ is invariant under $\gamma$ to infinite order at the origin $0 \in \mathbb{R}^n$ if and only if*

$$\frac{\partial^{\alpha+\beta}\gamma''}{\partial(x')^\alpha \partial t^\beta}\big((0,0),0\big) = 0 \ \text{ for every } \alpha, \beta.$$

The conclusion may be rephrased: $M$ is invariant under $\gamma$ to infinite order at the origin if and only if

(9.15) $$\gamma(x,t) \sim \big(x' + O(t), \ x'' + O(x'')O(t)\big).$$

*Proof.* The distance from $\gamma(x,t)$ to $M$ is comparable to $|\gamma''(x,t)|$. Thus the restriction to $M \times \mathbb{R}^k$ of the function $(x,t) \mapsto$ distance $(\gamma(x,t), M)$ is comparable to the function $(x',t) \mapsto |\gamma''\big((x',0),t\big)|$. □

The next result will be used in Section 21 to prove a theorem on differentiation of integrals.

LEMMA 9.8. *Suppose that $\gamma$ is real analytic near $(x_0, 0)$. Then a small neighborhood $U$ of $x_0$ may be decomposed as $E \cup V$ where $E$ is a closed set of Lebesgue measure zero, while the open set $V$ is foliated by real analytic leaves with the two properties*: (i) *each leaf is invariant under $\gamma$, and* (ii) *the restriction to each leaf of $\gamma$ satisfies $(\mathcal{C}_\mathfrak{g})$, relative to that leaf.*

By (ii) we mean that the Lie algebra generated by $\{X_\alpha\}$ spans the tangent space to each leaf.

*Proof.* If $\gamma$ satisfies $(\mathcal{C}_\mathfrak{g})$ at $x_0 \in \mathbb{R}^n$ then let $E$ be the empty set and $V = U$. Otherwise fix a small neighborhood $U$ of $x_0$ and $\delta > 0$ such that the restriction of $\gamma$ to $U \times \{|t| < \delta\}$ is analytic.

Let $\{X_\alpha\}$ be the collection of vector fields associated to $\gamma$, and let $\mathfrak{g}$ be the Lie algebra it generates. Define $r \leq n$ to be the maximum, over all $x \in U$, of the dimension of the subspace of $T_x \mathbb{R}^n$ spanned by $\mathfrak{g}$. Define $E$ to be the set of all $x \in U$ at which this subspace has dimension strictly less than $r$. Then $E$ is an analytic variety, hence is closed and has measure zero.

If $r = n$ then there is a trivial foliation of $V = U \setminus E$ in which every leaf is an open subset of $\mathbb{R}^d$. Otherwise in a neighborhood of any $x \in V$, $\mathfrak{g}$ spans a subspace of dimension $r$. By the Frobenius theorem, $V$ is foliated by integral leaves of $\mathfrak{g}$. Because $\mathfrak{g}$ is spanned by analytic vector fields, each leaf is analytic.

Each vector field $X_\alpha$ is everywhere tangent to each leaf. Hence for any $N$ and any leaf $L$,
$$\exp\Big(\sum_{|\alpha| \leq N} t^\alpha X_\alpha / \alpha!\Big)(x) \in L$$
for every $x \in L$ and every sufficiently small $t$. Therefore for any $x \in L$, the distance from $\gamma(x,t)$ to $L$ is $O(|t|^N)$ for every $N$. Because $L$ and $\gamma$ are analytic, by letting $N \to \infty$ we deduce that $L$ is invariant under $\gamma$.

Fix a leaf $L$ and regard $\gamma$ as a mapping from $L \times \mathbb{R}^k$ to $L$. The vector fields occurring in the exponential representation for this restricted mapping are equal to the restrictions of the $X_\alpha$ to $L$, by the uniqueness of the exponential representation. The Lie algebra generated by these is the restriction to $L$ of $\mathfrak{g}$, which spans the tangent space to $L$ at each of its points by construction. Therefore the restriction to $L$ of $\gamma$ satisfies curvature condition $(\mathcal{C}_\mathfrak{g})$. □

9.5. *Double fibration formulation.* Smoothly varying families of $k$ dimensional submanifolds of $\mathbb{R}^n$ are sometimes described in a nonparametrized form, as follows. For $j = 1, 2$ denote by $\tilde{\pi}_j : \mathbb{R}^n \times \mathbb{R}^n \mapsto \mathbb{R}^n$ the projection onto the $j^{\text{th}}$ factor, which we sometimes denote by $\mathbb{R}^n_j$ for clarity. Suppose given a smooth submanifold $\Lambda \subset \mathbb{R}^n \times \mathbb{R}^n$ of codimension $d = n - k$, containing the diagonal $\Delta = \{(x, y) : x = y\}$. Denote by $\pi_j : \Lambda \mapsto \mathbb{R}^n$ the restrictions to $\Lambda$ of the projections $\tilde{\pi}_j$.

Suppose that $D\pi_j : T\Lambda \mapsto T\mathbb{R}^n$ has rank $n$ at every point of a neighborhood of $(x_0, x_0)$. Throughout the discussion, all points are implicitly assumed to be sufficiently close to $x_0$ or to $(x_0, x_0)$. For $x, y \in \mathbb{R}^n$ define

$$M_x = \{y : (x, y) \in \Lambda\} \text{ and } M'_y = \{x : (x, y) \in \Lambda\}.$$

By the implicit function theorem, the assumption on the rank of $D\pi_j$ implies that $M_x, M'_y$ are smooth submanifolds of dimension $k$, depending smoothly on $x, y$ respectively.

There exist parametrizations $t \mapsto \gamma(x, t)$, such that

(9.16) $$\Lambda = \{(x, \gamma(x, t))\}$$

as $(x, t)$ varies over a small neighborhood of $(x_0, 0) \in \mathbb{R}^{n+k}$. Conversely, given any mapping $\gamma$ satisfying our usual hypothesis that $\gamma(x, 0) \equiv x$ and the additional requirement that the differential of $t \mapsto \gamma(x, t)$ has rank $k$, the set $\Lambda$ defined by (9.16) has the above properties near $(x_0, x_0)$. For any set $E \subset \mathbb{R}^n$,

$$\{\gamma(x, t) : x \in E, t \in \mathbb{R}^k\} = \pi_2 \pi_1^{-1}(E).$$

The purpose of this subsection is to discuss an intrinsic reformulation[9] $(\mathcal{C}_\Lambda)$ of $(\mathcal{C})$ in terms of $\Lambda$ itself, which involves no parametrization. This material will not be used elsewhere in the paper.

*Definition* 9.9. $\mathcal{V}_1$ is the $C^\infty(\Lambda)$ module consisting of all sections of the tangent bundle $T\Lambda$ that are in the nullspace of $D\pi_2$. Likewise $\mathcal{V}_2$ is the $C^\infty(\Lambda)$ module consisting of all sections of the tangent bundle $T\Lambda$ that are in the nullspace of $D\pi_1$.

Thus $\mathcal{V}_1$ is the set of all vector fields defined on $\Lambda$, tangent to $\Lambda$, and expressible as linear combinations of $\{\partial/\partial x_j\}$, with coefficients depending on $(x, y) \in \Lambda$. $\mathcal{V}_1$ is closed under Lie brackets. At each point of $\Lambda$, the subspace of $T\Lambda$ spanned by $\mathcal{V}_1$ has dimension $k$, by the hypothesis on the rank of $D\pi_2$. The leaves of the foliation defined by $\mathcal{V}_1$ have dimension $k$, are contained in $\Lambda$, and project to points under $\pi_2$. Therefore they are simply the manifolds $M'_p \times \{p\}$ where $p \in \mathbb{R}^n$ is arbitrary. Corresponding remarks apply to $\mathcal{V}_2$, with the roles of the indices $1, 2$ reversed; its leaves are the manifolds $\{p\} \times M_p$.

---

[9] This version of the curvature condition has been introduced independently by Seeger [38], in a more refined form.

*Definition* 9.10. (i) $\mathfrak{G}(\Lambda)$ denotes the Lie subalgebra of $T\Lambda$ generated by $\mathcal{V}_\infty \cup \mathcal{V}_\in$, as a $C^\infty(\Lambda)$ module.

(ii) $\mathfrak{g}(\Lambda)$ is the $C^\infty(\mathbb{R}^n)$ submodule of sections of $T\mathbb{R}^n$ given by $D\pi_1(\mathfrak{G}(\Lambda)|_\Delta)$.

*Definition* 9.11. $\gamma$ satisfies curvature condition $(\mathcal{C}_\Lambda)$ at $x_0$ if $\mathfrak{G}(\Lambda)$ spans the entire tangent space to $\Lambda$ at $(x_0, x_0)$.

LEMMA 9.12. *$\Lambda$ satisfies $(\mathcal{C}_\Lambda)$ at $\delta \in \Delta$ if and only if $\mathfrak{g}(\Lambda)$ spans the tangent space to $\mathbb{R}^n$ at $\pi_1(\delta)$.*

*Proof.* Let $\delta \in \Delta$. The dimension $\rho$ of the span $\mathfrak{G}(\Lambda)(\delta)$ of $\mathfrak{G}(\Lambda)$ at $\delta$ equals the rank $r$ of the restriction to $\mathfrak{G}(\Lambda)(\delta)$ of $D\pi_1$, plus the dimension $\nu$ of the nullspace of $D\pi_1$ on $\mathfrak{G}(\Lambda)(\delta)$. By the definition of $\mathfrak{g}(\Lambda)$, $r$ equals the dimension of its span at $\pi_1(\delta)$. On the other hand $\nu = k$, because the span at $\delta$ of $\mathcal{V}_2$ is contained in $\mathfrak{G}(\Lambda)(\delta)$, and equals the nullspace of $D\pi_1$ as a transformation from $T_\delta \Lambda$ to $\mathbb{R}^n$; this last nullspace has dimension $k$ by hypothesis on $\Lambda$. Thus $\rho = r + k$, and consequently $\rho = n + k$ if and only if $r = n$. □

Let $\gamma : \mathbb{R}^{n+k} \mapsto \mathbb{R}^n$ be a $C^\infty$ mapping satisfying $\gamma(x, 0) \equiv x$, whose differential with respect to $t$ has rank $k$ at $x_0$. To $\gamma$ is associated $\Lambda_\gamma = \{(x, \gamma(x, t)\} \subset \mathbb{R}^{n+n}$, where $(x, t)$ ranges over a small neighborhood of $(x_0, 0)$. The main result of this section is another equivalence of curvature conditions.

PROPOSITION 9.13. *Let $\gamma : \mathbb{R}^{n+k} \mapsto \mathbb{R}^n$ be a smooth mapping satisfying $\gamma(x, 0) \equiv x$, whose differential with respect to $t$ has rank $k$ at $(x_0, 0)$. Then $\gamma$ satisfies $(\mathcal{C})$ at $x_0$ if and only if the associated manifold $\Lambda_\gamma$ satisfies $(\mathcal{C}_\Lambda)$ at $(x_0, x_0)$.*

The proof of one half of this proposition is deferred to Lemma 10.5.

LEMMA 9.14. $(\mathcal{C}_\mathfrak{g}) \Rightarrow (\mathcal{C}_\Lambda)$.

*Proof.* Consider the diffeomorphic correspondence $F$ from a neighborhood of $(x_0, 0) \in \mathbb{R}^{n+k}$ to a neighborhood of $(x_0, x_0) \in \Lambda = \Lambda_\gamma$ defined by

$$\mathbb{R}^{n+k} \ni (\gamma^{-1}(x, t), t) \stackrel{F}{\mapsto} (\gamma^{-1}(x, t), x) \in \Lambda.$$

Thus $F(z, s) = (z, \gamma(z, s))$. Recall the vector fields $V_j$ in $\mathbb{R}^{n+k}$ associated to $\gamma^{-1}$, rather than to $\gamma$, by the relation

$$(\gamma^{-1}(x, t), t) = \exp(\sum_j t_j V_j)(x, 0),$$

as discussed in the proof of existence of the exponential representation.

In $\mathbb{R}^{n+k}$ there is a foliation by $k$-dimensional leaves parametrized by $t \mapsto (\gamma^{-1}(x,t), t)$, with a unique leaf passing through each point $(x, 0)$. The $C^\infty(\mathbb{R}^{n+k})$ module spanned by $\{V_j\}$ consists of all vector fields tangent to the leaves of this foliation. The correspondence $F$ maps each leaf to a set of the form $\{(\gamma^{-1}(x,t), x)\}$ where $x$ is constant. Therefore span $\{DF(V_j)\} \equiv \mathcal{V}_1$, as $C^\infty(\Lambda)$ modules.

On the other hand, denoting by $(z, s)$ coordinates in $\mathbb{R}^{n+k}$, we see that each vector field $DF(\partial/\partial s_i)$ in $\Lambda$ is tangent to the fibers of the projection $\pi_1$, since $F(z, s) = (z, \gamma(z, s))$. Thus span $\{DF(\partial/\partial s_i)\} \equiv \mathcal{V}_2$ at each point of $\Lambda$.

Combining these two equalities, we may identify $\mathfrak{G}(\Lambda)$ via $F^{-1}$ with the Lie algebra in $\mathbb{R}^{n+k}$ generated as a $C^\infty(\mathbb{R}^{n+k})$ module from all vector fields $V_j$ together with all $\partial/\partial s_i$. Thus $\mathfrak{g}(\Lambda)$ is identified with the $C^\infty$ module of vector fields on $\mathbb{R}^n$ obtained by restriction of the $\mathbb{R}^n$ component of each element of $\mathfrak{G}(\Lambda)$ to $s = 0$.[10] Each of the vector fields $X_\alpha$ constructed in the proof of existence of the exponential representation in subsection 9.1 was a finite linear combination of iterated commutators of the $V_j$ and $\partial/\partial s_i$, restricted to $s = 0$; that is to say, $X_\alpha$ belongs to $\mathfrak{g}(\Lambda)$. Consequently if $(\mathcal{C}_\Lambda)$ fails to hold, then $(\mathcal{C}_\mathfrak{g})$ likewise fails. □

For the purpose of performing explicit computations, it is sometimes convenient to describe $\Lambda$ implicitly by an equation $g = 0$, where $g$ maps $\mathbb{R}^n \times \mathbb{R}^n$ to $\mathbb{R}^{n-k}$ and is a submersion at $(x_0, x_0)$. As in Lemma 9.7, adopt coordinates $x = (x', x'') \in \mathbb{R}^k \times \mathbb{R}^{n-k}$ in which $x_0 = 0$ and $\{y : (0, y) \in \Lambda\} = \mathbb{R}^k \times \{0\}$. Then $g$ may be taken to have the form

$$(9.17) \qquad g(x, y) = y'' - x'' - h(x, y' - x')$$

where $h(x, 0) \equiv 0$, $h(0, y') \equiv 0$, and $\nabla h(0, 0) = 0$. A corresponding parametric representation of $\Lambda$ would be via $\gamma(x, t) = (x' + t, x'' + h(x, t))$. Vector fields in $\mathcal{V}_1, \mathcal{V}_2$ may then be regarded as restrictions to $\Lambda$ of those vector fields defined in some neighborhood of $(x_0, x_0)$ in $\mathbb{R}^n \times \mathbb{R}^n$ that annihilate $g$, and that also annihilate $\pi_2$ or $\pi_1$, respectively. Two such vector fields should be regarded as being equivalent if their restrictions to $\Lambda$ coincide.

We next analyze the meaning of $(\mathcal{C}_\Lambda)$ in more explicit terms in the hypersurface case $k = n - 1$. This is simpler than the case of higher codimension, because if there exists an invariant (to infinite order) submanifold at the origin, it must be $M_0$ (modulo an infinite order perturbation), whereas when $k < n-1$ any manifold of positive codimension containing $M_0$ is a potential candidate.

---

[10]This equals the module generated by differentiating each $V_j$ with respect to the coordinates $s$ an arbitrary number of times and by taking all possible Lie brackets, then restricting to $s = 0$.

PROPOSITION 9.15. *In the hypersurface case $k = n - 1$, $(\mathcal{C}_\Lambda)$ holds at $x_0$ if and only if $M_{x_0} \times M_{x_0} \subset \mathbb{R}^{n+n}$ is not tangent to $\Lambda$ to infinite order at $(x_0, x_0)$.*

Suppose that $\Lambda \subset \mathbb{R}^{n+n}$ satisfies our hypotheses and has codimension one. Choose coordinates $x = (x', x_n)$ in $\mathbb{R}^n$ with origin at $x_0$ such that $M_0 = \{y : (0, y) \in \Lambda\}$ coincides with the hyperplane $\{y : y_n = 0\}$. Then with the corresponding product coordinates in $\mathbb{R}^n \times \mathbb{R}^n$, the tangent space to $\Lambda$ at $(0, 0)$ is the subspace $\{(x, y) : x_n = y_n\}$. There exists a function $g : \mathbb{R}^{n+n} \mapsto \mathbb{R}$ of the form
$$g(x, y) = x_n - y_n - \phi(x, y')$$
such that $\Lambda = \{(x, y) : g(x, y) = 0\}$ near $(0, 0)$. Because $\{0\} \times M_0$ and the diagonal are contained in $\Lambda$, we have $\phi(0, 0) = 0$, and $\phi$ may be chosen to satisfy $\nabla \phi(0, 0) = 0$. The convention $\{0\} \times M_0 \subset \Lambda$ means that $g(0; y', 0) \equiv 0$; hence, equivalently, $\phi(0; y', 0) \equiv 0$.

We claim that in these coordinates, $(\mathcal{C}_\Lambda)$ holds at $0$ if and only if there exist multi-indices $\alpha, \beta$ such that

(9.18)
$$\frac{\partial^{\alpha+\beta}\phi}{\partial (x')^\alpha \partial (y')^\beta}(0, 0) \neq 0.$$

This is equivalent to the conclusion of Proposition 9.15.

*Proof.* For $i \leq n - 1$ define

(9.19)
$$\begin{aligned}
Y_i &= \frac{\partial}{\partial y_i} - \frac{\partial \phi}{\partial y_i}\frac{\partial}{\partial y_n}, \\
X_i &= \frac{\partial}{\partial x_i} + (1 - \partial \phi/\partial x_n)^{-1}\frac{\partial \phi}{\partial x_i}\frac{\partial}{\partial x_n}.
\end{aligned}$$

These vector fields are defined not merely on $\Lambda$, but on $\mathbb{R}^{n+n}$. They annihilate $g$, hence are tangent to its level sets, and in particular are tangent to $\Lambda$. Moreover $\mathcal{V}_1 = \text{span}\{X_i\}$ and $\mathcal{V}_2 = \text{span}\{Y_i\}$ as $C^\infty$ modules on $\Lambda$.

Define $f(x, y) = x_n$. Because $X_j(f)(0) = Y_j(f)(0)$ for all $j \leq n - 1$, $(\mathcal{C}_\Lambda)$ holds at $0$ if and only if there exists $V \in \mathfrak{G}(\Lambda)$ satisfying $V(f)(0) \neq 0$.

Fix a vector field $Z$ in $\mathbb{R}^{n+n}$ satisfying $Z(g) \equiv 0$ and $Z(f) \equiv 1$; then $Z$ does not belong to the span of $\mathcal{V}_1 + \mathcal{V}_2$ at $(0, 0)$, since every element of that span annihilates $f$ at the origin. Such a vector field exists, because the gradient of the restriction to $\Lambda$ of $f$ is nonzero at the origin. It is unique modulo multiplication by functions and addition of elements of the $C^\infty$ module $\mathcal{V}_1 + \mathcal{V}_2$. The vector fields annihilating $g$ are those in the module spanned by $\mathcal{V}_1 \cup \mathcal{V}_2 \cup \{Z\}$.

In order to compute $\mathfrak{G}(\Lambda)$, we express each of its elements as $V_1 + V_2 + \psi Z = V + \psi Z$, where $V_i \in \mathcal{V}_i$ and $\psi \in C^\infty$. Consider the $C^\infty(\Lambda)$ module $\mathcal{M}$ consisting of all functions $W(f)$ such that $W \in \mathfrak{G}(\Lambda)$, and let $\mathcal{M}^\dagger$ be the $C^\infty(\Lambda)$ module generated by all partial derivatives $\partial^{\alpha+\beta}\phi/\partial(x')^\alpha\partial(y')^\beta$. We claim that $\mathcal{M} \subset \mathcal{M}^\dagger$.

Note that $V \in \mathcal{V}_1 + \mathcal{V}_2$ implies $V(f) \in \mathcal{M}^\dagger$, by (9.19). Therefore $\mathcal{M} \subset \mathcal{M}^\dagger$ if and only if $\psi \in \mathcal{M}^\dagger$ for all vector fields $V + \psi Z \in \mathfrak{G}(\Lambda)$.

To prove the claim, it therefore suffices to show that $[V + \psi Z, V' + \psi' Z](f) \in \mathcal{M}^\dagger$ whenever $V, V' \in \mathcal{V}_1 + \mathcal{V}_2$ and $\psi, \psi' \in \mathcal{M}^\dagger$. That $[V, V'](f) \in \mathcal{M}^\dagger$ follows directly from (9.19). And $[\psi Z, V'](f) = \psi \cdot [Z, V'](f) + V'(\psi)$; the first term is a multiple of $\psi$ and hence belongs to $\mathcal{M}^\dagger$, while $V'(\psi) \in \mathcal{M}^\dagger$ for any $\psi \in \mathcal{M}^\dagger$ by (9.19).

The claim is therefore proved. Hence if $(\mathcal{C}_\Lambda)$ holds at the origin, then (9.18) must hold for some $\alpha, \beta$.

To prove the converse we do not claim that $\mathcal{M}^\dagger \subset \mathcal{M}$, but merely that for every $h \in \mathcal{M}^\dagger$ there exists $\psi \in \mathcal{M}$ such that $h - \psi$ is a finite $C^\infty$ linear combination of all $\partial^\delta \phi/\partial(y')^\delta$. Every such linear combination vanishes at the origin by our choice of coordinates, so the span of $\mathcal{M}$ at the origin (in the vector space $\mathbb{R}^1$) would thus contain the span of the set of all partial derivatives $\partial^{\alpha+\beta}\phi/\partial(x')^\alpha\partial(y')^\beta(0)$, which suffices to yield the converse.

Since $X_i \in \mathcal{V}_1 \subset \mathfrak{G}(\Lambda)$, $\partial \phi/\partial x_i \in \mathcal{M}$ for every $i \leq n-1$. Since $[V, \psi Z](f) = V(\psi)$ modulo a multiple of $\psi$ for any $V$, a simple induction shows that for any $(\alpha, \beta)$ with $\alpha \neq 0$, $\partial^{\alpha+\beta}\phi/\partial(x')^\alpha\partial(y')^\beta \in \mathcal{M}$, modulo addition of a $C^\infty$ multiple of the span of all $\partial^\delta\phi/\partial(y')^\delta$. □

*Remark* 9.2. The diagonal $\Delta \subset \Lambda$ is distinguished in this theory, because our primary aim is the analysis of operators of order zero, whose distribution kernels are carried by submanifolds, and whose restrictions to the submanifolds $M_x$ passing through $x$ are singular at $y = x$. But a pivotal part of the analysis, formalized in Theorem 8.11, concerns simpler variants more closely related to Radon transforms, whose distribution kernels have no additional singularities on the diagonal; for these, the diagonal should not be distinguished.

A related class of Radon type transforms can be associated to any $\Lambda \subset \mathbb{R}^{n_1+n_2}$, where $n_1$ need not equal $n_2$. We assume that both differentials $D\pi_i : T\Lambda \mapsto T\mathbb{R}^{n_i}$ are submersions, though this hypothesis might perhaps also be relaxed. To any measure $\mu$ on $\Lambda$ having a $C^\infty$ density, associate an operator $T$ by

$$\langle Tf_2, f_1 \rangle = \int_\Lambda f_2^*(y) f_1^*(x) \, d\mu(x,y),$$

where $f_i^* = f_i \circ \pi_i$. When $\Lambda$ is parametrized as $\{(x, \gamma(x,t))\}$, this is the class of all operators

$$Tf(x) = \int K(x,t) f(\gamma(x,t)) dt,$$

where $K \in C^\infty$.

Proposition 9.13 supplies a generalization of $(\mathcal{C}_\Lambda)$ to this situation. Indeed, the Lie algebras $\mathcal{V}_i$, and hence $\mathfrak{G}(\Lambda)$, may be defined as before. We say that $\Lambda$ satisfies $(\mathcal{C}_\Lambda)$ at a point $\lambda$ if $\mathfrak{G}(\Lambda)$ spans the tangent space to $\Lambda$ at $\lambda$. This generalization of $(\mathcal{C}_\Lambda)$ is invariant under the natural action of the product group Diff $(\mathbb{R}^{n_1}) \times$ Diff $(\mathbb{R}^{n_2})$.

Our methods lead to a generalization of Theorem 8.11: This curvature condition is necessary and sufficient for (all) nonsingular Radon transforms associated to $\Lambda$ to be smoothing of some positive order, in the general situation where $n_1, n_2$ are not necessarily equal.

## 10. Equivalence of curvature conditions

In this section we will demonstrate that $(\mathcal{C}_J) \Leftrightarrow (\mathcal{C}_\mathfrak{g}) \Leftrightarrow (\mathcal{C}_M)$. We will also introduce a variant $(\mathcal{C}_J)'$ and prove it to be equivalent to the other curvature conditions. Most implications are more easily demonstrated by proving their contrapositives.

10.1. *Invariant submanifolds and deficient Lie algebras.*

LEMMA 10.1.  $(\mathcal{C}_\mathfrak{g}) \Rightarrow (\mathcal{C}_M)$.

*Proof.* Suppose that $(\mathcal{C}_M)$ does not hold at $x_0$. Then there exists a smooth manifold $M$ of some dimension $0 \leq d < n$ containing $x_0$ that is invariant under $\gamma$ to infinite order at $x_0$. The condition $(\mathcal{C}_M)$ is diffeomorphism invariant, by Proposition 8.10, so we may switch to coordinates $y = (y', y'') \in \mathbb{R}^d \times \mathbb{R}^{n-d}$, with origin at $y_0 = 0$, so that $M = \{y : y'' = 0\}$ near 0. Write

$$\gamma(x', x'', t) = (\gamma'(x, t), \gamma''(x, t))$$

where $\gamma''(x, t) \in \mathbb{R}^{n-d}$.

By Lemma 9.7 and (9.15), invariance of $M$ to infinite order means that the Taylor expansion of $\gamma''(x', x'', t)$ with respect to $(x', x'', t)$ at $(0, 0, 0)$ takes the form

(10.1) $$\gamma''(x, t) \sim x'' + O(x'') \cdot O(t).$$

With respect to the coordinate system $(y', y''; s)$ in $\mathbb{R}^{n+k}$, any vector field $W$ in $\mathbb{R}^{n+k}$ may be decomposed uniquely as $W = W' + W''$ plus an element of the span of $\{\partial/\partial s_i\}$, where $W'$ belongs at every point to the span of $\{\partial/\partial y_i : i \leq d\}$ and $W''$ to the span of $\{\partial/\partial y_i : i > d\}$. There is a similar decomposition of any $W$ defined in $\mathbb{R}^n$.

Let $\mathcal{V}$ be the set of all vector fields $W$ defined near $0 \in \mathbb{R}^{n+k}$ that are everywhere in the span of $\{\partial/\partial y_j\}$, such that the Taylor expansion of $W''$ with respect to $(y, s)$ at the origin is $O(y'')$. Then $\mathcal{V}$ is closed under Lie brackets, and moreover, $[\partial_{s_i}, W] \in \mathcal{V}$ for every $W \in \mathcal{V}$ and every index $i$.

In the proof of Theorem 8.5 certain vector fields $Z_j$ were defined in (9.5) by $H_*(\partial/\partial t_j) = V_j = Z_j + \partial/\partial s_j$, where $H(x,t) = (\gamma(x,t), t)$. The hypothesis (10.1) means that $Z_j''(x,t) \sim O(x'')$.

The vector fields $Z_\beta$ defined in the proof of Theorem 8.5 were in the vector space spanned by iterated commutators of the $Z_j$ and of $\partial/\partial s_i$, and by the $Z_j$ themselves. Each $Z_j$ belongs to $\mathcal{V}$; hence so do all $Z_\beta$. The vector fields $X_\alpha$ defined in (9.7) by restricting $Z_\alpha$ to $s = 0$ therefore take the form

$$X_\alpha = X_\alpha' + X_\alpha'' \text{ where } X_\alpha''(x) \sim O(x'').$$

Once again, the set of all vector fields in $\mathbb{R}^n$ whose second components take this form is closed under Lie brackets. In particular, all iterated commutators of the $X_\alpha$ have this form, hence are tangent to $M$ to infinite order at $x = 0$. Since $M$ has positive codimension, they fail to span $T\mathbb{R}^n$ at 0. □

COROLLARY 10.2. *If a smooth manifold $M \ni y$ is invariant under $\gamma$ to infinite order at $y$, then $M$ is also invariant under $\gamma^{-1}$ to infinite order at $y$.*

*Proof.* Let the coordinate system be as in the preceding proof, so that $M = \{(x', x'') : x'' = 0\}$, and each vector field $X_\alpha = X_\alpha' + X_\alpha''$ in the exponential representation satisfies

$$X_\alpha''(x', x'') = O(x'') + O(|x'|^M) \text{ for every } M.$$

For any small $x', t$ and large $N$ let

$$\Phi(s) = \exp(-s \sum_{|\alpha| \leq N} t^\alpha X_\alpha/\alpha!)(x', 0).$$

The second component of $\Phi$ satisfies an ordinary differential equation

$$d\Phi''/ds = O(\Phi''(s)) + O(\Phi'(s))^M \text{ for every } M.$$

Moreover $d\Phi/ds = O(t)$, so that $\Phi(s) = O(st)$. Consequently the second component of $\Phi(1)$ is $O(|x'|^M + |t|^M)$ for every $M$ and $N$. Letting $N \to \infty$ we conclude that the second component of $\gamma_t^{-1}(x', 0)$ is $O(|t|^N + |x'|^M)$ for every $M$, as was to be proved. □

LEMMA 10.3. *Suppose that $\gamma$ is an analytic function of $x, t$ in a neighborhood of the origin in $\mathbb{R}^n \times \mathbb{R}^k$ and that curvature condition $(\mathcal{C}_J)$ does not hold at $x = 0$. Then there exists an analytic submanifold $M \subset \mathbb{R}^n$ containing 0, of positive codimension, such that $\gamma(x,t) \in M$ for all $(x,t) \in M \times \mathbb{R}^k$ sufficiently close to 0.*

*Proof.* Let $c_i$ be the maximum, as $\tau$ varies over a small neighborhood of 0, of the rank of the matrix $D\Gamma^i(0, \tau)$ of first partial derivatives with respect to $\tau \in \mathbb{R}^{ki}$ of $\Gamma^i$. If $c_1 = 0$ then $\Gamma^1(0, \tau) = \gamma_t(0) = 0$ for all sufficiently small $t$, and $M = \{0\}$ has the desired invariance property. Therefore we may suppose that $c_1 > 0$.

The hypothesis that curvature condition $(\mathcal{C}_J)$ fails to hold at $x = 0$ means that every term in the Taylor expansion with respect to $\tau$ of each Jacobian determinant $J_\xi(0, \tau)$ vanishes. Since $\gamma$ is analytic, each $J_\xi(0, \tau)$ vanishes identically. Therefore the rank of $D\Gamma^n(0, \tau)$ is strictly less than $n$ for every $\tau$, that is, $c_n < n$. Since $c_i \le c_{i+1}$ for all $i$, there exists a smallest integer $1 \le m < n$ such that $c_m = c_{m+1}$.

Let $\ell = c_m$. Fix $\tau_0$ so that $D\Gamma^m(0, \tau_0)$ has rank $\ell$, and consequently the rank is identically equal to $\ell$ at all $\tau$ near $\tau_0$. Therefore,[11] $\tau \mapsto \Gamma^m(0, \tau)$ traces out an $\ell$-dimensional real analytic submanifold $N$ of $\mathbb{R}^n$, as $\tau$ varies over a small neighborhood of $\tau_0$ in $\mathbb{R}^{km}$.

Consider the mapping $(s, \tau) \mapsto \gamma_s(\Gamma^m(0, \tau))$ for $s$ near 0. Since $c_m = c_{m+1}$, its matrix of first partial derivatives with respect to $\tau, s$ has rank identically equal to $\ell$, for $(s, \tau)$ near $(0, \tau_0)$. Therefore the image still traces out only an $\ell$-dimensional submanifold as $(s, \tau)$ varies over a small neighborhood of $(0, \tau_0)$. Hence $\gamma_s(y) \in N$ for all $s$ near 0 and $y$ near $\Gamma^m(0, \tau_0)$. Clearly then $N$ is also invariant under $\gamma^{-1}$ in the same sense.

This manifold $N$ need not contain the origin. To rectify the situation define $\Gamma^{-m}$ by $m$ iterations of $\gamma^{-1}$. Consider
$$G(\sigma, \tau) = \Gamma^{-m}(\Gamma^m(0, \tau), \sigma).$$
From the invariance of $N$ under $\gamma^{-1}$ it follows that the matrix of first partial derivatives of $G$ with respect to $\sigma, \tau$ has rank identically equal to $\ell$ in a neighborhood of $(0, \tau_0)$. Hence by analytic continuation its rank never exceeds $\ell$ for any $\sigma, \tau$. On the other hand, the rank of $\partial G/\partial \tau$ is $\ell$ for any $\sigma \in \mathbb{R}^{km}$ and for any $\tau$ sufficiently near $\tau_0$, since $y \mapsto \Gamma^{-m}(y, \sigma)$ is a diffeomorphism of a neighborhood of 0 with a neighborhood of $\Gamma^{-m}(y, \sigma)$, and the rank of the differential of the map $\tau \mapsto \Gamma^m(0, \tau)$ is already $\ell$.

Therefore by the inverse function theorem as above, the image under $G$ of a small neighborhood of $(\tau_0, \tau_0)$ in $\mathbb{R}^{kn+kn}$ traces out a submanifold $M$ of $\mathbb{R}^n$, of dimension $\ell < n$, and $0 = G(\tau_0, \tau_0)$ belongs to $M$. By the same reasoning as applied to the map $(\sigma, \tau) \mapsto G(\sigma, \tau)$ in the preceding paragraph, the rank of the differential of $(s, \sigma, \tau) \mapsto \gamma_s(G(\sigma, \tau))$ is everywhere $\le \ell$. Consequently $M$ is locally invariant under $\gamma$, by the same reasoning as applied to $N$ two paragraphs above. □

---

[11]If $F: M \mapsto N$ is a smooth mapping whose differential has constant rank $r$ in a neighborhood of $0 \in M$, and if $M$ has dimension $m$, construct a mapping $G = (F, g) : M \mapsto N \times \mathbb{R}^{m-r}$ whose differential is injective at 0. $G(M)$ is therefore a submanifold $Z$ of $N \times \mathbb{R}^{m-r}$ of dimension $m$. The implicit function theorem and constant rank hypothesis imply that the intersection of $Z$ with $N \times \{0\}$ is a submanifold of dimension $r$. This reasoning is also valid in the real analytic category since an implicit function theorem holds in that setting as well.

LEMMA 10.4.  $(\mathcal{C}_M) \Rightarrow (\mathcal{C}_{\mathfrak{g}})$.

*Proof.* Suppose that $\gamma$ fails to satisfy $(\mathcal{C}_{\mathfrak{g}})$ at $x_0$. Let $\{X_\alpha\}$ be the vector fields in the exponential representation for $\gamma$, let $\mathfrak{g}$ be the Lie algebra that they generate, and let $r$ be the dimension of the subspace of the tangent space spanned by $\mathfrak{g}$ at $x_0$. By hypothesis, $r < n$. If $r = 0$ then the manifold $\{x_0\}$ is invariant under $\gamma$ to infinite order, for

$$\exp\Big(\sum_{|\alpha|<N} t^\alpha X_\alpha/\alpha!\Big)(x_0) \equiv x_0$$

for every $N$, because each vector field $\sum_{|\alpha|<N} t^\alpha X_\alpha/\alpha!$ vanishes at $x_0$.

If $0 < r < n$ fix a maximal subset $\{V_j : 1 \le j \le r\}$ of $\mathfrak{g}$ that is linearly independent at $x_0$. For $v \in \mathbb{R}^r$ write $v \cdot \vec{V} = \sum_j v_j V_j$. The differential of the map $v \mapsto \exp(v \cdot \vec{V})(x_0)$ has full rank $r$ at $v = 0$; hence its image near $v = 0$ is a submanifold $M$ of $\mathbb{R}^n$ of dimension $r$. We claim that $M$ is invariant under $\gamma$ to infinite order at $x_0$.

To prove that $M$ is invariant to some arbitrarily high order $N$ write $y = (t, v) \in \mathbb{R}^{k+r}$ and consider the map $\phi$ defined by

$$\phi(y) = \gamma_t\Big(\exp(v \cdot \vec{V})(x_0)\Big) = \exp\Big(\sum_{0<|\alpha|<N} t^\alpha X_\alpha/\alpha!\Big)\exp(v \cdot \vec{V})(x_0) + O(|y|^N).$$

We claim that for some $\mathbb{R}^r$-valued polynomial $P$ depending on $N$,

(10.2) $$\phi(y) = \exp(P(y) \cdot \vec{V})(x_0) + O(|y|^N).$$

In this identity and similar ones below, the addition sign denotes the usual group operation in $\mathbb{R}^n$, in some fixed coordinate system.

By the Baker-Campbell-Hausdorff formula,

$$\phi(y) = \exp\Big(v \cdot \vec{V} + \sum_{1 \le |\beta| < N} y^\beta W_\beta^1\Big)(x_0) + O(|y|^N)$$

for certain vector fields $W_\beta^1$ belonging to $\mathfrak{g}$. Since every element of $\mathfrak{g}$ belongs to the span of $\{V_j\}$ at $x_0$, there exists for each $\beta$ a decomposition

$$W_\beta^1 = \sum_j c_{\beta,j} V_j + \tilde{W}_\beta^1$$

where each $c_{\beta,j}$ is a real constant, $\tilde{W}_\beta^1 \in \mathfrak{g}$, and $\tilde{W}_\beta^1(x_0) = 0$. Let $E_1 = E_1(y)$ be the vector field in $\mathbb{R}^n$, depending on the parameter $y$, defined by

$$E_1(y) = \sum_{1 \le |\beta| < N} y^\beta \tilde{W}_\beta^1.$$

Then

$$\phi(y) = \exp(P_1(y) \cdot \vec{V} + E_1(y))(x_0) + O(|y|^N)$$

for some vector valued polynomial $P_1$ satisfying $P_1(0) = 0$.

Since $E_1(y)(x_0) = 0$ for every $y$,
$$\exp(-E_1(y))(x_0) \equiv x_0.$$
Substituting this into the preceding identity and invoking the Baker-Campbell-Hausdorff formula give
$$\begin{aligned}\phi(y) &= \exp(P_1(y) \cdot \vec{V} + E_1(y))\exp(-E_1(y))(x_0) + O(|y|^N) \\ &= \exp(P_1(y) \cdot \vec{V} + \sum_{2 \le |\beta| < N} y^\beta W_\beta^2)(x_0) + O(|y|^N)\end{aligned}$$
where each $W_\beta^2 \in \mathfrak{g}$. Only multi-indices satisfying $|\beta| \ge 2$ arise, because $P_1(y)$ and $E_1(y)$ both vanish when $y = 0$. As above, this may be rewritten as
$$\phi(y) = \exp(P_2(y) \cdot \vec{V} + E_2(y))(x_0) + O(|y|^N),$$
where again $E_2(y)(x_0) \equiv 0$ and $E_2(y) \in \mathfrak{g}$ for every $y$, and now $E_2(y) = O(|y|^2)$ near $y = 0$. Again $\exp(-E_2(y))(x_0) \equiv x_0$. Substituting this and applying the Baker-Campbell-Hausdorff formula once more leads to a remainder term $E_3(y)$ which is $O(|y|^3)$ for small $y$, since $E_2$ is $O(|y|^2)$ and $P_2(y)$ is $O(|y|)$.

Repeating this argument $N$ times yields
$$\phi(y) = \exp(P_N(y) \cdot \vec{V} + E_N(y))(x_0) + O(|y|^N),$$
where $E_N(y) = O(|y|^N)$ and $P_N$ is a polynomial. Deleting the term $E_N$ from the right-hand side now merely results in an additional error term $O(|y|^N)$; so (10.2) is proved.

Any $x \in M$ close to $x_0$ may be expressed as $\exp(v \cdot \vec{V})(x_0)$, and $|v| \sim |x - x_0|$. Consequently
$$\gamma_t(x) = \phi(t, v) + O(|t|^N).$$
Since $\exp(P_N(y) \cdot \vec{V})(x_0) \in M$, the distance from $\gamma_t(x)$ to $M$ is $O(|t| + |v|)^N$, as was to be proved. □

The following lemma is needed to complete the proof of Proposition 9.13.

LEMMA 10.5. $(\mathcal{C}_\Lambda) \Rightarrow (\mathcal{C}_M)$.

*Proof.* Suppose that $(\mathcal{C}_M)$ fails to hold at $x_0$, so that there exists a submanifold $M$ of positive codimension invariant under $\gamma$ to infinite order at $x_0$. Choose coordinates $(x', x'')$ as in Lemma 9.7, (9.14) and (9.15), in which $x = 0$ and $M = \{x'' = 0\}$.

When we let $H(x,t) = (\gamma(x,t),t)$ and employ the notation of (9.14), $\gamma''(x,t) \sim O(x'')$ and consequently every $Z_j''(x,t) \sim O(x'')$. Therefore the $x''$ component of any partial derivative of any $Z_j$ with respect to $t$, evaluated at $t = 0$, is likewise $\sim O(x'')$. Just as in the proof that $(\mathcal{C}_\mathfrak{g})$ implies $(\mathcal{C}_M)$, this means that the span of the Lie algebra generated by all these partial derivatives is contained in the tangent space to $M$ at $0$. We have seen in the course of the proof of Lemma 9.14 that $\mathfrak{g}(\Lambda)$ is the $C^\infty$ module generated by all such partial derivatives and by forming iterated Lie brackets. Therefore $\mathfrak{g}(\Lambda)$ does not span, so $(\mathcal{C}_\Lambda)$ fails to hold. □

10.2. *Vanishing Jacobians.* In order to prove that $(\mathcal{C}_J)$ implies $(\mathcal{C}_\mathfrak{g})$ we require two subsidiary lemmas. Fix any coordinate system in a neighborhood of $x_0$. Then given any $n$ vector fields $V_1, \ldots V_n$ we may form the determinant $\det(V_1, \ldots V_n)(x)$ by writing $V_i(x)$ as $\sum a_{ij}(x)\partial/\partial x_j$ and taking the determinant of the $n \times n$ matrix $\big(a_{ij}(x)\big)$. Denote by $\nabla \cdot V$ the divergence of $V$.

By a polynomial in $\{X_\alpha\}$ we will mean any vector field expressible as a finite linear combination, with constant coefficients, of iterated commutators of the $X_\alpha$. By a Taylor series in $\{X_\alpha\}$ we will mean any formal series $V(x,\tau) = \sum_{|\sigma|>0} \tau^\sigma U_\sigma(x)$, where each $U_\sigma$ is a polynomial in $\{X_\alpha\}$ and $\tau = (\tau_{ij}) \in \mathbb{R}^{kn}$.

LEMMA 10.6.  *There exist Taylor series $W_{ij}$ in $\{X_\alpha\}$ such that*
$$\frac{\partial \Gamma(x,\tau)}{\partial \tau_{ij}} \sim W_{ij}(\Gamma(x,\tau))$$
*as Taylor series in $\tau$.*

*Proof.* By the Baker-Campbell-Hausdorff formula,
$$\Gamma(x,\tau) \sim \exp\Big(\sum_\sigma \tau^\sigma U_\sigma/\sigma!\Big)(x)$$
for certain polynomials $U_\sigma$ in $\{X_\alpha\}$. So
$$\Gamma(x,\tau+\Delta\tau) \sim \exp\Big(\sum \frac{\tau^\sigma}{\sigma!}U_\sigma + \Delta\tau \cdot W\Big)(x) + O(\Delta\tau)^2$$
as Taylor series in $(\tau, \Delta\tau)$, for a certain $kn$-tuple $W$ of Taylor series in $\{X_\alpha\}$. Suppose that $\Delta\tau$ has $i,j^{\text{th}}$ entry equal to $s$ and all other entries zero. Then by another application of the Baker-Campbell-Hausdorff formula we may rewrite
$$\begin{aligned}\Gamma(x,\tau+\Delta\tau) &\sim \exp(sW_{ij})\exp(\sum \tau^\sigma U_\sigma/\sigma!)(x) + O(s)^2 \\ &\sim \exp(sW_{ij})(\Gamma(x,\tau)) + O(s)^2\end{aligned}$$
where $W_{ij}$ is a Taylor series in $\{U_\sigma\}$, and hence is a Taylor series in $\{X_\alpha\}$. Differentiating with respect to $s$ and evaluating at $s=0$ conclude the proof. □

LEMMA 10.7.  *As a Taylor series in $\tau$, each $J_\xi(x,\tau)$ has the form*
$$\sum_{|\beta|>0} \tau^\beta V_\beta(x)$$

*where each function $V_\beta$ is a finite linear combination of terms*

$$\varphi(x) \cdot \det(Z_1(x), \ldots Z_n(x)),$$

*each $\varphi$ belongs to $C^\infty$, and the vector fields $Z_j$ are iterated commutators of the $X_\alpha$.*

*Proof.* By its definition and by the preceding lemma,

$$J_\xi(x,\tau) = \det\Big(U_1(\Gamma(x,\tau),\tau), \ldots U_n(\Gamma(x,\tau),\tau)\Big)$$

where each $U_k$ is one of the $W_{ij}$, depending on $\xi$. Our aim is to express each partial derivative of $J_\xi$ with respect to $\tau$ as a finite linear combination of $C^\infty$ functions of $x$ times determinants of similar $n$-tuples, each component of which is another Taylor series in $\{X_\alpha\}$, evaluated at $(\Gamma(x,\tau),\tau)$. Setting $\tau = 0$ then yields the result desired.

Set

$$F(x,\tau) = \det\Big(U_1(x,\tau), \ldots U_n(x,\tau)\Big)$$

so that $J_\xi(x,\tau) = F(\Gamma(x,\tau),\tau)$. When $\partial J_\xi(x,\tau)/\partial \tau_{ij}$ is calculated by the chain rule, the derivative falls either on $\Gamma(x,\tau)$, or on the second occurrence of $\tau$ in $F(\Gamma(x,\tau),\tau)$. In the latter event what results is equal, by the multilinearity of the determinant, to

$$\det\left(\frac{\partial U_1}{\partial \tau_{ij}}(\Gamma(x,\tau),\tau), U_2(\Gamma(x,\tau),\tau), \ldots U_n(\Gamma(x,\tau),\tau)\right)$$

plus $n-1$ more terms of the same form, in each of which some other $U_k$ is differentiated instead of $U_1$.

In the former event we may invoke the preceding lemma to express the result as $(W_{ij}F)(\Gamma(x,\tau),\tau)$. It is shown in [29, Lemma 2.6], that for any vector fields $B, A_1, \ldots A_n$,

$$B(\det(A_1, \ldots A_n))$$
$$= (\nabla \cdot B) \cdot \det(A_1, \ldots A_n) + \sum_{j \le n} \det(A_1, \ldots A_{j-1}, [B, A_j], A_{j+1}, \ldots A_n).$$

Thus $\partial J_\xi / \partial \tau_{ij}$ equals $(\nabla \cdot W_{ij}) J_\xi$ plus $2n$ terms of the general form

$$\det\Big(A_1(\Gamma(x,\tau),\tau), \ldots A_n(\Gamma(x,\tau),\tau)\Big)$$

where each $A_i(x,\tau)$ is a Taylor series in $\{X_\alpha\}$. Each such term is of the same form as $J_\xi$ itself. Therefore combining the two cases that resulted from application of the chain rule we find, by induction on $|\beta|$, that every partial derivative $\partial^\beta J_\xi / \partial \tau^\beta$ is of the desired form. □

LEMMA 10.8. $(\mathcal{C}_J) \Rightarrow (\mathcal{C}_{\mathfrak{g}})$.

*Proof.* Suppose that $(\mathcal{C}_\mathfrak{g})$ does not hold at $x$; that is, the $X_\alpha$ and all their commutators fail to span the tangent space to $\mathbb{R}^n$ at $x$. Then each $J_\xi(x,\tau)$ vanishes to infinite order at $\tau = 0$ by the preceding lemma, so that $(\mathcal{C}_J)$ also fails to hold at $x$. □

It is comparatively easy to show that $(\mathcal{C}_\mathfrak{g})$ implies a weaker form $(\mathcal{C}_J)'$ of $(\mathcal{C}_J)$, which suffices for our analytic applications. Recall that $\Gamma^r$ denotes the $r$-fold iterate of $\gamma$. Now $(\mathcal{C}_J)$ is the special case $r=n$ of the following definition.

*Definition* 10.9.   $\gamma$ satisfies $(\mathcal{C}_J)'$ at $x_0$ if there exist an integer $r$ and an $n$-tuple $\xi \in \{1,2,\ldots rk\}$ such that the Jacobian determinant of the $n$ by $n$ matrix of first partial derivatives of $\Gamma^r(x_0,\tau)$ with respect to $(\tau_{\xi_1},\ldots \tau_{\xi_n})$ does not vanish to infinite order at $\tau = 0$.

Note that $(\mathcal{C}_J)' \Rightarrow (\mathcal{C}_\mathfrak{g})$, by the proof of Lemma 10.8.

LEMMA 10.10.   $(\mathcal{C}_\mathfrak{g}) \Rightarrow (\mathcal{C}_J)'$.

*Proof.* This time we do not prove the contrapositive; suppose that $\gamma$ satisfies $(\mathcal{C}_\mathfrak{g})$ at $x_0 \in \mathbb{R}^n$. Let $S$ be the set of all $k$-tuples of nonnegative integers. For $\alpha = (\alpha^1,\ldots \alpha^k) \in S$, set $|\alpha| = \sum \alpha^i$. Let $S_m$ be the set of all $\alpha \in S$ satisfying $|\alpha| \le m$.

Denote by
$$X_I = [X_{\alpha_1}, [X_{\alpha_2}, \cdots [X_{\alpha_{q-1}}, X_{\alpha_q}]\cdots]$$
an arbitrary iterated commutator of $\{X_\alpha\}$, where each $\alpha_i \in S$, and $I = (\alpha_1,\ldots \alpha_q) \in S^q$ for arbitrary $q \ge 1$. Define $|I| = \sum_1^q |\alpha_i|$. Thus each $X_\alpha$ is an $X_I$, but not vice versa.

Fix an integer $m$ such that the set of all $X_I$ satisfying $|I| \le m$ spans $T\mathbb{R}^n$ at $x_0$. Let $\mathcal{N} = \mathcal{N}_m^{a_1,\cdots}$ be the free graded nilpotent Lie algebra introduced in Section 2, with one generator $Y_\alpha$ for each $\alpha \in S_m$. Let $\{Y_I : I \in \mathcal{I}\}$ be the basis for $\mathcal{N}$ discussed in Section 2, where $\mathcal{I}$ denotes the set of all basic multi-indices $I$ satisfying $|I| \le m$. Throughout this proof $I$ will always denote an element of $\mathcal{I}$. Denote by $N$ the connected, simply connected Lie group whose Lie algebra is $\mathcal{N}$. Identify the $Y_I$ with left invariant vector fields on $N$, via the exponential map $\sum u_I Y_I \mapsto \exp(\sum u_I Y_I)(0)$. Let $d$ be the dimension of $\mathcal{N}$.

To each $Y_I$ associate a vector field $X_I$ in a fixed neighborhood of $x_0 \in \mathbb{R}^n$, by replacing each $Y_{\alpha_j}$ in the commutator representation for $Y_I$ by the corresponding vector field $X_{\alpha_j}$. Then $\{X_I : I \in \mathcal{I}\}$ spans $T\mathbb{R}^n$ at $x_0$, by hypothesis and the choice of $m$.

Define
$$\tilde{\gamma}(x,t) = \exp\Big(\sum_{0<|\alpha|\le m} t^\alpha Y_\alpha/\alpha!\Big)(x) \text{ for } (x,t) \in N \times \mathbb{R}^k.$$

Denote by $\tilde{\Gamma}^\rho(x,\tau)$ its $\rho$-fold iterate, $\tau \in \mathbb{R}^{\rho k}$.

At the origin $\tilde{\gamma}$ satisfies $(\mathcal{C}_J)$. For with the coordinates on $N$ defined by the exponential map and the basis $\{Y_I\}$ for $\mathcal{N}$, $\tilde{\gamma}$ is a polynomial in $(x,t)$. Hence by Lemma 10.3, if $(\mathcal{C}_J)$ were to fail to hold then there would exist a manifold $M \ni 0$ of positive codimension invariant under $\tilde{\gamma}$. Each of the vector fields $Y_\alpha$, $\alpha \in S_m$, would be tangent to $M$ at every point of $M$, as follows from an induction on $|\alpha|$ and the proof of uniqueness of the exponential representation. Hence so would be all their commutators. Therefore $\{Y_\alpha\}$ would fail to generate $\mathcal{N}$, a contradiction.

Since $(\mathcal{C}_J)$ holds, we may fix $r$ and $\xi \in \{1, 2, \ldots rk\}^d$ such that the Jacobian determinant $J_\xi(\tau)$ of the $d$ by $d$ matrix of first partial derivatives of $\tilde{\Gamma}^r(0,\tau)$ with respect to $\tau_{\xi_1}, \ldots \tau_{\xi_d}$ does not vanish to infinite order as a function of $\tau \in \mathbb{R}^{rk}$ at $\tau = 0$. By the Baker-Campbell-Hausdorff formula, $\tilde{\Gamma}^r(0,\tau)$ takes the form $\exp(\sum_I P_I(\tau)Y_I)(0)$, where each $P_I$ is a homogeneous polynomial of degree $|I| \leq m$.

Denote by $\Gamma^r$ the corresponding $r$-fold iterate of $\gamma$, in a neighborhood of $x_0$ in $\mathbb{R}^n$. By the iteration process in the proof of Lemma 10.4, since $\{X_I : I \in \mathcal{I}\}$ spans $T\mathbb{R}^n$ at $x_0$, the Taylor expansion of $\Gamma^r(0,\tau)$ at $\tau = 0$ may be expressed as

$$\exp(\sum_{I \in \mathcal{I}} (P_I + R_I)(\tau)X_I)(x_0),$$

as a formal Taylor series, where the $P_I$ are the same polynomials as above, and each formal series $R_I$ vanishes at least to order $m+1$ at $\tau = 0$. This representation is of course not unique, unless $\{X_I : I \in \mathcal{I}\}$ is linearly independent at $x_0$.

Returning to $N$, we consider the formal Taylor series defined by

$$\hat{\Gamma}^r(\tau) \sim \exp(\sum_{I \in \mathcal{I}} (P_I + R_I)(\tau)Y_I)(0).$$

We claim that the Jacobian determinant $\hat{J}_\xi$ for $\hat{\Gamma}^r$, which is a formal Taylor series in $\tau$ rather than a function, does not vanish to infinite order at $\tau = 0$. Indeed, $J_\xi$ is a homogeneous polynomial in $\tau$ of some degree $Q$, and each $R_I$ vanishes at least to order $m+1$ whereas each $P_I$ is homogeneous of some degree $\leq m$. Therefore $\hat{J}_\xi - J_\xi$ vanishes at least to order $Q+1$ at $\tau = 0$, from which the claim follows.

Consider the map $F : N \mapsto \mathbb{R}^n$ defined by

$$F : \exp(\sum_I u_I Y_I)(0) \mapsto \exp(\sum_I u_I X_I)(x_0).$$

By their definitions,

(10.3) $$\Gamma^r(x_0, \tau) \sim F \circ \hat{\Gamma}^r(0,\tau)$$

as formal Taylor series. Note that $F$ is a well defined $C^\infty$ map, taking $0$ to $x_0$, which is a submersion at $0$ because $\{X_I : I \in \mathcal{I}\}$ spans $T\mathbb{R}^n$ at $x_0$.

Therefore by the implicit function theorem there exists a diffeomorphism $\Psi$ of a neighborhood of $0 \in N$ with a neighborhood of $0 \in \mathbb{R}^d$, mapping 0 to 0, so that $F \equiv \pi \circ \Psi$ where $\pi : \mathbb{R}^d \mapsto \mathbb{R}^n$ is the projection $\pi(x', x'') = x'$ for $(x', x'') \in \mathbb{R}^n \times \mathbb{R}^{d-n} = \mathbb{R}^d$. By (10.3),

(10.4) $$\Gamma^r(x_0, \tau) \sim \pi \circ \check{\Gamma}^r$$

where $\check{\Gamma}^r(0, \tau) = \Psi(\hat{\Gamma}^r(0, \tau))$.

Consider the corresponding Jacobian determinant $\check{J}_\xi(\tau)$ of $\check{\Gamma}^r$ with respect to $\tau_{\xi_1}, \ldots \tau_{\xi_d}$; this is another Taylor series in $\tau$. Formally it equals $\hat{J}_\xi$ times a nonvanishing function of $\tau$, so it does not vanish to infinite order at $\tau = 0$. By elementary matrix operations, the $d$ by $d$ matrix determinant $\check{J}_\xi(\tau)$ may be expanded as a finite linear combination of determinants of its $n$ by $n$ minors of the form

$$J^{(\gamma)}(\tau) = \partial[\pi \circ \check{\Gamma}^r]/\partial[\tau_{\xi_{\gamma_1}}, \ldots \tau_{\xi_{\gamma_n}}],$$

where each $J^{(\gamma)}$ and each coefficient in this linear combination is itself a formal Taylor series in $\tau$, and $\{\xi_{\gamma_1}, \ldots \xi_{\gamma_n}\}$ runs over all subsets of $\{\xi_1, \ldots \xi_d\}$ of cardinality $n$.

Since $\check{J}_\xi$ does not vanish to infinite order at $\tau = 0$, there must exist $\gamma$ such that $J^{(\gamma)}$ likewise does not vanish to infinite order. But by (10.4), $J^{(\gamma)}(\tau)$ is the Taylor series of the Jacobian determinant of $\Gamma^r(x_0, \tau)$ with respect to $\{\tau_{\xi_i} : i = \gamma_1, \ldots \gamma_n\}$. Hence the original $\gamma$ in $\mathbb{R}^n$ satisfies $(\mathcal{C}_J)'$ at $x_0$. $\square$

Both halves of the following result have now been established.

COROLLARY 10.11.  $(\mathcal{C}_J)' \Leftrightarrow (\mathcal{C}_\mathfrak{g})$.

*Remark* 10.1. We have shown that $(\mathcal{C}_\mathfrak{g})$ implies a more precise form of $(\mathcal{C}_J)'$, in which appears the iterate $\Gamma^r$ of $\gamma$ of a specific order $r = d$, equal to the dimension of a certain free nilpotent Lie algebra determined by the Lie algebra in $\mathbb{R}^n$ generated by $\{X_\alpha\}$.

10.3. *Construction of invariant submanifolds.* In order to complete the proof of equivalence of $(\mathcal{C}_\mathcal{J})$ with the other curvature conditions, we need the following implication.

LEMMA 10.12.  $(\mathcal{C}_M) \Rightarrow (\mathcal{C}_J)$.

This lemma and its lengthy proof are not needed for our analytic applications, so this section may safely be skipped over by the reader primarily interested in the analysis of operators. To prove it we will show how submanifolds invariant under $\gamma$ to infinite order may be constructed, given that $(\mathcal{C}_J)$ fails to hold.

In the next two lemmas we reduce the general $C^\infty$ case to that of polynomial mappings $\gamma$, so that Lemma 10.3 on the analytic case may be applied.

Some care is required in doing so, for the hypothesis that $(\mathcal{C}_J)$ fails to hold at $x_0$ is not stable under truncation of the Taylor expansion for $\gamma(x,t)$ at $(x_0, 0)$, even though the vanishing of any particular Taylor coefficient of any Jacobian determinant $J_\xi$ is stable under truncations of sufficiently high order.

Let there be given an integer $n \geq 1$ and weights $0 < w_1 \leq \ldots w_n < \infty$. For any $m \leq n$ define mappings $\delta_r : \mathbb{R}^m \times \mathbb{R}^k \mapsto \mathbb{R}^m \times \mathbb{R}^k$ by

$$\delta_r(x_1, \ldots x_m; t) = (r^{w_1} x_1, \ldots r^{w_m} x_m; rt).$$

The same notation $\delta_r$ will be used to denote the related mappings

$$\delta_r(x_1, \ldots x_m) = (r^{w_1} x_1, \ldots r^{w_m} x_m).$$

We say that a function $P : \mathbb{R}^{m+k} \mapsto \mathbb{R}$ is homogeneous of weight $v$ if

$$P(\delta_r(x,t)) \equiv r^v P(x,t) \qquad \text{for all } r \in \mathbb{R}^+ \text{ and all } (x,t) \in \mathbb{R}^{m+k}.$$

Any such homogeneous function that is $C^\infty$ is necessarily a polynomial. The analogous definition applies to functions with domain $\mathbb{R}^m$. More generally, some weights will be allowed to equal $+\infty$. If $w_j < \infty$ for all $j \leq q \leq n$, then the dilations $\delta_r$ are still defined as above on $\mathbb{R}^m \times \mathbb{R}^k$ and on $\mathbb{R}^m$ for each $m \leq q$, and homogeneity of weight $v$ is still defined for functions with domains $\mathbb{R}^{m+k}$ or $\mathbb{R}^m$, for all $m \leq q$.

The strategy of our construction of invariant submanifolds is to simultaneously construct an associated system of coordinates and weights[12] in which $\gamma$ is expressed as an irreducible leading order part, plus a less important higher order perturbation. The weights serve as measures of the extent to which $\gamma$ looks flat, that is, not curved, in a given coordinate system. Larger weights correspond to less apparent curvature, and to a higher degree of invariance for a subspace $x'' = 0$, in some coordinate system split as $x = (x', x'')$.

As an illustration consider the example

$$\gamma(x,t) = (x_1 + t, x_2 + 2x_1 t + t^2),$$

to which weights $(w_1, w_2) = (1, 2)$ would naturally be assigned, because $\gamma(rx_1, r^2 x_2, rt) = \delta_r \gamma(x, t)$ where $\delta_r(x) = (r^{w_1} x_1, r^{w_2} x_2)$. This family of curves may be transformed by a change of coordinates to $\hat\gamma(x,t) = (x_1 + t, x_2)$, for which the weights are $(1, \infty)$ and the line $x_2 = 0$ is an invariant submanifold. Now $\hat\gamma$ is homogeneous with respect to any weights $(1, w_2)$, and in such a case we choose the largest possible $w_2$.

A second example is

$$\gamma(x,t) = (x_1 + t, x_2 + 2x_1 t + t^2 + t^3).$$

---

[12]This construction was originally devised as the basis for a lifting procedure which dealt directly with a mapping $\gamma$, rather than with the vector fields $X_\alpha$ in its exponential representation.

This mapping is not homogeneous with respect to any system of weights, but the system of weights $(1,2)$ is naturally associated to $\gamma$, for the dilation group defined by these weights satisfies

$$\gamma(\delta_r(x,t)) = \delta_r \gamma(x,t) + \left(0, O(r^3)\right),$$

and the exponent 3 in the remainder term is greater than the weight $w_2 = 2$. Thus $\gamma$ is homogeneous modulo a higher order remainder term. This mapping $\gamma$ may be transformed by conjugation with an appropriate diffeomorphism into $(x_1+t, x_2+t^3)$, which is homogeneous with respect to the weights $(1,3)$. After conjugation no weight has decreased, and one weight has increased.

In the next lemma $t$ will as usual be in $\mathbb{R}^k$, but $x \in \mathbb{R}^{p+1}$ for some $p \geq 1$, and we write $x = (x'; x_{p+1}) \in \mathbb{R}^p \times \mathbb{R}$. The map $\gamma$ is assumed to send $\mathbb{R}^{p+1} \times \mathbb{R}^k$ to $\mathbb{R}^{p+1}$.

LEMMA 10.13. *Suppose that in $\mathbb{R}^{p+1} = \mathbb{R}^p \times \mathbb{R}^1$, $\gamma$ takes the special form*

$$\gamma_t(x) = (\gamma'_t(x'); x_{p+1} + P_{p+1}(x',t))$$

*where for each $j$,*

$$\begin{aligned}
\gamma'_t(x') &= (x_1 + P_1(x',t), \ldots x_p + P_p(x',t)) \in \mathbb{R}^p, \\
P_j(x',0) &\equiv 0 \quad \text{for all } j, \\
P_j &\text{ depends only on } x_1, \ldots x_{j-1}, t, \\
1 &\leq w_1, \ldots w_{p+1} \in \mathbb{Z}, \\
P_j &\text{ is a homogeneous polynomial of weight } w_j,
\end{aligned}$$

*and where*

(10.5) $\quad\quad\quad\quad \gamma'$ satisfies $(\mathcal{C}_J)$ at $0 \in \mathbb{R}^p$,

$\quad\quad\quad\quad \gamma$ does not satisfy $(\mathcal{C}_J)$ at $0 \in \mathbb{R}^{p+1}$.

*Then there exists a homogeneous polynomial $h : \mathbb{R}^p \mapsto \mathbb{R}$ of weight $w_{p+1}$, such that*

(10.6) $\quad\quad\quad\quad P_{p+1}(x',t) \equiv \left(h \circ \gamma'_t - h\right)(x').$

*Moreover, if a diffeomorphism $\Phi$ of $\mathbb{R}^{p+1}$ is defined to be $\Phi(x) = (x'; x_{p+1} - h(x'))$, then the conjugated family of mappings $\tilde{\gamma}_t(x) = \Phi \gamma_t(\Phi^{-1} x)$ takes the form*

$$\tilde{\gamma}_t(x) \equiv (\gamma'_t(x'); x_{p+1}).$$

In the second to last hypothesis (10.5), $\gamma'$ is regarded as a mapping from $\mathbb{R}^{p+k}$ to $\mathbb{R}^p$, and is assumed to satisfy $(\mathcal{C}_J)$ as such a mapping.

*Proof.* Let $\Gamma'^j(x',\tau)$ and $\Gamma^j(x,\tau)$ denote the $j^{\text{th}}$ iterates of $\gamma'$ and $\gamma$, respectively, and let $\pi : \mathbb{R}^{p+1} \mapsto \mathbb{R}^p$ be the projection onto the first $p$ coordinates.

Since $\pi \circ \gamma_t \equiv \gamma'_t \circ \pi$, $\pi \Gamma^j(x,\tau) \equiv \Gamma'^j(\pi x, \tau)$. Now apply the construction of Lemma 10.3 to the real analytic mapping $\gamma$, and let $m, \ell$ be the integers defined there.

We assert that $\ell = p$. Indeed, some iterate $\Gamma'^j(0,\tau)$ has a surjective differential at some $\tau_0$, since $\gamma'$ satisfies $(\mathcal{C}_J)$. The rank of the differential of $\Gamma^j(0,\tau)$ must be at least as large as that of $\pi \Gamma^j(0,\tau) = \Gamma'^j(0,\tau)$, so is at least $p$ at $\tau_0$. On the other hand it never exceeds $p$ since $(\mathcal{C}_J)$ is assumed not to hold for $\gamma$ at $x = 0$. Thus in the reasoning of Lemma 10.3 we have $c_p = c_{p+1} = p$, and the differential of $\Gamma'^p(0,\tau)$ is surjective at $\tau_0$.

Now in a small neighborhood of $(\tau_0, \tau_0)$, the differential of
$$(\sigma, \tau) \mapsto \Gamma'^{-k}(\Gamma'^k(0,\tau), \sigma)$$
has rank $p$. In other words, the map $\pi|_M$ from the manifold $M$ constructed in Lemma 10.3 to $\mathbb{R}^p$ has a surjective differential in a neighborhood of the origin. Therefore $M$ may be expressed locally as the graph of a smooth function:
$$M = \{(x', h(x')) : |x'| < \varepsilon\}.$$

Define $\delta_r x = (r^{w_1} x_1, \ldots r^{w_{p+1}} x_{p+1})$. Then $\gamma_{rt}(\delta_r x) \equiv \delta_r \gamma_t(x)$ by the homogeneity hypothesis. The same identity holds for $\gamma^{-1}$, by uniqueness of the inverse. By iteration the same follows for $\Gamma^{-p} \circ \Gamma^p$. Therefore the manifold $M$ is invariant under $\delta$ in a neighborhood of the origin, so $h$ must satisfy
$$h(r^{w_1} x_1, \ldots r^{w_p} x_p) \equiv r^{w_{p+1}} h(x_1, \ldots x_p).$$
Then $h$ must be a polynomial.

Now we may define $M$ globally as the graph of $h$ over all of $\mathbb{R}^p$. By analytic continuation, $\gamma_t(M) \subset M$ for all $t \in \mathbb{R}^k$. In terms of the equation for $M$, this says that for all $x', t$,
$$\gamma_t(x'; h(x')) = \Big(\gamma'_t(x'); h(\gamma'_t(x'))\Big).$$
But the left-hand side also equals
$$\Big(\gamma'_t(x'); h(x') + P_{p+1}(x', t)\Big),$$
so $h$ satisfies the cocycle identity (10.6). The second conclusion of the lemma now follows directly from the definition of $\Phi$. □

In the next lemma we work in $\mathbb{R}^n$ with a fixed coordinate system and an associated ordered $n$-tuple $(w_1, \ldots w_n)$ of weights, for which each $w_j$ is either a positive integer, or equals $+\infty$. Each component of $t$ is always assigned weight 1. With respect to such a system of weights we make the following definition.

*Definition* 10.14. The weight of a monomial $t^\beta \prod x_j^{\alpha_j}$ with respect to the system of weights $(w_1, \ldots w_n)$ is defined to be
$$(10.7) \qquad |\beta| + \sum \alpha_j w_j .$$

The weight of a $C^\infty$ real-valued function $F$ of $(x,t)$ equals the minimum of the weights of all monomials $x^\alpha t^\beta$ such that $\partial^{\alpha+\beta} F/\partial x^\alpha \partial t^\beta(0,0) \neq 0$, provided there exists at least one such $(\alpha, \beta)$. It equals $+\infty$ if the Taylor series of $F$ with respect to $(x,t)$ about $(0,0)$ vanishes identically.

Observe that the sum (10.7) may equal either a nonnegative integer, or $+\infty$. If some weight $w_j$ equals $+\infty$, then $F$ may have weight $+\infty$, even if its Taylor expansion does not vanish identically.

LEMMA 10.15. *Let an integer $q < n$ and weights $(w_1, \ldots w_n)$ be given. Suppose that $w_j < \infty$ for all $j \leq q$, and $w_j = +\infty$ for all $j > q$. Assume that $\gamma$ takes the form*

$$(10.8) \qquad \gamma_t(x) = x + \big(Q_1(x,t), \ldots Q_n(x,t)\big)$$

*where $Q_j(x,0) \equiv 0$ for every $j$, and where for each $j \leq q$, $Q_j$ takes the form*

$$(10.9) \qquad Q_j(x,t) = P_j(x,t) + R_j(x,t)$$

*with the following properties:*

(10.10)
$$\begin{aligned}
P_j(x,0) &\equiv R_j(x,0) \equiv 0, \\
P_j(x,t) &\text{ depends only on } t, x_1, \ldots x_{j-1}, \\
P_j &\text{ is a homogeneous polynomial of weight } w_j, \\
\gamma'_t(x') &= (x_1 + P_1(x,t), \ldots x_q + P_q(x,t)) \text{ satisfies } (\mathcal{C}_J) \text{ at } 0 \in \mathbb{R}^q, \\
R_j &\text{ is of weight strictly greater than } w_j, \\
Q_j &\text{ is of weight greater than or equal to } w_q \text{ for every } 1 \leq j \leq n.
\end{aligned}$$

*Suppose further that*

$$(10.11) \qquad \gamma \text{ fails to satisfy } (\mathcal{C}_J) \text{ at } 0 \in \mathbb{R}^n.$$

*Then either*

(A) *All the $Q_j$ have infinite weight, or*

(B) *There exists a diffeomorphism $\Phi$ of a neighborhood of $0 \in \mathbb{R}^n$, fixing 0, such that $\tilde\gamma_t(x) = \Phi \gamma_t \Phi^{-1}(x)$ again takes the general form (10.8),(10.9) but either*

   (B1) *with $q$ increased by 1, or*

   (B2) *with $q$ unchanged, none of the weights of the $Q_j$ decreased, and the weight of one of those $Q_i$, $i > q$, formerly having minimal weight increased by at least one.*

In the special case $q = 0$ we mean the following. The only hypothesis[13] is that $\gamma \in C^\infty$ satisfies $\gamma(x, 0) \equiv x$. We write $\gamma_t(x) = x + (Q_1(x, t), \ldots Q_n(x, t))$, and assign infinite weight to each component of $x$ and weight 1 to each component of $t$.

*Proof.* Consider first the case $q = 0$. If $\partial^\beta Q_j/\partial t^\beta(0,0) = 0$ for every $j$ and every multi-index $\beta$, then alternative (A) holds and the proof is complete. If not, select $\beta$ and $j$ for which this partial derivative is nonzero, such that $|\beta|$ is minimal. Let $\Phi$ be the diffeomorphism of $\mathbb{R}^n$ which interchanges the $j^{\text{th}}$ and first coordinates. Define new weights $w_1 = |\beta|$, and $w_j = +\infty$ for all $j > 1$, and define a new index $q = 1$. Then the conjugated mapping $\tilde\gamma$ takes the form (10.8), where $Q_1(0, t)$ has weight $w_1 < \infty$. Thus alternative (B1) holds, and the proof is again complete.

Suppose now that $1 \leq q < n$. Write $x' = (x_1, \ldots x_q)$. If all the $Q_j$ have infinite weight then alternative (A) holds. Otherwise choose some $i_0 > q$ for which $Q_{i_0}$ has minimal weight, $v$, among all $Q_i$ for $i > q$. Permute the last $n - q$ coordinates so that $i_0 = q + 1$. Let $S(x, t)$ be the sum of all monomials of weight $v$ in $Q_{q+1}$; since $x_{q+1}, \ldots x_n$ all have infinite weight, $S$ can depend only on $x', t$ and will henceforth be denoted by $S(x', t)$. In $\mathbb{R}^{q+1}$ define

$$\tilde\gamma_t(x', x_{q+1}) = (\gamma'_t(x'), x_{q+1} + S(x', t)).$$

There are two cases to analyze, depending on whether or not $\tilde\gamma$ satisfies $(\mathcal{C}_J)$ at $0 \in \mathbb{R}^{q+1}$. If it does, define $P_{q+1} = S$, $R_{q+1} = Q_{q+1} - S$, and $w_{q+1} = v$. Then alternative (B1) holds.

If it does not, invoke Lemma 10.13 to obtain a diffeomorphism $\tilde\Phi$ of $\mathbb{R}^{q+1}$ which conjugates $\tilde\gamma$ to the form $(\gamma'_t(x'), x_{q+1})$, and which takes the form $\tilde\Phi(x', x_{q+1}) = (x', x_{q+1} + h(x'))$ where $h$ is homogeneous of weight $v$. Define a diffeomorphism of $\mathbb{R}^n$ by

$$\Phi(x) = (\tilde\Phi(x_1, \ldots x_{q+1}), x_{q+2}, \ldots x_n).$$

Then $\tilde\gamma_t(x) = \Phi\gamma_t(\Phi^{-1}x)$ satisfies alternative (B2), for

$$\tilde\gamma_t(x) = \left( x_1 + P_1 + R_1, \ldots x_q + P_q + R_q, x_{q+1} + \tilde Q_{q+1}(x,t), \ldots, x_n + \tilde Q_n(x,t) \right)$$

where $\tilde Q_{q+1}(x, t) = Q_{q+1}(x, t) - S(x', t)$, and $\tilde Q_j$ has weight $\geq v$ for each $j \geq q + 2$. Since $Q_{q+1} - S$ has weight strictly greater than $v$, alternative (B2) holds. □

The next lemma is the $C^\infty$ analogue of Lemma 10.3.

LEMMA 10.16. *Suppose that $\gamma$ fails to satisfy the curvature condition $(\mathcal{C}_J)$ at $0 \in \mathbb{R}^n$. Then there exist a nonnegative integer $p < n$ and a coordinate*

---

[13]The hypothesis (10.11) will not be used in the case $q = 0$.

*system $x = (x', x'')$ near the origin in $\mathbb{R}^n = \mathbb{R}^p \times \mathbb{R}^{n-p}$ in which $\gamma$ takes the form*

(10.12) $$\gamma_t(x) = (x_1 + S_1(x,t), \ldots, x_n + S_n(x,t))$$

*where for each $j > p$,*

(10.13) $$S_j(x', 0, t) = O(|x'| + |t|)^N \quad \text{for every } N.$$

Thus the manifold $M = \mathbb{R}^p \times \{0\} \subset \mathbb{R}^n$ is invariant under $\gamma$, to infinite order at 0.

*Proof.* Both curvature conditions $(\mathcal{C}_M)$ and $(\mathcal{C}_J)$ have been proved to be invariant under conjugation of $\gamma$ with diffeomorphisms of $\mathbb{R}^n$, so we are free to change coordinates in the proof.

Certainly the hypotheses of Lemma 10.15 are satisfied with $q = 0$. Apply that lemma. If alternative (A) holds then we are finished. Otherwise make the change of coordinates indicated for alternative (B). If alternative (A) then holds in these new coordinates, we are again finished. If not, repeat the process, stopping if and only if alternative (A) ever arises.

Suppose now that the process were never to stop. Observe that the index $q$ in the lemma can never become equal to $n$, for if it did, (10.10) would say that the leading order part $\tilde{\gamma} = (x_1 + P_1, \ldots, x_n + P_n)$ of $\gamma$ did satisfy $(\mathcal{C}_J)$ at 0. This in turn would imply that $\gamma$ must satisfy $(\mathcal{C}_J)$ at 0, because the Jacobians associated to $\gamma, \tilde{\gamma}$ are related: For each $\xi$, $\tilde{J}_\xi(0, \tau)$ is a homogeneous polynomial of some degree in $\tau$, and $J_\xi(0, \tau)$ equals $\tilde{J}_\xi(0, \tau)$ modulo terms of strictly higher degree. Consequently if $\partial_\tau^\alpha \tilde{J}_\xi(0, \tau) \neq 0$ when $\tau = 0$, then the same holds for $J$, with the same $\xi, \alpha$.

Upon each application of Lemma 10.15, the index $q$ either increases, or remains unchanged. Therefore $q$ eventually stabilizes at some value strictly less than $n$. Thereafter alternative (B2) arises every time that Lemma 10.15 is applied.

In this event each iteration yields a diffeomorphism $\Phi_i$, and we wish to compose them all to obtain the desired change of coordinates. This is possible because if the permutations of coordinates made above for notational convenience are omitted, then each $\Phi_i$ agrees with the identity map at 0, up to an order which tends to infinity with the number of iterations $N$ of this process. Therefore any particular term in the Taylor expansion about 0 of the composition of the first $N$ maps $\Phi_i$ becomes independent of $N$, once $N$ becomes sufficiently large. Therefore the infinite composition is well defined as a formal Taylor series about $x = 0$. By a theorem of E. Borel, there exists a $C^\infty$ function $\Phi_\infty : \mathbb{R}^n \mapsto \mathbb{R}^n$ whose Taylor expansion about the origin concides with this formal series. Then conjugation with $\Phi_\infty$ brings $\gamma$ into the form desired. □

*Remark* 10.2. Implicit in the proofs of Lemmas 10.13 and 10.15 is a construction which applies to any mapping $\gamma$, whether or not it satisfies $(\mathcal{C})$. When $(\mathcal{C})$ is satisfied it produces a system of coordinates in which $\gamma(x,t)$ appears in a sense to be maximally flat near $x = 0$; the higher the order to which all $J_\xi(0,\tau)$ vanish at $\tau = 0$, the larger are the weights $w_j$, roughly speaking. This coordinate system can be used in the same spirit as in the proof of Theorem 8.11 to derive restrictions on pairs of exponents $p,q$ for which an operator $Tf(x) = \int f(\gamma(x,t))K(x,t)dt$, with $K \in C^0$ and $K(x,0) \neq 0$, can map $L^p$ to $L^q$. See Section 20.

## Part 3. Analytic theory

## 11. Statements and reduction to the free case

Let $\gamma$ be a $C^\infty$ mapping $(x,t) \mapsto \gamma(x,t) = \gamma_t(x)$ defined in a neighborhood of the point $(x_0, 0) \in \mathbb{R}^n \times \mathbb{R}^k$, with range in $\mathbb{R}^n$. Note that $\gamma$ is assumed to satisfy the equivalent curvature conditions discussed in Section 8.

Let $K$ be a Calderón-Zygmund kernel in $\mathbb{R}^k$. By this we mean that $K \in C^1(\mathbb{R}^k \setminus \{0\})$ is homogeneous of degree $-k$, and satisfies $\int_{|t|=1} K(t) d\sigma(t) = 0$. We also choose a nonnegative $C^\infty$ cut-off function $\psi$, supported near $x_0$, and a small positive constant $a$. Then we form the singular Radon transform $T$, defined initially for compactly supported $C^1$ functions by

$$(11.1) \qquad T(f)(x) = \psi(x) \operatorname{pv} \int_{|t| \leq a} f(\gamma_t(x)) K(t) dt$$

where $\operatorname{pv} \int_{|t| \leq a} g(t) dt = \lim_{\varepsilon \to 0} \int_{\varepsilon \leq |t| \leq a} g(t) dt$.

The cutoff function $\psi$ and the small constant $a$ serve to localize matters to a small neighborhood of $x_0$, so that everything is well defined. Throughout the argument we implicitly assume $a$ to be chosen to be sufficiently small for various purposes. In particular, we assume always that the map $x \mapsto \gamma_t(x)$ is a diffeomorphism from a neighborhood of the support of $\psi$ to an open subset of $\mathbb{R}^n$, uniformly for every $|t| \leq a'$, for some constant $a' > a$. Then the operator

$$f \mapsto \psi(x) \int_{a \leq |t| \leq a'} f(\gamma_t(x))K(t)dt$$

is bounded on $L^p$ for every $p \in [1, \infty]$, by the very elementary Lemma 14.1, below. This holds provided merely that $K \in L^1$, without any curvature hypothesis on $\gamma$.

Our first main theorem is then:

THEOREM 11.1. *Suppose that $\gamma$ satisfies $(\mathcal{C})$ and that $K$ is as above. Then the operator $T$ defined by (11.1) extends to a bounded operator from $L^p(\mathbb{R}^n)$ to itself, for every $1 < p < \infty$.*

A more general formulation of this theorem is given in Section 21 below.

There is a corresponding maximal theorem. For every continuous $f$ with compact support we define $M(f)$ by

$$(11.2) \qquad M(f)(x) = \sup_{0<r<a} r^{-k} |\psi(x) \int_{|t|\leq r} f(\gamma_t(x))\, dt|.$$

THEOREM 11.2. *Suppose that $\gamma$ satisfies $(\mathcal{C})$ in a neighborhood of the support of $\psi$. Then $M$ extends to a bounded operator from $L^p(\mathbb{R}^n)$ to itself, for every $1 < p \leq \infty$.*

We shall show first that these results can be obtained from the corresponding statements in the lifted situation where the vector fields are free in the sense of Section 5.

Recall that by the exponential Taylor formula (Proposition 8.1), $\gamma_t \sim \exp\left(\sum t^\alpha X_\alpha/\alpha!\right)$. More precisely, since $\gamma$ satisfies $(\mathcal{C})$, there exists $m > 0$ such that

$$(11.3) \qquad \gamma(x,t) = \exp\left(\sum_{0<|\alpha|\leq m} t^\alpha X_\alpha/\alpha!\right)(x) + R(x,t),$$

where $R(x,t) = O(|t|^{m+1})$, and where the $X_\alpha$ and their commutators of degree $\leq m$ span the tangent space at $x_0$. Throughout the discussion, $X_\alpha$ is assigned the degree $|\alpha|$.

Let $X_1, X_2, \ldots X_p$ be an enumeration of $\{X_\alpha : |\alpha| \leq m\}$, and denote by $a_i$ the degree $|\alpha|$ already assigned to $X_i$. Let $\tilde{X}_i$ denote corresponding lifted vector fields, satisfying the conclusions of Proposition 6.1. In the extended space $\{(x,z)\} = \mathbb{R}^n \times \mathbb{R}^{d-n}$ we have

$$\tilde{X}_\alpha = X_\alpha + \sum_{k=1}^{d-n} b_\alpha^k(x,z) \frac{\partial}{\partial z_k}.$$

Denote by $\pi : \mathbb{R}^n \times \mathbb{R}^{d-n} \mapsto \mathbb{R}^n$ the projection by $\pi(x,z) = x$. In conjunction with (11.3), when $t$ is small we define $\tilde{\gamma}$ by

$$(11.4) \qquad \tilde{\gamma}(x,z,t) = \exp\left(\sum_{0<|\alpha|\leq m} t^\alpha \tilde{X}_\alpha/\alpha!\right)(x,z) + \left(R(x,t), 0\right).$$

Then $\tilde{\gamma}(x,z,0) \equiv (x,z)$, and

$$(11.5) \qquad \pi(\tilde{\gamma}(x,z,t)) = \gamma(x,t).$$

Because $R(x,t) = O(|t|^{m+1})$, and because the exponential representation is unique, the vector fields in the representation

$$\tilde{\gamma}(x,z,t) \sim \exp(\sum_\alpha t^\alpha \hat{X}^\alpha/\alpha!)(x,z)$$

must satisfy $\hat{X}_\alpha \equiv \tilde{X}_\alpha$ for every $|\alpha| \leq m$. Since the vector fields $\tilde{X}_\alpha$ with $|\alpha| \leq m$ together with their iterated commutators span the tangent space, $\tilde{\gamma}$ therefore also satisfies curvature condition $(\mathcal{C}_{\mathfrak{g}})$.

Next let $\eta$ be a cut-off function in the space $\mathbb{R}^{d-n}$ which is $C^\infty$, has compact support, and is $= 1$ when $|z| \leq 1$. Lift each function $f$, defined on $\mathbb{R}^n$, to a function $\tilde{f}$ defined on $\mathbb{R}^n \times \mathbb{R}^{d-n}$ by $\tilde{f}(x,z) = \eta(z)f(x)$. Note that

$$\| \tilde{f} \|_{L^p(\mathbb{R}^n \times \mathbb{R}^{d-n})} = C \| f \|_{L^p(\mathbb{R}^n)} .$$

In keeping with (11.1) we define $\tilde{T}$ by

$$(11.1') \qquad \tilde{T}(F)(x,z) = \tilde{\eta}(z)\,\psi(x)\,\text{pv}\int_{|t|\leq a} F(\tilde{\gamma}_t(x,z))\,K(t)dt$$

where $\tilde{\eta}$ is another cut-off function, which equals 1 for $z$ near 0. Then because of (11.5) we have

$$\tilde{T}(\tilde{f})(x,z) \equiv T(f)(x),$$

for $z$ near the origin, if $a$ is small enough. The $L^p$ boundedness of $\tilde{T}$, once established, then implies the corresponding result for $T$. The same argument applies to the maximal operator $M$. We see therefore that it suffices to prove Theorems 11.1 and 11.2 in the "free" case, and for that reason we shall simplify notation henceforth by writing $X_\alpha$ instead of $\tilde{X}_\alpha$. Moreover, if we consider the $X_\alpha$, $|\alpha| \leq m$, together with all their commutators of degree $\leq m$, we obtain a collection $\{X_I\}$ that is a basis for the tangent space at $x_0$; here the multi-index $I$ ranges over those $|I| \leq m$ that are basic, as defined at the end of Section 2. We also write $\tilde{\gamma} = \gamma$, and $d = n$.

## 12. The multiple mapping $\widetilde{\Gamma}$

In accordance with our previous reduction, we assume that for $x$ near $x_0$

$$(12.1) \qquad \gamma_t(x) = \gamma(x,t) = \exp\left(\sum_{0<|\alpha|\leq m} t^\alpha X_\alpha/\alpha!\right)(x) + R(x,t),$$

with $R(x,t) = O(|t|^{m+1})$ as $t \to 0$, where the vector fields $X_\alpha$ and their commutators of degree $\leq m$ span $\mathbb{R}^n$, and are free up to degree $m$. In particular, $\gamma$ satisfies $(\mathcal{C})$.

We next form a variant of the compositions $\Gamma^j$ of several of the $\gamma$'s defined earlier in (8.5). This variant involves both $\gamma$ and $\gamma^{-1}$. Choose a fixed integer $N \geq n$. Define $\widetilde{\Gamma}(x,\tau) = \widetilde{\Gamma}_\tau(x)$ by

$$(12.2) \qquad \widetilde{\Gamma}(x,\tau) = \gamma^{-1}_{t^{2N}} \cdot \gamma_{t^{2N-1}} \ldots \gamma^{-1}_{t^2} \cdot \gamma_{t^1}(x) .$$

Here $\tau = (t^1, t^2, \ldots t^{2N}) \in \mathbb{R}^{2Nk}$ with each $t^j \in \mathbb{R}^k$, and $x$ is near $x_0$, while $\tau$ is small.

We also bring in the local dilation $\delta_r^x$ centered at $x$, described in Section 5. Define $\widetilde{\Gamma}^{(j)}$ by

$$(12.3) \qquad \widetilde{\Gamma}^{(j)}(x,\tau) = \delta_{2^j}^x\left(\widetilde{\Gamma}(x, 2^{-j}\tau)\right).$$

Note that the superscript $j$ does not have the same meaning as in the notation $\Gamma^j$ employed in (8.5).

PROPOSITION 12.1. *For $(x,\tau)$ near $(x_0, 0)$, the mappings $\tau \mapsto \widetilde{\Gamma}^{(j)}(x,\tau)$ from a neighborhood of the origin in $\mathbb{R}^{2Nk}$ to a neighborhood of $x$ in $\mathbb{R}^n$ satisfy the hypotheses of Proposition 7.2, uniformly in $j$ and $x$ for all sufficiently large $j$.*

The possibility of scaling and the resulting uniformity constitute the principal advantage conferred by the lifting procedure.

*Proof.* Consider first the special case arising from the group $N$ whose Lie algebra is $\mathcal{N}_m^{a_1,\ldots a_p}$ (see §§2 and 4). This Lie algebra has generators $Y_1, \ldots, Y_p$, where each index $i$ corresponds to a unique multi-index $\alpha$ satisfying $0 < |\alpha| \leq m$, and where the degree $a_i$ equals $|\alpha|$ for the corresponding multi-index $\alpha$. In this case we take

$$\gamma(t) = \exp\left(\sum_{0<|\alpha|\leq m} t^\alpha Y_\alpha / \alpha!\right)$$

where exp is the exponential mapping from $\mathcal{N}$ to $N$ (as in §4), and we set $\gamma_t(x) = \gamma(x,t) = x \cdot \gamma(t)$, where the product $\cdot$ is the group multiplication. See also example 8.5. We restrict our attention to $x_0 = 0$, the group identity element.

Using the group multiplication formula (4.1), we see that $\widetilde{\Gamma}(x_0, \tau) = \gamma_{t^{2N}}^{-1} \ldots \gamma_{t^1}(x_0)$ is given exactly by

$$(12.4) \qquad \exp\left(\sum_{I \text{ basic}} Q_I(\tau) Y_I\right),$$

and each $Q_I(\tau)$ is a polynomial in $\tau$ which is homogeneous of degree $|I|$. Since $\delta_r^x(y) = \delta_r(y)$ when $x = 0$ in the group case, the homogeneity of the $Q_I$'s is equivalent in that case with the identity

$$\delta_r^x(\widetilde{\Gamma}(x, r^{-1}\tau)) = \widetilde{\Gamma}(x, \tau).$$

Next, since the vector fields $\{Y_\alpha\}$ discussed in Section 4 and their commutators span the Lie algebra $\mathcal{N}_m^{a_1,\ldots a_p}$, curvature condition $(\mathcal{C}_\mathfrak{g})$ is satisfied on $N$. Thus Theorem 8.8 on the equivalence of curvature conditions guarantees the existence of $\beta$ and $\xi$ such that

$$(12.5) \qquad \partial_\tau^\beta J_\xi(0,\tau)\bigg|_{\tau=0} \neq 0,$$

where $\tau$ is restricted to be of the form $\tau = (t^1, 0, t^2, 0, \ldots t^n, 0)$ in the definition of $J_\xi$. Since $Q_I(\tau)$ is homogeneous of degree[14] $|I|$ in $\tau$, it is easy to see that the determinant $\triangle(\tau) = J_\xi(0, \tau)$ is homogeneous of degree $\sum_{I \text{ basic}}(|I|-1) = Q-n$, where $Q$ is the homogeneous dimension of $\mathcal{N}$, and $n$ is its Euclidean dimension. Thus the multi-index in (12.5) must satisfy $|\beta| = Q - n$.

We now pass to the case of the vector fields arising in (12.1), which are free up to degree $m$. Because of (12.1) we also know that

$$\gamma_t^{-1}(x) = \exp\left(\sum_{0<|\alpha|\leq m} -t^\alpha X_\alpha/\alpha!\right)(x) + R'(x,t)$$

with $R'(x,t) = O(|t|^{m+1})$. Hence in this case, by the Baker-Campbell-Hausdorff formula (see Corollary 3.2),

$$\begin{aligned}\widetilde{\Gamma}(x,\tau) &= \exp\left(\sum_{I \text{ basic}} (Q_I(\tau) + O(|\tau|^{m+1}))X_I\right)(x) + O(|\tau|^{m+1}) \\ &= \exp\left(\sum_{I \text{ basic}} (Q_I(\tau) + O(|\tau|^{m+1}))X_I\right)(x),\end{aligned}$$

where the $Q_I$ are the same polynomials as in (12.4). The first equality, as we have seen in identity (5.7), follows because all commutator identities among the $X_I$ having total degrees $\leq m$ are identical with the corresponding commutator identities for the $Y_I$. Moreover the $X_I$ give a basis for the tangent space at $x$; the second identity follows from the first and the fact that the exponential map is a local diffeomorphism.

Applying the local dilations $\delta_r^x$ to the above gives
(12.6)
$$\widetilde{\Gamma}^{(j)}(x,\tau) = \delta_{2^j}^x(\widetilde{\Gamma}(x, 2^{-j}\tau)) = \exp\left(\sum_{I \text{ basic}} \left[Q_I(\tau) + 2^{j|I|}O(|2^{-j}\tau|)^{m+1}\right]X_I\right)(x),$$

since each $Q_I(\tau)$ is homogeneous of degree $|I|$. Note that (12.6) means that in local coordinates given by the exponential map involving the $X_I$ and centered at $x$,

$$\widetilde{\Gamma}^{(j)}(x,\tau) = (Q_I(\tau) + 2^{j|I|}O(|2^{-j}\tau|^{m+1}))_I.$$

Therefore since every multi-index satisfies $|I| \leq m$, the mappings $\tau \mapsto \widetilde{\Gamma}^{(j)}(x,\tau)$ are uniformly of class $C^r$, as $j \to \infty$, for each $r$, as long as $(x, \tau)$ is sufficiently close to $(x_0, 0)$. Moreover, the Jacobian determinant of these mappings equals

$$\triangle(\tau) + O(2^{-j}|\tau|^{Q-n+1}),$$

where $\triangle(\tau)$ is homogeneous of degree $Q-n$, and equals the determinant arising in the case of the Lie algebra $\mathcal{N}_m^{a_1,\ldots a_p}$.

---

[14] That is, $Q_I(r\tau) \equiv r^{|I|}Q_I(\tau)$ for $r > 0$; each $Q_I$ is homogeneous with respect to the Euclidean dilation structure.

Denote by $J_\xi^{(j)}(x,\tau)$ the Jacobian determinant with respect to certain components of $\tau$, formed from $\widetilde{\Gamma}^{(j)}(x,\tau)$. We have now shown that in the free situation, $\gamma$ satisfies $(\mathcal{C}_J)$ uniformly at all points and all scales:

LEMMA 12.2. *For the same $\beta$ and $\xi$ as for $\mathcal{N}_m^{a_1,\ldots a_p}$, for small $\tau$ or for large $j$, and for $x$ sufficiently near $x_0$,*

$$(12.7) \qquad |\partial_\tau^\beta J_\xi^{(j)}(x,\tau)| \geq C > 0$$

*uniformly in $x, \tau, j$.*

This completes the proof of Proposition 12.1. □

## 13. The space $L_\delta^1$

Recall the space $L_\delta^1(\mathbb{R}^n)$ that arose in Section 7. It consists of all functions $h \in L^1(\mathbb{R}^n)$ that satisfy

$$(13.1) \qquad \int_{\mathbb{R}^n} |h(x-y) - h(x)|dx \leq A|y|^\delta \qquad \text{for all } y \in \mathbb{R}^n.$$

The norm on $L_\delta^1$ is defined by $\|h\|_{L_\delta^1} = \|h\|_{L^1} + A$ where $A$ is the smallest constant for which (13.1) holds.

Recall also the map $\Theta$ treated in Section 5. We have $\Theta_x(y) = \Theta(x,y) = u$, which is defined for $y$ near $x$ by the canonical coordinates $u = \{u_I\}$ given by $y = \exp(\sum_{I \text{ basic}} u_I X_I)(x)$. A related quantity is the quasi-distance defined by $d(x,y) = \rho(\Theta(x,y))$ with $\rho(u) = \sum_{I \text{ basic}} |u_I|^{1/|I|}$. The sum $\sum_{I \text{ basic}} |I|$ equals $Q$, the homogeneous dimension of the Lie group $N$.

Fix constants $C_1, C_2 < \infty$. Consider any nonnegative measure $\mu$ defined in $\mathbb{R}^n$, with the following two properties. First, there exist $x \in \mathbb{R}^n$ and integers $j, \nu \geq 0$, such that the support of $\mu$ is contained in the set of all $(y,z)$ satisfying

$$(13.2) \qquad d(y,z) \leq C_1 2^{-j-\nu}, \ d(x,z) \leq C_1 2^{-j}, \ \text{and } d(x,y) \leq C_1 2^{-j}.$$

Second, there exist bounded nonnegative functions $m_1, m_2$ so that for each $f \in C^0(\mathbb{R}^n)$,

$$\iint f(y)\, d\mu(y,z) = \int f(y)\, m_1(y)dy$$

and

$$\iint f(z)\, d\mu(y,z) = \int f(z)\, m_2(z)dz,$$

with

$$(13.3) \qquad m_1(y), m_2(z) \leq C_2.$$

PROPOSITION 13.1. *Suppose $h \in L^1_\delta$ and that the measure $\mu$ is as described above. For each $x \in \mathbb{R}^n$ let*

$$(13.4) \qquad I(h) = 2^{jQ} \int |h(\delta^x_{2^j}(y)) - h(\delta^x_{2^j}(z))| \, d\mu(y,z).$$

*Then provided that $C_1$ is sufficiently small, there exist $\delta', A \in \mathbb{R}^+$ such that*

$$(13.5) \qquad I(h) \leq A 2^{-\nu\delta'} \| h \|_{L^1_\delta}, \quad \text{for all } j, \nu \geq 0$$

*uniformly in $x$, where $\delta'$ depends only on $\delta, m, n$ and $A$ depends only on these quantities and on the constants $C_1, C_2$ in (13.2) and (13.3).*

The functions $h$ to which this result will later be applied depend also on the point $x$, but satisfy (13.2) and (13.3) with uniform bounds, hence satisfy (13.5) with a fixed $\delta', A$.

*Proof.* We observe first that

$$(13.6) \qquad I(h) \leq A \| h \|_{L^1}.$$

In fact,

$$I(h) \leq 2^{jQ} \int |h(\delta^x_{2^j}(y))| \, d\mu(y,z) + 2^{jQ} \int |h(\delta^x_{2^j}(z))| \, d\mu(y,z).$$

The first term equals $2^{jQ} \int |h(\delta^x_{2^j}(y))| m_1(y) dy$ and so is majorized by $C 2^{jQ} \int |h(\delta^x_{2^j}(y))| \, dy$.

Now we make the change of variables $y' = \delta^x_{2^j}(y)$, recalling that on the support of $\mu$, and hence for all $y$ in the support of $m_1$, we have $d(x,y) \leq C_1 2^{-j}$, which shows that the mapping $y \mapsto y' = \delta^x_{2^j}(y)$ is well defined there. To compute its Jacobian, we observe that in the coordinates $u = \Theta_x(y)$ centered at $x$, the mapping $y \mapsto \delta^x_{2^j}(y)$ becomes $u \mapsto \delta_{2^j}(u)$. The Jacobian determinant of the latter is exactly $2^{jQ}$. Thus, $dy \approx 2^{jQ} dy'$, and hence the integral $2^{jQ} \int |h(\delta^x_{2^j}(y))| dy$ is majorized by $C \int |h(y')| dy' = C \| h \|_{L^1}$. The same analysis applies to $2^{jQ} \int |h(\delta^x_{2^j}(z)| d\mu(y,z)$ and so (13.6) is proved.

We next show that

$$(13.7) \qquad I(h) \leq A \left[ \| h \|_{L^1} + \| \nabla^{n+1} h \|_{L^1} \right] \cdot 2^{-\nu/m}.$$

Observe that if $h$ and all of its partial derivatives of order $n+1$ belong to $L^1$, then $h$ is a Lipschitz function, satisfying

$$|h(y) - h(z)| \leq A_h |y - z|, \quad \text{where } A_h = A \cdot (\| h \|_{L^1} + \| \nabla^{n+1} h \|_{L^1}).$$

Now make the change of variable $u = \Theta_x(y) = \Theta(x,y)$. Then

$$\begin{aligned}
\left| h(\delta^x_{2^j}(y)) - h(\delta^x_{2^j}(z)) \right| &\leq A_h |\delta^x_{2^j}(y) - \delta^x_{2^j}(z)| \leq C A_h |\delta_{2^j} \Theta(x,y) - \delta_{2^j} \Theta(x,z)| \\
&\leq A'_h \rho(\delta_{2^j} \Theta(x,y) - \delta_{2^j} \Theta(x,z)) \\
&= A'_h 2^j \rho(\Theta(x,y) - \Theta(x,z)),
\end{aligned}$$

by the identity $\rho(\delta_{2^j} u) = 2^j \rho(u)$, and by the inequality $|u| \le C\rho(u)$, which holds as long as $u$ is bounded.

We now invoke the basic estimate (5.12), which in the present notation states that

(13.8) $\qquad \rho(\Theta(x, y) - \Theta(x, z)) \le C\Big(d(y, z) + d(y, z)^{1/m} d(x, y)^{1-1/m}\Big)$ .

Since $d(y, z) \le C_1 2^{-j-\nu}$ and $d(x, y) \le C_1 2^{-j}$, combining this with the right-hand side of the preceding displayed inequality yields

$$\begin{aligned} |h(\delta^x_{2^j}(y)) - h(\delta^x_{2^j}(z))| &\le CA_h \cdot 2^j \left[2^{-j-\nu} + 2^{-(j+\nu)/m} \cdot 2^{-j(1-1/m)}\right] \\ &\le CA_h 2^{-\nu/m}. \end{aligned}$$

Therefore

(13.9) $\qquad |h(\delta^x_{2^j}(y)) - h(\delta^x_{2^j}(z))| \le CA_h 2^{-\nu/m}.$

Now $\iint d\mu(y, z) = \int m_1(y) dy \le C 2^{-jQ}$, since $m_1$ is supported on the set where $d(x, y) \le C_1 2^{-j}$, and $m_1$ is bounded. Inserting these in the definition (13.4) of $I(h)$ yields (13.7).

Let $h \in L^1_\delta$. Given any $\sigma \ge 1$, $h$ may be decomposed as $h^0_\sigma + h^1_\sigma$ where

(13.10) $\qquad \begin{cases} \| h^0_\sigma \|_{L^1} \le A \| h \|_{L^1_\delta} \cdot \sigma^{-\delta} \\ \\ \| h^1_\sigma \|_{L^1} + \| \nabla^{n+1} h^1_\sigma \|_{L^1} \le A \| h \|_{L^1_\delta} \cdot \sigma^{n+1-\delta}. \end{cases}$

The construction is standard: fix $\Phi \in C^\infty_0$ satisfying $\int_{\mathbb{R}^n} \Phi(x) dx = 1$, and set $\Phi_\sigma = \sigma^n \Phi(\sigma x)$. The functions $h^1_\sigma = h * \Phi_\sigma$ and $h^0_\sigma = h - h * \Phi_\sigma$ then have the properties desired.

Now $I(h) \le I(h^0_\sigma) + I(h^1_\sigma)$. For $I(h^0_\sigma)$ we use the estimate (13.6), and for $I(h^1_\sigma)$ we apply (13.7). The result is

$$I(h) \le A \| h \|_{L^1_\delta} \left[\sigma^{-\delta} + 2^{-\nu/m} \sigma^{n+1-\delta}\right].$$

Setting $\sigma = 2^{\nu/m(n+1)}$ yields (13.5) with $\delta' = \delta/m(n+1)$, concluding the proof. $\square$

*Remark* 13.1. The exponent $\delta' = \delta/m(n+1)$ is not optimal, but suffices for our purpose.

We now formulate as a separate proposition a step in the proof of (13.5). Its proof is the same as that of (13.7).

PROPOSITION 13.2. *Suppose $h$ is a Lipschitz function on $\mathbb{R}^n$. Then*

(13.11) $\qquad I(h) \le C 2^{-\nu/m}$ ,

*where $C$ depends on the Lipschitz norm of $h$.*

The following simple observation regarding our mappings $\gamma$ will be needed.

LEMMA 13.3.

$$d(y, \gamma_t(y)) = O(|t|) \text{ and } d(y, \gamma_t^{-1}(y)) = O(|t|), \tag{13.12}$$

uniformly in $t, y$, for $y$ sufficiently close to $x_0$ and $t$ sufficiently near $0$.

*Proof.* For all $t$,

$$\gamma_t(y) = \exp\left(\sum_{|\alpha| \leq m} t^\alpha X_\alpha/\alpha!\right)(y) + O\left(|t|^{m+1}\right) = \exp\left(\sum_{I \text{ basic}} u_I X_I\right)(y),$$

where

$$u_I = \begin{cases} \dfrac{t^\alpha}{\alpha!} + O(|t|^{m+1}) & \text{when } I = \alpha \\ O(|t|^{m+1}) & \text{otherwise.} \end{cases}$$

Such an expression is possible because the map from $u$ to $\gamma_t(y)$ is a diffeomorphism. Note that $|I| = |\alpha|$ in the first case, and $|I| \leq m$, when $I$ is basic. Thus

$$d(y, \gamma_t(y)) = \sum_{I \text{ basic}} |u_I|^{1/|I|} \leq C|t|.$$

The estimate for $d(y, \gamma_t^{-1}(y))$ is a consequence after we replace $y$ by $\gamma_t^{-1}(y)$. □

*Example* 13.1. Proposition 13.1 will be applied to measures $\mu$ defined by

$$\int f(y, z) \, d\mu(y, z) = 2^{k(j+\nu)} \int_{B(x, c2^{-j})} \left\{\int_{|t| \leq c2^{-j-\nu}} f(y, \gamma_t^{-1}(y)) \, dt\right\} dy \tag{13.13}$$

where $B(x, r) = \{y : d(x, y) \leq r\}$, and $c$ is a small constant. By Lemma 13.3, it is clear that $\mu$ is supported where $d(y, z) \leq c'2^{-j-\nu}$, and the other hypotheses of Proposition 13.1 are easily verified.

## 14. The almost orthogonal decomposition

We shall decompose the operator $T$ arising in Theorem 11.1 as a sum $T = \sum_{j=0}^\infty T_j$, where the summands $T_j$ are almost orthogonal in the appropriate sense. To do this we begin with some preliminary observations. Recall the cut-off function $\psi$ and the small constant $a$ appearing in the statement of Theorem 11.1.

LEMMA 14.1. *The sublinear operator*

$$f \to \psi(x) \int_{|t| \leq a'} |f(\gamma_t(x))| \, dt$$

*is bounded from $L^p(\mathbb{R}^n)$ to itself, for all $p \in [1, \infty]$.*

*Proof.* In fact the operators $f \mapsto \psi \cdot f \circ \gamma_t$ are bounded for all such $p$, uniformly for all $|t| \leq a'$, simply because $x \mapsto \gamma_t(x)$ is a diffeomorphism from a neighborhood of the support of $\psi$ into $\mathbb{R}^n$ for all such $t$, and is invertible, uniformly in $t$. Lemma 14.1 then follows by applying Minkowski's integral inequality. □

Let $1 = \sum_{-\infty}^{\infty} \eta(2^j t)$ be a standard partition of unity on $\mathbb{R}^k \setminus \{0\}$, where $\eta \in C^\infty$ is radial and is supported where $a/2 \leq |t| \leq 2a$. Set $\eta_j(t) = \eta(2^j t)$. Then if $|t| \leq a$,

$$K(t) = \sum_{j=0}^{\infty} K(t)\eta_j(t) = \sum_{j=0}^{\infty} K_j(t).$$

Note that $K_j(t) = K_0(2^j t) 2^{jk}$, because $K(t)$ is homogeneous of degree $-k$; also
Let

(14.1)   $T_j(f)(x) = \psi(x) \int f(\gamma_t(x)) K_j(t)\, dt = \psi(x) \int f(\gamma_{2^{-j}t}(x)) K_0(t)\, dt .$

Then for any compactly supported $f \in C^1$

(14.2)   $$T(f) = \sum_{j=0}^{\infty} T_j(f)$$

where the sum on the right converges uniformly. Convergence can be seen by noting that

$$T_j(f)(x) = \psi(x) \int [f(\gamma_{2^{-j}t}(x)) - f(x)] K_0(t)\, dt = O(2^{-j}), \text{ for every } f \in C^1.$$

PROPOSITION 14.2.   (a) $\| T_j \|_{L^p \to L^p} \leq A < \infty$, for all $1 \leq p \leq \infty$ with $A$ independent of $j$.
(b) *The adjoint $T_j^*$ of $T_j$ takes the form $T_j' + R_j$, where*

$$T_j' f(x) = \int f(\gamma_t^{-1}(x)) \overline{K}_j(t) \overline{\psi}(\gamma_t^{-1}(x))\, dt,$$

and

$$\| R_j \|_{L^p \to L^p} \leq A 2^{-j}, \text{ for all } 1 \leq p \leq \infty.$$

*Proof.* Part (a) is proved in the same way as Lemma 14.1, recalling that $\int |K_j(t)|\, dt = \int |K_0(t)|\, dt < \infty$. And

$$T_j^*(f)(x) = \int f(\gamma_t^{-1}(x))\, \bar\psi(\gamma_t^{-1}(x))\, J(x,t) \bar K_j(t)\, dt,$$

where $J(x,t)$ is the Jacobian determinant of the transformation $x \to \gamma_t^{-1}(x)$. Since $J(x,t) = 1 + O(t)$ and $\bar\psi(\gamma_t^{-1}(x)) = \bar\psi(x) + O(t)$, assertion (b) follows because $\int |t|\, |K_j(t)|\, dt = A' 2^{-j}$. □

With these preliminaries out of the way we come to a pivotal point of the $L^2$ theory of the singular Radon transformation $T$. It is the proof of the almost orthogonal inequalities

$$\text{(14.3)} \qquad \| T_j^* T_i \| + \| T_j T_i^* \| \leq A 2^{-\varepsilon |i-j|} \qquad \text{for some } \varepsilon > 0.$$

Throughout the discussion from here through Section 19, the notation $\| \cdot \|$ with no subscript will denote the $L^2$ operator norm. As is well-known (see [43, Ch. VII, Th. 1]), (14.3) suffices to prove the $L^2$ part of Theorem 11.1.

In proving (14.3), a key observation (see [5]) is that to prove e.g. $\| T_j^* T_i \| \leq A 2^{-\varepsilon |i-j|}$, it suffices to show that there exists some fixed $N$ such that

$$\text{(14.4)} \qquad \| (T_j T_j^*)^N T_i \| \leq A' 2^{-\varepsilon' |i-j|},$$

for some $\varepsilon' > 0$; then $T_j^* T_i$ satisfies (14.3) with $\varepsilon = \varepsilon'/2N$. We indicate this reduction when $N$ is of the form $2^\ell$. Then

$$\| T_j^* T_i \|^2 = \| T_i^* T_j T_j^* T_i \| \leq A \| T_j T_j^* T_i \|$$

by the uniform boundedness of $T_i^*$ guaranteed by Proposition 14.2. This proves the case $N = 1$. For $N = 2$ the same reasoning gives

$$\| T_j T_j^* T_i \|^2 = \| T_i^* T_j T_j^* T_j T_j^* T_i \| \leq A \| (T_j T_j^*)^2 T_i \|;$$

we then can proceed to general $N$ of the form $2^\ell$ by induction.

We shall use the reduction to (14.4) when $i \geq j$. The second inequality $\| T_j T_i^* \| \leq A 2^{-\varepsilon |i-j|}$ in (14.3) can be reduced when $i \geq j$ to

$$\| (T_j^* T_j)^N T_i^* \| \leq A' 2^{-\varepsilon' |i-j|},$$

which will be handled in the same way as (14.4). Similar reductions can be made when $i \leq j$. So we now turn to the task of establishing (14.4).

## 15. Kernel of $(T_j T_j^*)^N$; the $L^2$ theorem

The first step in proving (14.4) is to obtain the following kernel representation of the operator $(T_j T_j^*)^N$.

PROPOSITION 15.1.  *Suppose $N \geq n$. Then there exists $\varepsilon \in (0,1)$ such that for every $j \geq 0$,*

$$\text{(15.1)} \qquad (T_j T_j^*)^N(f)(x) = \int f(\delta_{2^{-j}}^x(y)) \mathcal{K}(x,y) \, dy,$$

*where for each $x$ sufficiently near $x_0$, the function $y \mapsto \mathcal{K}(x,y)$ is supported in a fixed bounded set and is in $L_\varepsilon^1$ uniformly in $x, j$.*

Although $\mathcal{K}$ actually depends on $j$, we have suppressed $j$ in the notation to avoid further encumbering some of the formulas below.

*Proof.* Looking back at Section 14, we see that

$$T_j(f)(x) = \psi(x) \int f(\gamma_{2^{-j}t}(x)) K_0(t) \, dt$$

and

$$T_j^*(f)(x) = \int f(\gamma_{2^{-j}t}^{-1}(x)) \bar{\psi}(\gamma_{2^{-j}t}^{-1}(x)) J(x, 2^{-j}t) \bar{K}_0(t) \, dt.$$

Therefore

$$(15.2) \quad (T_j T_j^*) f(x) = \int f(\gamma_{2^{-j}t_2}^{-1} \cdot \gamma_{2^{-j}t_1}(x)) \psi(x, t_1, t_2) \, dt_1 \, dt_2$$

where

$$\psi(x, t_1, t_2) = \psi(x) \bar{\psi} (\gamma_{2^{-j}t_2}^{-1} \cdot \gamma_{2^{-j}t_1}(x)) J\big(\gamma_{2^{-j}t_1}(x), 2^{-j}t_2\big) K_0(t_1) \bar{K}_0(t_2).$$

Continuing this way we obtain

$$(15.3) \quad (T_j T_j^*)^N (f)(x) = \int f(\widetilde{\Gamma}(x, 2^{-j}\tau)) \psi_j(x; \tau) \, d\tau,$$

where $\tau = (t^1, t^2, \ldots, t^{2N}) \in \mathbb{R}^{2Nk}$, $\widetilde{\Gamma}(x, \tau) = \gamma_{t^{2N}}^{-1} \cdot \gamma_{t^{2N-1}} \cdots \gamma_{t^2}^{-1} \cdot \gamma_{t^1}(x)$, and the functions $\psi_j(x, \tau)$ have compact support, and are in $C^\infty$ with respect to $x$ and $\tau$, uniformly in $j$. Note also that the support of $\psi_j$ as a function of $\tau$ can be restricted to a small neighborhood of the origin, uniformly in $j$, by making the constant $a$ small enough, since each $K_0(t)$ is supported in the set $|t| \leq 2a$.

(15.3) may be rewritten as

$$(15.4) \quad (T_j T_j^*)^N (f)(x) = \int f(\delta_{2^{-j}}^x(\widetilde{\Gamma}^{(j)}(x, \tau))) \psi(x; \tau) \, d\tau,$$

where $\widetilde{\Gamma}^{(j)}(x, \tau) = \delta_{2^j}^x(\widetilde{\Gamma}(x, 2^{-j}\tau))$ is as defined in (12.3); we have also suppressed the dependence on $j$ of $\psi_j$ in the notation. Thus

$$(T_j T_j^*)^N (f)(x) = \int (f \circ \delta_{2^{-j}}^x)(y) \, d\mu(y)$$

where $d\mu = \widetilde{\Gamma}^{(j)}(x, \tau)_*(\psi(x, \tau) \, d\tau)$.

By Proposition 12.1, each mapping $\tau \mapsto \widetilde{\Gamma}^{(j)}(x, \tau)$ satisfies the hypotheses of Proposition 7.2. Therefore $d\mu(y) = \mathcal{K}(x, y) \, dy$, where $\mathcal{K}(x, \cdot) \in L^1_\varepsilon$ for some $\varepsilon > 0$, for each $x, j$. The uniformities asserted in those two propositions establish the uniformity claimed in the present proposition, whose proof is then complete. □

Denote by $\|S\|_{L^p \to L^q}$ the norm of any operator $S$, as a mapping from $L^p$ to $L^q$.

LEMMA 15.2. *There exist $\delta' > 0$ and $C < \infty$ such that*

$$(15.5) \quad \| (T_j T_j^*)^N T_i \|_{L^\infty \to L^\infty} \leq C 2^{-\delta'(i-j)}, \qquad \text{for all } i \geq j \geq 0.$$

*Proof.* By (15.1), for any $f \in L^\infty$,

$$(15.6) \quad (T_j T_j^*)^N T_i(f)(x) = \int f(\gamma_t \delta_{2^{-j}}^x(y')) \mathcal{K}(x, y') \psi(\delta_{2^{-j}}^x(y')) K_i(t) \, dy' dt .$$

In this formula we make the change of variables $y = \delta_{2^{-j}}^x(y')$. Since $\mathcal{K}(x, y')$ is supported in a fixed bounded region independent of $j$, $y$ is restricted to $B(x, C2^{-j})$ for a certain constant $C$. Recalling the mapping $\Theta_x$ defined in Section 5, we can express $y$ as $y = \Theta_x^{-1} \delta_{2^{-j}} \Theta_x(y')$. Thus the Jacobian determinant of the mapping $y \mapsto y'$ is the product of three determinants; since that corresponding to $\delta_{2^{-j}}$ is exactly $2^{-jQ}$, we see that $dy' = 2^{jQ} J(x, y') dy$ where $J \in C^\infty$ may depend on $j$, but all its derivatives are bounded uniformly in $j$. Consequently the right-hand side of (15.6) may be rewritten as

$$(15.7) \quad 2^{jQ} \int f(\gamma_t(y)) \mathcal{K}_j(x, \delta_{2^j}^x(y)) K_i(t) \, dy \, dt$$

where $\mathcal{K}_j(x, y') = \mathcal{K}(x, y') \psi(\delta_{2^{-j}}^x(y')) J(x, y')$, and $y' \mapsto \mathcal{K}_j(x, y')$ belongs to $L_\delta^1$, uniformly in $j$ and in $x$.

Substituting $z = \gamma_t(y)$ in (15.7) and noting that $\det(\partial y / \partial z) = 1 + O(t)$,

$$(T_j T_j)^* T_i(f)(x) = 2^{jQ} \int K(x, y) f(y) + E(f)$$

where

$$(15.8) \quad K(x, y) = \int \mathcal{K}_j(x, \delta_{2^j}^x \gamma_t^{-1}(y)) K_i(t) dt$$

and the error operator $E$ satisfies

$$|E(f)(x)| \leq C 2^{jQ} \int \left| \mathcal{K}_j(x, \delta_{2^j}^x \gamma_t^{-1}(y)) f(y) K_i(t) \right| |t| \, dt \, dy .$$

Now $\int |t| |K_i(t)| dt \leq C 2^{-i}$, because $K_i$ is supported where $|t| \leq a 2^{-i}$, and $|K_i(t)| \leq C|t|^{-k}$. Therefore uniformly in $x$,

$$(15.9) \quad |Ef(x)| \leq C 2^{-i} \|f\|_\infty.$$

The kernel in (15.8) can be rewritten as

$$(15.10) \quad K(x, y) = \int \left[ \mathcal{K}_j(x, \delta_{2^j}^x \gamma_t^{-1}(y)) - \mathcal{K}_j(x, \delta_{2^j}^x(y)) \right] K_i(t) dt,$$

since $\int K_i(t) dt = 0$. Therefore it is dominated by

$$C 2^{ki} \int_{|t| < c 2^{-i}} \left| \mathcal{K}_j(x, \delta_{2^j}^x \gamma_t^{-1}(y)) - \mathcal{K}_j(x, \delta_{2^j}^x(y)) \right| dt .$$

We claim that as a result

(15.11) $$2^{jQ} \int |K(x,y)|dy \leq C 2^{-(i-j)\delta'}.$$

To see this, we need only to recall that the functions $y \mapsto \mathcal{K}_j(x,y)$ are uniformly in $L^1_\delta$, and to utilize Proposition 13.1 with $h(y) = \mathcal{K}_j(x,y)$, $\nu = i - j$; the measure $d\mu$ is as detailed in Example 13.1. If $a$ is chosen to be sufficiently small then the constant $C_1$ in (13.2) will be as small as required. Thus (15.11), together with (15.9), establishes (15.5). □

We next combine (15.5) with the easy estimate

$$\| (T_j T_j^*)^N T_i \|_{L^1 \to L^1} \leq C,$$

which follows directly from Proposition 14.2. A consequence, either by a straightforward argument or by the Riesz convexity theorem, is that

$$\| (T_j T_j^*)^N T_i \|_{L^2 \to L^2} \leq C 2^{-(\delta'/2)(i-j)}.$$

This is (14.4) with $\varepsilon' = \delta'/2$, and with it we conclude the proof of Theorem 11.1 in the case $p = 2$.

## 16. The $L^p$ argument; preliminaries

We turn to the $L^p$ inequalities for the operator $T$ stated in Theorem 11.1. An earlier method [45] used for the case of operators invariant under the Euclidean group structure of $\mathbb{R}^n$ proceeded by embedding $T$ in an analytic family of operators $T(s)$, with $T(0) = T$. The operators $T(s)$ were more singular than $T(0)$ when $\mathrm{Re}(s) > 0$, but were nonetheless bounded on $L^2$; while for $\mathrm{Re}(s) < 0$, $L^p$ bounds could be obtained for all $1 < p < \infty$, because the operators $T(s)$ were of a more standard type, treatable by Calderón-Zygmund theory. An interpolation theorem led to the desired result for $T(0) = T$.

We proceed instead[15] as follows, using [5],[34]. First we construct $\{S_j\}$, a family of averaging operators with smooth densities adapted to the geometry underlying our situation, so that $S_j$ tends to the identity operator $I$ as $j \to \infty$. Define $R_j = S_{j+1} - S_j$, and note that formally $I = S_j + \sum_{\nu=0}^{\infty} R_{j+\nu}$.

Together with the decomposition $T = \sum_{j \geq 0} T_j$, this yields $T_j = T_j S_j + \sum_{\nu \geq 0} T_j R_{j+\nu}$ and so

(16.1) $$T = \sum_{\nu=0}^{\infty} U_\nu + \sum_{j=0}^{\infty} T_j S_j,$$

---

[15] A close analogue of the analytic family of operators mentioned above could be constructed by setting $T(s) = \sum_{\nu \geq 0} 2^{\nu s} U_\nu + \sum_{j \geq 0} T_j S_j$.

where

$$U_\nu = \sum_{j=0}^{\infty} T_j R_{j+\nu}. \tag{16.2}$$

These series converge in the strong operator topology on the space of bounded linear operators from $C_0^1$ to $L^2$, as follows from uniform boundedness of the $S_j$, the almost orthogonality inequality (14.3), and the estimate $\|f - S_j f\|_2 = O(2^{-j})$ for every $f \in C_0^1$, which will be an easy consequence of the definition of $S_j$. Our aim is to show that

(16.3) $\quad \| U_\nu \|_{L^2 \to L^2} \leq A 2^{-\nu \varepsilon}$, for some $\varepsilon > 0$.

(16.4) $\quad \| U_\nu \|_{L^p \to L^p} \leq A_{p,\varepsilon'} 2^{\nu \varepsilon'}$, for every small $\varepsilon' > 0$ and $1 < p \leq 2$.

(16.5) $\quad \| \sum_j T_j S_j \|_{L^p \to L^p} \leq A_p$, for $1 < p \leq 2$.

Interpolation between (16.3) and (16.4) then gives the better estimate

$$\| U_\nu \|_{L^p \to L^p} \leq A_p 2^{-\nu \varepsilon_p}, \text{ for some } \varepsilon_p > 0, \text{ for every } 1 < p \leq 2.$$

This, in view of (16.1) and (16.5), implies the desired conclusion for $T$, when $1 < p \leq 2$; a duality argument then leads to the same result for $2 < p < \infty$.

We begin by defining the averaging operators $S_j$. Fix a $C_0^\infty$ nonnegative function $\Phi$, which is even, is identically equal to 1 near the origin, is supported where $|u| < a$, and satisfies $\int_{\mathbb{R}^n} \Phi(u)\,du = 1$. Let $\Phi_j(u) = 2^{jQ} \Phi(\delta_{2^j}(u))$. Fix a cutoff function $\mathcal{X} \in C_0^\infty(\mathbb{R}^n)$ which is $\equiv 1$ near $x_0$, and has appropriately small support. Let $\mathcal{X}_0(x) = \mathcal{X}(x) J(x,x)^{1/2}$, where

$$J(x,y) = \left| \det \left( \partial \Theta_x(y)/\partial y \right) \right|,$$

and define $S_j$ by

$$S_j(f)(x) = \mathcal{X}_0(x) \int \Phi_j(\Theta(x,y)) \mathcal{X}_0(y) f(y)\,dy.$$

Denote by $S_j(x,y) = \mathcal{X}_0(x) \Phi_j(\Theta(x,y)) \mathcal{X}_0(y)$ the associated distribution kernel function, and by $S_j^*$ the adjoint operator.

PROPOSITION 16.1. *There exists $C < \infty$ such that for every $j$, $S_j(x,y) \equiv 0$ whenever $d(x,y) \geq C 2^{-j}$. Moreover*

(a) $\int S_j(x,y)\,dy = \mathcal{X}^2(x) + O(2^{-j})$ as $j \to \infty$.

(b) $S_j(f)(x) \to \mathcal{X}^2 f(x)$, uniformly as $j \to \infty$ for any continuous $f$.

(c) $\| S_j \|_{L^p \to L^p} \leq A < \infty$, for all $j$ and all $p \in [1, \infty]$.

(d) $S_j^* = S_j$.

(e) *For all $x, y_1, y_2$,*
$$|S_j(x, y_1) - S_j(x, y_2)| \leq A 2^{j(Q+1/m)} d(y_1, y_2)^{1/m}.$$

The next result is a direct consequence.

COROLLARY 16.2. *Let $R_j = S_{j+1} - S_j$, and let $R_j(x, y)$ be its distribution kernel. Then conclusions* (c), (d), *and* (e) *hold with $R_j$ substituted for $S_j$. Moreover $\int R_j(x, y) dy = O(2^{-j})$, and $\int |R_j(x, y)| dy \leq A < \infty$, uniformly in $x, j$.*

*Proof of Proposition* 16.1. To prove conclusion (a), note that with $u = \Theta(x, y) = \Theta_x(y)$, we have $dy = J(x, y)^{-1} du$. Thus
$$\int S_j(x, y) dy = \mathcal{X}_0(x) \int \Phi_j(u) \left[ \mathcal{X}_0(y) J(x, y)^{-1} - \mathcal{X}_0(x) J(x, x)^{-1} \right] du$$
$$+ \mathcal{X}_0^2(x) J(x, x)^{-1} \int \Phi_j(u) du.$$

The second term on the right equals $\mathcal{X}^2(x)$, because $\mathcal{X}_0(x) = \mathcal{X}(x)(J(x,x))^{1/2}$, and $\int \Phi_j(u) du = 1$. The first term is $O(2^{-j})$, since
$$|\mathcal{X}_0(y) J(x, y)^{-1} - \mathcal{X}_0(x) J(x, x)^{-1}| \leq A|x - y| \leq A'|u| \leq A''\rho(u)$$
and $\int \Phi_j(u) \rho(u) du = c 2^{-j}$, because $\rho(\delta_{2^{-j}}(u)) = 2^{-j} \rho(u)$. Thus (a) is proved.

Since $S_j(x, y) \geq 0$, (a) implies that the $S_j$ are uniformly bounded as operators on $L^\infty$. Conclusion (b) follows from this, because the support of $y \mapsto S_j(x, y)$ shrinks to $x$ as $j \to \infty$. Indeed, $S_j(x, y)$ is supported where $\rho(\Theta(x, y)) \leq C 2^{-j}$, that is, where $d(x, y) \leq C 2^{-j}$.

Since (c) holds for $p = \infty$ and $S_j(x, y) \equiv S_j(y, x)$, it holds also for $p = 1$, and thus also for all $p \in [1, \infty]$. The symmetry of $S_j(x, y)$ also gives (d).

Since $|S_j(x, y)| \leq C 2^{jQ}$, the inequality in (e) holds whenever $d(y_1, y_2) \geq c 2^{-j}$, for any fixed small constant $c$. Assuming henceforth that $d(y_1, y_2) \leq c 2^{-j}$, it holds also whenever $d(x, y_1) \geq C' 2^{-j}$, provided $C'$ is chosen to be a sufficiently large constant, by the quasi-triangle inequality. Thus we may also assume $d(x, y_1) \leq C' 2^{-j}$. Under these conditions
$$(16.6) \quad |S_j(x, y_1) - S_j(x, y_2)| \leq A |\Phi_j(\Theta(x, y_1)) - \Phi_j(\Theta(x, y_2))|$$
$$+ A' |\Phi_j(\Theta(x, y_1))| \cdot |y_1 - y_2|.$$

Since $\Phi_j(u) = 2^{jQ} \Phi(\delta_{2^j} u)$,[16]

---

[16] The range of $\Theta$ is $\mathbb{R}^n$, equipped with a distinguished coordinate system $(u_I)$; addition and subtraction of elements of this target space are defined with respect to the Euclidean group structure associated to this coordinate system, not with respect to the nilpotent group $N$. The dilations $\delta_r$ are automorphisms of both the abelian group $\mathbb{R}^n$ and $N$.

$$\begin{aligned}|\Phi_j(\Theta(x,y_1))-\Phi_j(\Theta(x,y_2))| &\leq 2^{jQ}|\delta_{2^j}\Theta(x,y_1)-\delta_{2^j}\Theta(x,y_2)|\\ &\leq A2^{jQ}\rho\Big(\delta_{2^j}(\Theta(x,y_1)-\Theta(x,y_2))\Big)\\ &= A2^{j(Q+1)}\rho\Big(\Theta(x,y_1)-\Theta(x,y_2)\Big).\end{aligned}$$

Invoking (13.8) yields

$$\begin{aligned}\rho(\Theta(x,y_1)-\Theta(x,y_2)) &\leq C\left[(d(y_1,y_2))+d(y_1,y_2)^{1/m}d(x,y_1)^{1-1/m}\right]\\ &\leq C'd(y_1,y_2)^{1/m}2^{-j(1-1/m)},\end{aligned}$$

since $d(y_1,y_2) \leq c2^{-j}$ and $d(x,y_1) \leq C2^{-j}$. Therefore the first term on the right side of (16.6) is majorized by

$$C2^{j(Q+1)}\,2^{-j(1-1/m)}\,d(y,y_2)^{1/m} = C2^{j(Q+1/m)}\,d(y_1,y_2)^{1/m}.$$

The second term in (16.6) also contributes $O\Big(2^{j(Q+1/m)}d(y_1,y_2)^{1/m}\Big)$, because

$$|\Phi_j| \leq C2^{jQ},\quad \text{and}\quad |y_1-y_2| \leq Cd(y_1,y_2) \leq Cd(y,y_2)^{1/m}\,2^{-j(1-1/m)}.$$

Hence assertion (e) is established, and the proposition is proved. $\square$

We now adjust the cutoff function $\mathcal{X}$ in the definition of $S_j$, and the constant $a$ in the definition (11.1) of $T$, so that on the support of $\psi$, $\mathcal{X}(\gamma_t(x)) \equiv 1$, for all $|t| \leq a$. Then because of Proposition 16.1 (b),

$$T_j(f) = \lim_{i\to\infty} T_j S_i(f) = T_j S_j(f) + \sum_{\nu=0}^{\infty} T_j R_{j+\nu}(f)$$

for all $f \in C_0^1$, and so the decomposition (16.1) is established.

## 17. Further $L^2$ estimates

We begin with a reprise of the $L^2$ theory of Sections 14 and 15, involving some easy modifications. We consider $U_\nu = \sum_{j\geq 0} T_j R_{j+\nu}$ and prove (16.3), i.e.

(17.1) $$\|U_\nu\|_{L^2\to L^2} \leq A2^{-\varepsilon\nu} \quad \text{for some } \varepsilon > 0.$$

By the almost orthogonality argument already used, it suffices to see that

(17.2) $$\|(T_j R_{j+\nu})^*(T_i R_{i+\nu})\|_{L^2\to L^2} + \|(T_j R_{j+\nu})(T_i R_{i+\nu})^*\|_{L^2\to L^2}$$
$$\leq A2^{-\nu\varepsilon}2^{-\varepsilon|i-j|}.$$

Each term on the left can be majorized using the geometric mean of two of the following three estimates, if we recall that $T_j$ and $R_j$ have uniformly bounded norms, and that $R_j^* = R_j$.

(17.3) $\qquad \|T_j^* T_i\| \leq A 2^{-\varepsilon'|i-j|} \qquad$ for some $\varepsilon' > 0$,

(17.4) $\qquad \|R_j R_i\| \leq A 2^{-|i-j|/m}$,

(17.5) $\qquad \|T_j R_{j+\nu}\| \leq A 2^{-\varepsilon' \nu} \qquad$ for some $\varepsilon' > 0$.

The first was already proved in Section 14. To prove the second, note that by taking adjoints if necessary, we may assume that $i \geq j$. We write $R_j R_i(f)(x) = \int K(x,y) f(y) dy$, where $K(x,y) = \int R_j(x,z) R_i(z,y) dz$. Thus

$$K(x,y) = \int [R_j(x,z) - R_j(x,y)] R_i(z,y) dz + O(2^{-i}|R_j(x,y)|)$$

by Corollary 16.2. This splits $K(x,y)$ as $K'(x,y) + E(x,y)$. The operator with kernel $E(x,y)$ has a norm which is clearly $O(2^{-i}) = O(2^{-|i-j|})$, since $i \geq j$. To estimate

$$K'(x,y) = \int [R_j(x,z) - R_j(x,y)] R_i(z,y) dz,$$

observe that $K'(x,y)$ is supported where $d(x,y) \leq c 2^{-j}$. Moreover, by conclusion (e) of Corollary 16.2,

$$|R_j(x,z) - R_j(x,y)| \leq C 2^{j(Q+1/m)} d(z,y)^{1/m}.$$

However

$$\int d(z,y)^{1/m} |R_i(z,y)| dz \leq C 2^{-i/m} \int |R_i(z,y)| dz \leq C 2^{-i/m},$$

since $\int |R_i(z,y)| dz \leq C$, and $R_i(z,y)$ is supported where $d(z,y) \leq C 2^{-i}$. It follows that

$$\int |K'(x,y)| dx \leq C 2^{jQ} 2^{-(i-j)/m} |B(x, C 2^{-j})| \leq C' 2^{-|i-j|/m}.$$

The same bound holds for $\int |K'(x,y)| dy$, and (17.4) follows directly.

We now come to (17.5). As in the proof of (17.3) in Section 14, it suffices to show

(17.6) $\qquad \|(T_j^* T_j)^N R_{j+\nu}\| \leq A 2^{-\varepsilon \nu}, \quad$ for some $\varepsilon > 0$.

As in (15.1), the kernel $K(x,y)$ associated to $(T_j^* T_j)^N R_{j+\nu}$ is

(17.7) $$K(x,y) = \int \mathcal{K}(x, z') R_{j+\nu}(\delta_{2^{-j}}^x (z'), y) dz'$$

where $z' \mapsto \mathcal{K}(x, z')$ is uniformly in $L_\delta^1$, for $x$ near $x^0$ and has uniformly compact support in $z'$. Here we are dealing with $(T_j^* T_j)^N$ instead of $(T_j T_j^*)^N$, but this distinction amounts merely to interchanging the roles of $\gamma$ and $\gamma^{-1}$ in (15.1).

When we make the change of variables $z' = \delta_{2^j}^x(z)$, (17.7) becomes

$$K(x,y) = 2^{jQ} \int \mathcal{K}(x, \delta_{2^j}^x(z)) J(x, \delta_{2^j}^x(z)) R_{j+\nu}(z,y) dz$$

where
$$J(x, z') = 2^{-jQ}\left|\det\left(\partial z'/\partial z\right)\right| = 2^{-jQ}\left|\det\left(\partial \delta^x_{2^{-j}}(z')\big/\partial z'\right)\right|^{-1} = 1 + O(2^{-j})$$

in any $C^N$ norm as a function of $z'$, uniformly in $x, j$ provided that $|z'| = O(1)$ and $j \geq 0$. Therefore, fixing $x$ and letting $h(z) = \mathcal{K}(x, z)J(x, z)$, we see that $K(x, y)$ splits as

$$K'(x, y) + E(x, y) = 2^{jQ} \int \left[h(\delta^x_{2^j}(z)) - h(\delta^x_{2^j}(y))\right] R_{j+\nu}(z, y) dz$$
$$+ O(2^{-j-\nu} 2^{jQ} |h(\delta^x_{2^j}(y))|),$$

where the bound for $E(x, y)$ holds because $\int R_{j+\nu}(z, y)\, dz = O(2^{-j-\nu})$.

Since $2^{jQ} \int |h(\delta^x_{2^j}(y))| dy = O(1)$, the operator whose kernel is $E(x, y)$ has $L^\infty$ operator norm bounded by $C 2^{-j-\nu} \leq C 2^{-\nu}$. For the main term, $\int |K'(x, y)| dy$ is estimated by means of Proposition 13.1. Since $h(z) = \mathcal{K}(x, z)J(x, z)$ has small compact support as a function of $z$ (provided the constant $a$ is chosen to be sufficiently small) and belongs to $L^1_\delta$, uniformly in the parameters $j, x$, we can estimate $\int |K'(x, y)| dy$ by a constant multiple of $I(h)$, as defined in (13.4). Here $d\mu(y, z) = 2^{j+\nu} \chi(y, z)\, dy\, dz$, where $\chi$ is the characteristic function of the set where $d(y, z) \leq C_1 2^{-(j+\nu)}$, $d(x, y) \leq C_1 2^{-j}$, and $d(x, z) \leq C_1 2^{-j}$.

The conclusion is that $\int |K(x, y)| dy \leq C 2^{-\nu \delta'}$, and therefore the $L^\infty$ operator norm of $(T_j^* T_j)^N R_{j+\nu}$ is $O(2^{-\nu \delta'})$. Since this operator is a composition of operators bounded on $L^1$, uniformly in $j, \nu$, we deduce (17.5) with $\varepsilon' = \delta'/2$ by interpolation. Therefore (17.1) is proved.

We next turn to the $L^2$ estimate for $\sum_{j \geq 0} T_j S_j$. Using almost orthogonality again, it suffices to show that

$$\| (T_j S_j)^* T_i S_i \| + \| (T_j S_j)(T_i S_i)^* \| \leq A 2^{-\varepsilon |i-j|} \qquad \text{for some } \varepsilon > 0.$$

The bound for the first term is a consequence of the estimates of $T_j^* T_i$ already proved in Sections 14 and 15.

To treat the second, we may suppose $i \geq j$, and revisit the proof of (17.4). The $S$'s satisfy the same size and difference inequalities as the $R$'s, but not their cancellation property. The key observation is that the composed operators $S_i T_i^*$ nonetheless do effectively share this cancellation property. To begin, we may replace $T_i^*$ by $T_i'$, as defined in Proposition 14.2, making an error whose operator norm is $O(2^{-i})$, by the use of part (b) of Proposition 14.2. Now the distribution kernel $S'_i(x, y)$ associated to $S_i T_i'$ is given by

$$S'_i(x, y) = \bar{\psi}(y) \int S_i(x, \gamma_t^{-1}(y))\, \bar{K}_i(t)\, dt.$$

Hence $\int S'_i(x, y)\, dx = O(2^{-i})$, since

$$\int K_i(t)\, dt = 0, \int |K_i(t)|\, dt \leq C < \infty, \text{ and } \int S_i(x, y)\, dx = \mathcal{X}^2(y) + O(2^{-i}).$$

Invoking Lemma 13.3, we see moreover that $S'_i(x,y)$ is supported where $d(x,y) \leq C'2^{-i}$, since $S_i(x,y)$ is supported where $d(x,y) \leq C2^{-i}$. The rest of the proof of the inequality $\|S_j S_i T_i^*\| \leq A2^{-\varepsilon|i-j|}$ then follows that of (17.4) closely.

## 18. The $L^p$ estimates; conclusion

We consider the operator $U_\nu = \sum_{j \geq 0} T_j R_{j+\nu}$. We shall see below that there exists a kernel function $U_\nu(x,y)$, satisfying

$$\int U_\nu(f)(x)\bar{g}(x)dx = \int U_\nu(x,y)f(y)\,\bar{g}(x)\,dy\,dx$$

whenever $f$ and $g$ are (say) continuous with *disjoint* compact supports. For this kernel we shall prove the main estimate needed to apply the Calderón-Zygmund theory.

PROPOSITION 18.1. *There exists $\bar{c} < \infty$ such that for every small $\varepsilon > 0$ there exists $A_\varepsilon < \infty$ satisfying*

(18.1)
$$\int_{d(x,y_1) \geq \bar{c}d(y_1,y_2)} |U_\nu(x,y_1) - U_\nu(x,y_2)|\,dx \leq A_\varepsilon 2^{\varepsilon\nu} \quad \text{for all distinct } y_1, y_2.$$

We shall need to examine the distribution kernels of the operators $T_j R_{j+\nu}$. Define $K_{j,\nu}(x,y)$ by

$$T_j R_{j+\nu}(f)(x) = \int K_{j,\nu}(x,y)f(y)dy.$$

LEMMA 18.2. *The kernels $K_{j,\nu}$ have the following properties.*

  (i) $K_{j,\nu}(x,y)$ *is supported where* $d(x,y) \leq C2^{-j}$.

  (ii) $\int |K_{j,\nu}(x,y)|dx \leq C < \infty$, *uniformly in $y$, $j$ and $\nu$.*

  (iii) $\int |K_{j,\nu}(x,y_1) - K_{j,\nu}(x,y_2)|\,dx \leq Cd(y_1,y_2)^\varepsilon 2^{(j+\nu)\varepsilon}$, *for every $\varepsilon \in [0, m^{-1}]$.*

*Proof.* By the definitions of $R_{j+\nu}$ and $T_j$,

(18.2) $$K_{j,\nu}(x,y) = \psi(x)\int R_{j+\nu}(\gamma_t(x),y)\,K_j(t)dt.$$

Since $d(\gamma_t(x),x) \leq A|t|$ by Lemma 13.3, and $R_{j+\nu}(x,y)$ is supported where $d(x,y) \leq C2^{-j-\nu}$, conclusion (i) follows. Next,

$$\int |K_{j,\nu}(x,y)|dx \leq C \iint |R_{j+\nu}(\gamma_t(x),y)K_j(t)|\,dt\,dx$$
$$\leq C' \int |R_{j+\nu}(x,y)|dx \cdot \int |K_j(t)|\,dt$$

and conclusion (ii) follows from Corollary 16.2 since $\int |K_j(t)|dt \equiv C < \infty$.

Finally,
$$|R_{j+\nu}(x,y_1) - R_{j+\nu}(x,y_2)| \leq c2^{(j+\nu)(Q+1/m)} d(y_1,y_2)^{1/m}$$

by conclusion (e) of Proposition 16.1. We may assume that $d(y_1,y_2) \leq 2^{-j-\nu}$, for otherwise (iii) follows directly from (ii). Then
$$\int |R_{j+\nu}(x,y_1) - R_{j+\nu}(x,y_2)| \, dx \leq C2^{(j+\nu)/m} d(y_1,y_2)^{1/m}$$

because the integration is limited to a ball of radius $\leq C2^{-j-\nu}$. This result, combined with (18.2), gives conclusion (iii) with $\varepsilon = 1/m$. Since by (ii) it holds for $\varepsilon = 0$, the full conclusion with $0 \leq \varepsilon \leq 1/m$ is an immediate consequence. □

We turn to the proof of the proposition. The almost orthogonality property (17.2) $\sum_j T_j R_{j+\nu}$ implies (see [43, p. 318]) that $\sum_j T_j R_{j+\nu}(f)$, converges in the $L^2$ norm, for every $f \in L^2$. Thus for any $f, g \in L^2$,
$$(U_\nu f, g) = \sum_{j=0}^{\infty} (T_j R_{j+\nu} f, g) = \sum_j \int f(y) K_{j,\nu}(x,y) \bar{g}(x) \, dx \, dy.$$

If $f$ and $g$ have disjoint compact supports, then only finitely many terms $K_{j,\nu}(x,y)$ are nonzero on the support of $f(y)\bar{g}(x)$, by (i) of Lemma 18.1. Thus $\sum_j K_{j,\nu}(x,y)$ is the kernel of $U_\nu$ in the sense claimed at the beginning of Section 18.

We have therefore $U_\nu(x,y) = \sum_{j=0}^{\infty} K_{j,\nu}(x,y)$, and hence
$$|U_\nu(x,y_1) - U_\nu(x,y_2)| \leq \sum_j |K_{j,\nu}(x,y_1) - K_{j,\nu}(x,y_2)|.$$

Now fix any index $j \geq 0$, and consider
$$\int_{d(x,y_1) \geq \bar{c} d(y_1,y_2)} |K_{j,\nu}(x,y_1) - K_{j,\nu}(x,y_2)| dx .$$

For this integral to be nonzero, there must be an $x$ in the support of $K_{j,\nu}(x,y_1)$ or $K_{j,\nu}(x,y_2)$, with $d(x,y_1) \geq \bar{c} d(y_1,y_2)$. Then in view of conclusion (i) of the lemma, either $d(x,y_1) \leq C2^{-j}$ or $d(x,y_2) \leq C2^{-j}$. In the latter case,
$$d(y_1,y_2) \leq \bar{c}^{-1} d(x,y_1) \leq C\bar{c}^{-1}\Big(d(x,y_2) + d(y_1,y_2)\Big),$$

so $d(y_1,y_2) \leq c' 2^{-j}$ provided that $\bar{c} > C^{-1}$; the same conclusion follows more directly in the former case.

Invoking conclusion (iii) of the lemma gives for any $0 < \varepsilon \leq 1/m$,
$$\int_{d(x,y_1) \geq \bar{c} d(y_1,y_2)} |U_\nu(x,y_1) - U_\nu(x,y_2)| \, dx \leq C_\varepsilon \sum_j d(y_1,y_2)^\varepsilon \, 2^{(j+\nu)\varepsilon},$$

where according to the conclusion of the preceding paragraph, the sum need be taken only over all $j$ satisfying $d(y_1, y_2) \leq c'2^{-j}$. Summing the geometric series establishes (18.1).

The $L^2$ boundedness of the operators $U_\nu$, guaranteed by (17.1), together with the inequality (18.1) regarding their kernels that we have just proved, allow us to appeal to the theory of singular integral operators to obtain $L^p$ estimates for $U_\nu$. Indeed, the underlying space $\mathbb{R}^n$ has been endowed with a quasi-distance $d(x, y)$ and associated balls $B(x, r) = \{y : d(x, y) < r\}$. The measures of these balls satisfy the doubling property (5.10): $|B(x, 2r)| \leq C|B(x, r)|$. Consequently by, for example, the case $q = 2$ of Theorem 3 in Chapter 1 of [43], $U_\nu$ extends to a bounded operator on $L^p$, for every $1 < p \leq 2$, with $\| U_\nu \|_{L^p \to L^p} \leq A_{p,\varepsilon} 2^{\nu \varepsilon}$, for the same arbitrarily small exponent $\varepsilon$ as in (18.1). This is assertion (16.4).

In the same way (16.5) follows: the $L^2$ estimate was proved in Section 17, and the $L^p$ estimate for this operator is deduced as before since its kernel also satisfies the condition (18.1), with $\nu = 0$. Hence taking into account the remarks after (16.4), we see that $T$ has been shown to extend to a bounded operator on $L^p$, for every $1 < p \leq 2$. To pass to exponents $p \in (2, \infty)$, we recall that by Proposition 14.2, $T_j^* = T_j' + R_j$, where $\| R_j \|_{L^p \to L^p} \leq C 2^{-j}$. The operators $T_j'$ are essentially like the $T_j$, except that $\gamma_t$ is replaced by $\gamma_t^{-1}$. The same analysis as for $\{T_j\}$ therefore implies that $\sum_{j \geq 0} T_j'$ is bounded on $L^p$, for every $1 < p \leq 2$. Therefore $T^* = \sum_{j \geq 0} T_j^*$ is also bounded on $L^p$ for the same range of exponents, and the boundedness of $T$ for $2 \leq p < \infty$ then follows by duality. The proof of Theorem 11.1 is complete.

## 19. The maximal function

To prove the maximal inequality, Theorem 11.2, fix a nonnegative auxiliary function $\zeta \in C_0^\infty(\mathbb{R}^k)$, strictly positive where $|t| < a/2$, supported where $|t| < a$, and satisfying $\int \zeta = 1$. Define $A_j$ by

$$A_j(f)(x) = \psi(x) \, 2^{jk} \int f(\gamma_t(x)) \, \zeta(2^j t) \, dt$$

for $j = 0, 1, 2, \ldots$. Define $A_*(f)(x) = \sup_j |A_j(f)(x)|$.

It suffices to prove that $f \mapsto A_*(f)$ is bounded on $L^p(\mathbb{R}^n)$, for all $1 < p \leq 2$. This boundedness being obvious for $p = \infty$, the assertion for $1 < p \leq 2$ would at once carry over to $1 < p \leq \infty$, by a standard interpolation argument. This in turn would imply the corresponding result for the maximal operator $M$ defined in (11.2), by the obvious majorization $M(f) \leq CA_*(f) + C \int_{|t| \leq a'} f(\gamma_t(x)) \, dt$, which is valid for all nonnegative $f$.

With $S_j$ as in Section 16, let

(19.1) $$\triangle_j = A_j - S_j,$$

and define the square function $\mathcal{G}$ by

$$(19.2) \qquad \mathcal{G}(f)(x) = \left( \sum_{j=0}^{\infty} |\triangle_j(f)(x)|^2 \right)^{1/2}.$$

It will be enough to prove that

$$(19.3) \qquad \| \mathcal{G}(f) \|_{L^p} \leq A_p \| f \|_{L^p}, \quad \text{for all } 1 < p \leq 2,$$

because

$$\sup_j |A_j(f)| \leq \mathcal{G}(f) + \sup_j |S_j(f)|,$$

and $f \mapsto \sup_j |S_j(f)|$ is majorized by a constant multiple of the standard maximal function associated to the quasi-distance $d$ and to the associated balls. This standard maximal function is bounded on $L^p$ for all $p \in (1, \infty]$; see e.g. [43, Th. 1 in Ch. 1].

To deduce (19.3), we bring in the Rademacher functions $\{r_j(\omega)\}$, and write

$$(19.4) \qquad T^{(\omega)} = \sum_{j=0}^{\infty} r_j(\omega) \triangle_j.$$

Here $\omega$ represents a point in the underlying probability space that labels the arbitrary choices of plus or minus sign in the definition $r_j(\omega) = \pm 1$. As is well known (see e.g. [39, Ch. IV, §5]), the inequality (19.3) can be obtained from the inequality

$$(19.5) \qquad \| T^{(\omega)} \|_{L^p \to L^p} \leq A_p < \infty \quad \text{for all } 1 < p \leq 2$$

provided $A_p$ is independent of $\omega$.

The proof of (19.5) is very similar to that for the operator $T$ in Theorem 11.1. We give a brief outline, indicating only the changes needed. First, define $U_\nu^{(\omega)}$ by

$$U_\nu^{(\omega)} = \sum_{j=0}^{\infty} r_j(\omega) \triangle_j R_{j+\nu}.$$

We claim that

$$\| U_\nu^{(\omega)} \|_{L^2 \to L^2} \leq A 2^{-\varepsilon \nu}, \quad \text{for some } \varepsilon > 0,$$

with bounds independent of $\omega$. By the almost orthogonality arguments used above, this can be established by showing in analogy with (17.3)–(17.5) that uniformly in $i, j$,

$$(19.6) \qquad \| \triangle_j^* \triangle_i \| + \| \triangle_j \triangle_i^* \| \leq A 2^{-\varepsilon |i-j|},$$

$$(19.7) \qquad \| \triangle_j R_{j+\nu} \| \leq A 2^{-\varepsilon \nu}.$$

(19.7) is the consequence of two inequalities: $\| A_j R_{j+\nu} \| \leq A 2^{-\varepsilon\nu}$, and $\| S_j R_{j+\nu} \| \leq A 2^{-\varepsilon\nu}$. The first of these two is an easy variant of (17.5), in the proof of which the cancellation property $\int K_j(t)\,dt = 0$ of $T_j$ did not enter, and so $T_j$ could just as well have been replaced by $A_j$. The inequality $\| S_j R_{j+\nu} \| \leq A 2^{-\varepsilon\nu}$ is a variant of (17.4), with $i$ replaced by $j+\nu$; again the cancellation property of the first factor $R_j$ played no role in the proof of (17.4), so that the reasoning applies equally well with $R_j$ replaced by $S_j$.

We will discuss only the bound for the first term on the left in (19.6), and only in the case $i \geq j$; the case $i < j$ and the bound for the second term may be obtained in the same manner. The desired majorization is a consequence of two inequalities: $\| A_j^* \triangle_i \| \leq A 2^{-\varepsilon|i-j|}$ and $\| S_j \triangle_i \| \leq A 2^{-\varepsilon|i-j|}$. The first is reducible to $\| (A_j A_j^*)^N \triangle_i \| \leq A 2^{-\varepsilon'|i-j|}$. Assuming $i \geq j$, one estimates the kernel of the operator $(A_j A_j^*)^N \triangle_i$ to obtain a bound of $A 2^{-\varepsilon'|i-j|}$ for the $L^\infty$ operator norm, by the same line of reasoning as used above in the proof of (17.6). This leads to the consideration of the kernel

$$\int [h(\delta_{2^j}^x(z)) - h(\delta_{2^j}^x(y))]\,d\lambda_i^{(y)}(z)$$

for each $x$, where the measure $\lambda_i^{(y)}$ is defined by

$$\int f(z)\,d\lambda_i^{(y)}(z) = \psi(x) \int f(\gamma_t(y))\,\zeta(2^i t)\,dt - \int f(z)\,S_i(z,y)\,dz.$$

We write $i = j + \nu$, and define the nonnegative measure $\mu$ by

$$\int f(y,z)\,d\mu(y,z) = \int_{B(x,c2^{-j})} \left\{ \int f(y,z)\,d|\lambda_i^{(y)}(z)| \right\} dy.$$

Proposition 13.1 can now be applied to prove that the $L^\infty$ operator norm of $(A_j^* A_j)^N \triangle_i$ is $O(2^{-\delta'|i-j|})$. There is a corresponding result for $S_j \triangle_i$. It uses the same measure $\mu$, but in this case Proposition 13.2 suffices.

Given these $L^2$ estimates, the arguments establishing the $L^p$ inequalities, in particular $\| U_\nu^{(w)} \|_{L^p \to L^p} \leq A 2^{\varepsilon\nu}$, are almost unchanged from those carried out in Section 18. This is because the kernels of the operators $\triangle_j R_{j+\nu}$ satisfy the same estimates as those for the operators $T_j R_{j+\nu}$, given by Lemma 18.1. It follows that $\sum_{j\geq 0} r_j(\omega) \triangle_j R_{j+\nu}$ satisfies the estimate (18.1) given for $\sum_j T_j R_{j+\nu}$. In a similar way, the treatments of $\sum_j \triangle_j S_j$ and $\sum_j T_j S_j$ are completely parallel. This concludes our sketch of the proof of Theorem 11.2.

## 20. The smoothing property

The purpose of this section is to prove Theorems 8.11 and 8.12. Assume that $\gamma : \mathbb{R}^{n+k} \mapsto \mathbb{R}^n$ satisfies $(\mathcal{C})$ at $x_0 \in \mathbb{R}^n$. Throughout this section we work directly in $\mathbb{R}^n$, rather than in the freed situation of Proposition 6.1.

Consider an operator $T(f)(x) = \int f(\gamma(x,t))\,\psi(x,t)dt$, where $\psi$ is a $C^1$ function supported in a sufficiently small neighborhood of $(x_0, 0)$. Associated to $T$ are the iterates $U_\ell = (TT^*)^{2^\ell}$.

LEMMA 20.1. *Suppose that $\gamma$ satisfies $(\mathcal{C})$ at $x_0$. If $\ell$ is chosen to be sufficiently large, then*

$$U_\ell(f)(x) = \int \mathcal{K}(x,y)\,f(y)\,dy$$

*where $\mathcal{K}(x,y)$ is compactly supported, and $y \mapsto \mathcal{K}(x,y)$ is in $L^1_\delta$, for some $\delta > 0$, uniformly in $x$. The same holds with the roles of $x$ and $y$ reversed.*

The proof of this lemma is a simpler reprise[17] of the proof of Proposition 15.1. The principal ingredients are the curvature condition $(\mathcal{C}_J)$ and Proposition 7.2. The details are left to the reader.

To exploit the lemma we utilize the spaces $\mathcal{L}^p_\alpha$ occurring in [39, Ch. 6]; see also [4]. Here $\mathcal{L}^p_\alpha = \{f : J_{-\alpha}(f) \in L^p\}$, where $J_\alpha(\hat{f}) = (1 + 4\pi^2|\xi|^2)^{-\alpha/2}\hat{f}(\xi)$. Note $H_s = \mathcal{L}^2_s$, and $\mathcal{L}^p_0 = L^p$. A known inclusion result[18] asserts $L^1_\delta \subset \mathcal{L}^p_{\delta_1}$, for some $p > 1$ and $\delta_1 > 0$. Since $\mathcal{K}(x,y)$ has compact support and belongs to $L^1_\delta$ as a function of $x$, uniformly in $y$, it follows that $U_\ell$ maps $L^p$ to $\mathcal{L}^p_{\delta_1}$ boundedly, i.e. $J_{-\delta_1}U_\ell$ is bounded from $L^p$ to itself. By duality, since $U_\ell^* = U_\ell$, it follows that $U_\ell J_{-\delta_1}$ is bounded from $L^{p'}$ to itself, where $1/p + 1/p' = 1$. Thus the analytic family $s \to J_{s-\delta_1/2}U_\ell J_{-s-\delta_1/2}$ is bounded on $L^p$, when $\mathrm{Re}(s) = -\delta_1/2$, and on $L^{p'}$ when $\mathrm{Re}(s) = \delta_1/2$. Thus by interpolation, $J_{-\delta_1/2}U_\ell J_{-\delta_1/2}$ is bounded on $L^2$.

However $(J_{-\delta_1/2}U_{\ell-1})(J_{-\delta_1/2}U_{\ell-1})^* = J_{-\delta_1/2}U_\ell J_{-\delta_1/2}$, so that $J_{-\delta_1/2}U_{\ell-1}$ is bounded on $L^2$, and by duality so is $U_{\ell-1}J_{-\delta_1/2}$. Another interpolation argument then shows that $J_{-\delta_1/4}U_{\ell-1}J_{-\delta_1/4}$ is bounded on $L^2$. Continuing this way we conclude that $J_{-\delta_\ell}U_0 J_{-\delta_\ell}$ is bounded on $L^2$ with $\delta_\ell = \delta/2^{\ell+1}$. This means that $J_{-\delta_\ell}TT^*J_{-\delta_\ell}$ is bounded on $L^2$, and hence $J_{-\delta_\ell}T$ is bounded. This shows that $T$ maps $L^2$ to $\mathcal{L}^2_{\delta_\ell}$, i.e. maps $L^2$ to $H_s$, with $s = \delta_\ell$. One half of Theorem 8.11 is therefore proved.

We turn to the converse half. When $\gamma$ is not assumed to satisfy $(\mathcal{C})$, our hypotheses permit very degenerate mappings such as $\gamma(x,t) \equiv x$. In that case the operator $Tf(x) = \int f(\gamma(x,t))K(x,t)dt$ is identically zero, hence is smoothing, whenever $\int K(x,t)dt \equiv 0$. To exclude such pathology we say that $\gamma$ is at worst mildly degenerate at $x_0$ if there exists a neighborhood $V$ of $(x_0, 0)$ such that $\gamma(x,t) \neq x$ whenever $(x,t) \in V$ and $t \neq 0$. This holds if the differential of the map $t \mapsto \gamma(x_0, t)$ has maximal rank $k$ at $t = 0$.

---

[17]In the proof of Proposition 15.1, properties of the "free" context were used to show that all estimates were suitably uniform in the index $j$. Here this uniformity is no longer an issue.

[18]This elementary result may be proved using Littlewood-Paley technique.

We will prove a more general version of the converse half of the theorem: If $\gamma$ is at worst mildly degenerate and fails to satisfy curvature condition $(\mathcal{C})$ at $x_0$, then for any kernel $K \in C^\infty$ supported in a sufficiently small neighborhood of $(x_0, 0)$ satisfying $K(x_0, 0) \neq 0$, $T$ cannot map $L^2$ to $H_s$ for any $s > 0$. The proof will show that when $K \geq 0$, the same conclusion holds without the auxiliary hypothesis of at worst mild degeneracy. It relies on an elementary lemma, whose proof is omitted.

LEMMA 20.2. *For any $p \geq 1$, $s > 0$ and $\sigma \in (0,1)$ there exists $c > 0$ such that for every $r \in (0,1]$ and every real-valued $f \in C_0^\infty(\mathbb{R}^p)$ satisfying $f(x) \equiv 0$ for all $|x| \geq r$ and $f(x) \geq 1$ for all $|x| \leq \sigma r$,*
$$\|f\|_{H_s(\mathbb{R}^p)} \geq c r^{-s} r^{p/2}.$$

*Proof of converse half of Theorem* 8.11. Suppose that $(\mathcal{C})$ does not hold at $x_0$. Then according to Lemma 9.7, for some $0 \leq d < n$ there exist coordinates $(x', x'') \in \mathbb{R}^d \times \mathbb{R}^{n-d}$ with origin at $x_0$ such that $\gamma(x', x'', t) = (\gamma'(x,t), \gamma''(x,t))$ takes the form

(20.1) $$\gamma(x,t) \sim (x' + O(t), x'' + O(x'')O(t)).$$

Set $p = n - d \geq 1$.

Let $N$ be a large parameter, to be chosen at the conclusion of the proof. Fix a nonnegative function $\varphi \in C_0^\infty(\mathbb{R}^n)$ satisfying $\varphi(x) = 1$ whenever $|x| \leq 1/2$, and $\varphi(x) \equiv 0$ when $|x| \geq 1$. For small $\delta > 0$ set $f_\delta(x) = \varphi(\delta^{-1} x', \delta^{-N} x'')$. Then $\|f_\delta\|_{L^2} \leq C \delta^{(d+Np)/2}$.

If $\varepsilon > 0$ is chosen to be sufficiently small, then for all sufficiently small $\delta > 0$, whenever $|x'| \leq \varepsilon \delta$, $|x''| \leq \varepsilon \delta^N$ and $|t| \leq \varepsilon \delta$, we have $|\gamma'(x,t)| \leq \delta/4$ and $|\gamma''(x,t)| \leq \delta^N/4$, by (20.1). Moreover, if $V$ is chosen to be sufficiently small but independent of $N, \delta$, then whenever $2\delta^N \leq |x''| \leq 3\delta^N$, we have $|\gamma''(x,t)| \geq \delta^N$.

Assume without loss of generality that $K(x_0, 0) > 0$. The hypothesis that $\gamma$ is at worst mildly degenerate at $x_0$ ensures that if $V, \delta$ are sufficiently small and $x$ is sufficiently close to $x_0$ then $K(x,t) > 0$ whenever $\gamma(x,t)$ belongs to the support of $f_\delta$. Therefore by the preceding paragraph, there exist $\varepsilon, \varepsilon' > 0$ such that $Tf_\delta(x) \geq \varepsilon' \delta^k$ whenever $|x'| \leq \varepsilon \delta$ and $|x''| \leq \varepsilon \delta^N$. On the other hand, $Tf(x) \equiv 0$ if $x \in V$ and $2\delta^N \leq |x''| \leq 3\delta^N$.

Fix any $|x'| \leq \varepsilon \delta$ and consider the function $g_{x'}(x'') = Tf_\delta(x', x'')$. Applying Lemma 20.2 to the function $\delta^{-k} g_{x'}$, with $r = 2\delta^N$, yields
$$\|g_{x'}\|_{H_s(\mathbb{R}^p)} \geq c \delta^{k-Ns} \delta^{Np/2}.$$

But
$$\|Tf_\delta\|_{H_s(\mathbb{R}^n)}^2 \geq \int_{|x'| \leq \varepsilon \delta} \|g_{x'}\|_{H_s(\mathbb{R}^p)}^2 \, dx' \geq c \delta^{2(k-Ns)} \delta^{d+Np}.$$

Therefore the ratio $\|Tf_\delta\|_{H_s} / \|f_\delta\|_{L^2}$ is bounded below by $c \delta^{k-Ns}$.

Let $s > 0$ be given, choose $N > k/s$ and let $\delta$ tend to zero. Now the ratio tends to infinity. Thus $T$ does not map $L^2$ to $H_s$ boundedly.

The same reasoning applies if $K$ is assumed to be nonnegative and $K(x_0, 0) > 0$, without the hypothesis that $\gamma$ be at worst mildly degenerate. □

An analogous result can be formulated in the $L^p$ context.

*Definition* 20.3. An operator $S$ is $L^p$ improving if for some $p \in (1, \infty)$, there exists an exponent $q > p$ such that $S$ is bounded from $L^p$ to $L^q$.

Note that if $K$ is bounded and supported in a sufficiently small neighborhood of $\{t = 0\}$, our operators $T$ are bounded from $L^p$ to itself for all $p \in [1, \infty]$, by Lemma 14.1. Consequently if such an improvement holds for one $p$, then by interpolation it holds for all $p \in (1, \infty)$, with $q = q(p)$. Theorem 8.12 thus asserts essentially that if $K \in C^0$ and $K(x_0, 0) \neq 0$, then $T$ is $L^p$ improving if and only if $\gamma$ satisfies $(\mathcal{C})$ at $x_0$.

*Proof of Theorem* 8.12. $|K|$ may be majorized by a function in $C_0^1$. Assuming $(\mathcal{C})$, we see that Theorem 8.11 asserts that the operator defined by replacing $K$ by such a majorant maps $L^2$ to the Sobolev space $H_s$ for some $s > 0$; $H_s$ in turn embeds in $L^q$ for some $q > 2$. The converse half follows from the pointwise lower bound derived for $Tf_\delta$ in the proof of Theorem 8.11 together with a small computation, which is left to the reader. □

*Remark* 20.1. To describe the optimal degree of smoothing for operators of this type is no simple matter. A number of results are known under various additional hypotheses; see [7], [9], [14], [30], [31], [32], [35], [37], [38].

## 21. Complements and remarks

In this section we discuss some further results related to Theorems 11.1 and 11.2.

(1) Our first remark concerns the almost everywhere convergence implied by the maximal theorem. Let $\gamma(x, t) = \gamma_t(x)$ satisfy the curvature condition $(\mathcal{C})$ described in Section 8, at $x = x_0$. Consider any exponent $p > 1$, and let $c_k$ equal the volume of the unit ball in $\mathbb{R}^k$. Then in any neighborhood of $x_0$ in which $(\mathcal{C})$ holds,

$$\lim_{r \to 0} c_k^{-1} r^{-k} \int_{|t| \leq r} f(\gamma_t(x))\, dt = f(x) \quad \text{a.e.,} \tag{21.1}$$

for $x$ near $x_0$, for every $f \in L^p(\mathbb{R}^n)$. Now (21.1) follows from Theorem 11.2 in the standard fashion.

The conclusion (21.1) also holds without the curvature hypothesis on $\gamma_t$, if we require instead that $(x,t) \mapsto \gamma_t(x)$ be real analytic in $(x,t)$. In fact, when $\gamma_t$ is real analytic, then the set of all points in a small connected neighborhood of $x_0$ at which $(\mathcal{C})$ holds is the complement of an analytic variety, hence either is an open set of full measure, or is empty. In the first case, the conclusion (21.1) follows for almost every $x$, because the assertion is purely local. In the second case, when the curvature condition fails identically, Lemma 9.8 asserts that outside a lower-dimensional variety, a neighborhood of $x_0$ can be fibered by (proper) submanifolds which are invariant under $\gamma_t$. An application of Fubini's theorem reduces the conclusion to the assertion (21.1) for each of those submanifolds. The result is then a consequence of the validity of (21.1) under hypothesis $(\mathcal{C})$, applied in a lower dimensional situation.

Bourgain [3] has shown that for arbitrary real analytic families of straight lines $\gamma(x,t) = x + tv(x)$ in the plane, the maximal function $M$ is bounded on $L^p$ for $p > 1$.

(2) The singular integrals in Theorem 11.1 may be generalized as follows. Suppose $K(x,\cdot)$ is for each $x \in \mathbb{R}^n$ a distribution on $\mathbb{R}^k$ which agrees with a function $K(x,t)$, when $0 \neq t \in \mathbb{R}^k$. Assume that $K(x,t)$ is supported where $x$ is near $x_0$ and $|t| \leq a$, and is $C^1$ in the open set $t \neq 0$. Suppose moreover that $K(x,t)$ satisfies the differential inequalities

$$(21.2) \qquad \sup_x |\partial_x^\alpha \partial_t^\beta K(x,t)| \leq A|t|^{-k-|\beta|}$$

when $0 \leq |\alpha| \leq 1$ and $0 \leq |\beta| \leq 1$. Finally, we require that for each $x$, the Fourier transform of the distribution $K(x,\cdot)$ is an $L^\infty$ function; we denote it by $\hat{K}(x,\xi)$. More particularly, we require that

$$(21.3) \qquad \sup_{x,\xi} |\partial_x^\alpha \hat{K}(x,\xi)| \leq A < \infty, \text{ for } 0 \leq |\alpha| \leq 1.$$

Now for each $f \in C_0^\infty(\mathbb{R}^n)$ we define $T(f)$ by

$$(21.4) \qquad T(f)(x) = (K(x,\cdot), F)$$

where $F$ is the test function given by $F(t) = f(\gamma_t(x))$.

COROLLARY 21.1. *Assume that $\gamma_t$ satisfies the curvature conditions in Section* 8. *The operator $T$ defined by* (21.4) *extends to a bounded operator on $L^p(\mathbb{R}^n)$, for all $1 < p < \infty$.*

The proof proceeds as follows. First (21.2) and (21.3) imply the cancellation condition

$$(21.5) \qquad \sup_x |(\partial_x^\alpha K(x,\cdot), \varphi_R)| \leq A, \quad 0 \leq |\alpha| \leq 1$$

for all normalized bump-functions $\varphi$, with $\varphi_R(t) = \varphi(t/R)$, with $A$ a constant independent of $R$. For this see [43, Prop. 3, §4.5, Ch. 6]. Next, using (21.5) and proceeding as in [43, Ch. 3, §5.33], we can decompose $K(x, \cdot)$ as

$$(21.6) \qquad K(x, \cdot) = c(x)\delta_0 + \sum_{j=0}^{\infty} K_j(x, t)$$

where $\delta_0$ is the delta function at the origin, in the $t$-variable, and $c$ is a $C^1$ function of compact support. Here the $K_j(x,t)$ are supported in $2^{-j} \cdot a/2 \leq |t| \leq 2^{-j} \cdot 2a$, and satisfy (21.3) uniformly in $j$. Moreover they also satisfy the cancellation conditions

$$(21.7) \qquad \int K_j(x,t)dt = 0, \quad \text{for all } j \text{ and } x.$$

With the decomposition (21.6) we can define $T_j$ by $T_j(f)(x) = \int K_j(x,t)f(\gamma_t(x))dt$, and we have $T(f)(x) = c(x)f(x) + \sum_{j=0}^{\infty} T_j(f)(x)$, for every $f \in C_0^1$. The cancellation conditions (21.7) suffice for the proof of the $L^p$ boundedness of $\sum_j T_j$ to proceed as before.

# Part 4. Appendix

## 22. Proof of the lifting theorem

To prove Proposition 6.1, we need to find a manifold $M$ of dimension $d$ (with a distinguished point $O \in M$), and a mapping $\pi$ of a neighborhood of $O$ to a neighborhood of the origin in $\mathbb{R}^n$, so that $\pi$ has maximal rank ($= n$); in addition we need to construct vector fields $\tilde{X}_i$, $1 \leq i \leq p$, defined on $M$ near $O$, so that $\pi_*(\tilde{X}_i) = X_i$, and $\tilde{X}_i$ are free relative to the exponents $a_1, \ldots, a_p$ up to order $m$.

In fact, by the maximal rank property there would be a coordinate system $\{(x,z)\}$ for $M$ centered at $O$, with $x \in \mathbb{R}^n$, $z \in \mathbb{R}^{n-d}$ so that $\pi(x,z) = x$. Then $\pi_*(\tilde{X}_i) = X_i$ for all $i$ implies that the $\tilde{X}_i$ take the form (6.3).

Recall the Lie algebra $\mathcal{N} = \mathcal{N}_m^{a_1,\ldots a_p}$ and its corresponding Lie group $N$; also the basis $\{Y_I\}$ of $\mathcal{N}$, where each $Y_I$ is an appropriate commutator of the elements $Y_1, Y_2, \ldots, Y_p$. The $Y_I$'s are (left-invariant) vector fields on $N$.

Consider now the space $N \times \mathbb{R}^n$, and the vector fields $U_j$ defined near the origin of this space by

$$(22.1) \qquad U_j = Y_j \oplus X_j,$$

(i.e. where $Y_j$ acts on the $N$-variable, and $X_j$ acts on the $\mathbb{R}^n$ variable).

We can also form the commutators $U_I = Y_I \oplus X_I$ for $I$ basic, and note that the $U_I$'s are linearly independent (because the $Y_I$ are linearly independent).

574    CHRIST, NAGEL, STEIN, AND WAINGERWe next define the submanifold $M$ of $N \times \mathbb{R}^n$ by a (local) diffeomorphism $\Phi$ with $N$

$$\Phi : N \to M$$

given by

(22.2) $$\Phi : \exp\left(\sum_{I \text{ basic}} y_I Y_I\right) \to \exp\left(\sum_{I \text{ basic}} \exp y_I U_I\right)(O),$$

which is defined whenever $(y_I)$ is small.

Let $\tilde{\pi}$ denote the coordinate projection of $N \times \mathbb{R}^n$. We then define the (local) mapping $\pi$ of $M$ to $\mathbb{R}^n$ by restricting $\tilde{\pi}$ to $M$.

It follows from (22.2) that $U_I(O)$ are tangent vectors to $M$ at the origin. Since $\pi_*(U_I(O)) = \tilde{\pi}_*(U_I(O)) = X_I(O)$, and the $X_I$ span, we see that $\pi$ has rank $= n$ near the origin.

Having defined $M$ and $\pi$, we need next to define the vector fields $\tilde{X}_i$. Let $V_i = \Phi_*(Y_i)$, $V_I = \Phi_*(Y_I)$. If the $V_i$'s equaled the $U_i$'s, the latter would be tangent to $M$ and we could take these to be the $\tilde{X}_i$. However in general this identity does not hold. As a substitute we prove

(22.3) $$V_i - U_i \quad \text{vanishes of order } > m - a_i.$$

To make this statement precise, we temporarily introduce a new coordinate system near the origin in $N \times \mathbb{R}^n$, adapted to the vector fields $\{U_I\}$. We choose $n$ vector fields, $Z_1, Z_2, \cdots Z_n$, so that $\{U_I, I \text{ basic, and } Z_j, 1 \leq j \leq n\}$ spans $N \times \mathbb{R}^n$ near the origin. With this choice fixed, we assign to the point $\exp\left(\sum_{I \text{ basic}} y_I U_I + \sum_j v_j Z_j\right)(O)$ the coordinate $[y, v]$ where $y = (y_I)$, $v = (v_j)$. Note that in this coordinate system, the submanifold $M$ is given by $v = 0$, while $(y_I)$ can be taken as a coordinate system for $M$ near $O$.

Now since $\Phi$ is a local diffeomorphism, the commutation relations satisfied by the $Y_I$ (i.e. those of the Lie algebra $\mathcal{N}$) are the same as those satisfied by the $V_I = \Phi_*(Y_I)$. Moreover, in view of $\Phi$, in the coordinates $y$, the vector field $V_i$ is expressed by

(22.4) $$\frac{d}{dt} f\left(\exp tY_i \cdot \exp\left(\sum_{I \text{ basic}} y_I Y_I\right)\right)\bigg|_{t=0},$$

while the expression of the $U_i$ in terms of the coordinates $[y, v]$ is given by

(22.5) $$\frac{d}{dt} f\left(\exp tU_i \cdot \exp\left(\sum_{I \text{ basic}} y_I U_I + \sum_{j=1}^n v_j Z_j\right)(O)\right)\bigg|_{t=0}$$

so we need to compare

$$\exp t Y_i \cdot \exp \sum_I y_I Y_I \quad \text{and} \quad \exp t U_i \cdot \exp\left(\sum_I y_I U_I + \sum_j v_j Z_j\right)$$

via the Baker-Campbell-Hausdorff formula. To make the comparison we define an order on monomials of $y = (y_I)$ (consistent with that of the vector fields $Y_I$), by assigning $y_j$ to have order $a_j$, and making the order multiplicative in the usual way. Then it is clear that in multiplying out $\exp tY_i \cdot \exp\left(\sum_I y_I Y_I\right)$, the coefficient of any $Y_{I_0}$ is a polynomial function of $t$ whose coefficients are polynomials in $y$, in which the coefficient of the linear term $t^1$ has order $|I_0| - a_i$. On $M$ we have $v = 0$, and so using the reasoning which follows (5.6), we can see that as long as we consider commutators of $U_i$ with various $U_I$ whose total order is $\leq m$, the effect is the same as with $Y_i$ and $Y_I$; thus only commutators which involve $U_i$ and the $U_I$ of total order $> m$ lead to a different answer. In calculating (22.5) we need only consider the coefficient of $t^1$, and from this we get (22.3).

Next let us denote by $T_\mu(M)$ the tangent space of $M$ at $\mu \in M$; we adopt a similar terminology for $T_\mu(N \times \mathbb{R}^n)$, and we note that $T_\mu(M) \subset T_\mu(N \times \mathbb{R}^n)$. We can then find a subspace $S$ of $T_O(N \times \mathbb{R}^n)$ so that

$$(22.6) \qquad T_O(M) + S = T_O(N \times \mathbb{R}^n), \text{ and } S \subset \text{ kernel } (\tilde{\pi})_*.$$

Indeed, $T_O(M)$ is spanned by the $U_I(O) = Y_I(0) + X_I(0)$ and the kernel of $\tilde{\pi}_*$ is $\mathcal{N}$, which is spanned by $Y_I(0)$. Since the $X_I(0)$ span $\mathbb{R}^n$, such an $S$ can be determined. With this $S$ fixed we also have

$$(22.7) \qquad T_\mu(M) + S = T_\mu(N \times \mathbb{R}^n) = T_O(N \times \mathbb{R}^n)$$

for $\mu$ near the origin. The decomposition (22.7) allows us to define a smoothly varying family of projections, $P_\mu$ of $T_\mu(N \times \mathbb{R}^n)$, onto its subspace $T_\mu(M)$, by setting $P_\mu$ to be the identity on $T_\mu(M)$, and $P_\mu$ to vanish on $S$. We are then in a position to define the $\tilde{X}_i$ by

$$(22.8) \qquad \tilde{X}_i(\mu) = P_\mu(U_i(\mu)).$$

Now first,

$$\pi_*(\tilde{X}_i) = \tilde{\pi}_*(\tilde{X}_i) = \tilde{\pi}_*(U_i) - \tilde{\pi}_*((1 - P_\mu)U_i(\mu)) = \tilde{\pi}_*(U_i) = X_i.$$

Therefore, $\pi_*(\tilde{X}_i) = X_i$. Also applying $P_\mu$ to (22.3) we see that since $P_\mu(V_i)(\mu) \equiv V_i(\mu)$ for $\mu \in M$, and because $P_\mu$ depends smoothly on $\mu$, $V_i - \tilde{X}_i$ vanishes of order $> m - a_i$. From this it is clear by induction on $r$ that the span $\{\tilde{X}_I(O), |I| \leq r\}$ = the span $\{V_I(O), |I| \leq r\}$, for $r \leq m$. Since the Lie algebra generated by the $V_I$'s is isomorphic to that of the $Y_I$'s (i.e. $\mathcal{N}$), we see as a result that the $\tilde{X}_i$ are free relative to $a_1, \ldots a_p$ up to order $m$, and the proposition is proved.

University of California, Berkeley, CA
*E-mail address*: mchrist@math.berkeley.edu

University of Wisconsin, Madison, WI
*E-mail address*: nagel@math.wisc.edu

Princeton University, Princeton, NJ
*E-mail address*: stein@math.princeton.edu

University of Wisconsin, Madison, WI
*E-mail address*: wainger@math.wisc.edu